\documentclass[reqno,oneside,a4paper]{amsart}
\usepackage[a4paper, left=2.5cm, right=2.5cm, top=2.5cm, bottom=2.5cm]{geometry}
\usepackage{amsmath,amssymb,latexsym}
\usepackage{amsthm}
\usepackage{amsfonts}
\usepackage{rotating}
\usepackage{graphicx,tikz,float,subfigure}
\usepackage{xcolor}
\usepackage{latexsym}
\usepackage{extarrows}
\usepackage{arydshln}
\usepackage[colorlinks=true,linkcolor=blue,citecolor=red,anchorcolor=green]{hyperref}
\allowdisplaybreaks
\theoremstyle{definition}
\newtheorem{Theorem}{Theorem}[section]
\newtheorem{lemma}[Theorem]{Lemma}
\newtheorem{Proposition}[Theorem]{Proposition}
\newtheorem{Corollary}[Theorem]{Corollary}
\newtheorem{remark}{Remark}
\newtheorem{claim}{Claim}
\newtheorem{Definition}{Definition}
\numberwithin{equation}{section}

\newcommand{\la}{\langle}
\newcommand{\ra}{\rangle}

\newcommand{\beq}{\begin{equation}}
	\newcommand{\eeq}{\end{equation}}
\newcommand{\bes}{\begin{equation*}}
	\newcommand{\ees}{\end{equation*}}

\newcommand{\Id}{{\rm Id}}

\newcommand{\A}{{\mathbb A}}
\newcommand{\B}{{\mathbb B}}
\newcommand{\Z}{{\mathbb Z}}
\newcommand{\R}{{\mathbb R}}

\newcommand{\C}{{\mathbb C}}

\newcommand{\T}{{\mathbb T}}
\newcommand{\J}{{\mathbb J}}

\newcommand{\N}{{\mathbb N}}

\newcommand{\CA}{{\mathcal A}}
\newcommand{\CC}{{\mathcal C}}

\newcommand{\CE}{{\mathcal E}}
\newcommand{\CF}{{\mathcal F}}
\newcommand{\CG}{{\mathcal G}}
\newcommand{\CB}{{\mathcal B}}
\newcommand{\CH}{{\mathcal H}}
\newcommand{\CI}{{\mathcal I}}
\newcommand{\CQ}{{\mathcal Q}}
\newcommand{\CL}{{\mathcal L}}
\newcommand{\CM}{{\mathcal M}}

\newcommand{\CX}{{\mathcal X}}

\newcommand{\CS}{{\mathcal S}}

\newcommand{\CU}{{\mathcal U}}
\newcommand{\CT}{{\mathcal T}}
\newcommand{\CK}{{\mathcal K}}

\newcommand{\CJ}{{\mathcal J}}
\newcommand{\CR}{{\mathcal R}}

\def\ad{{\rm ad}}

\def\Id{{\rm Id}}

\def\beq{\begin{equation}}
	\def\eeq{\end{equation}}

\def\det{\mathrm{det}\ }
\def\Ker{{\mathrm{Ker}}}

\def\N{{\Bbb N}}
\def\Z{{\Bbb Z}}
\def\R{{\Bbb R}}
\def\T{{\Bbb T}}
\def\C{{\Bbb C}}

\let\cal=\mathcal
\def\H{{\cal H}}
\hfuzz=\maxdimen
\def\N{{\Bbb N}}
\def\Z{{\Bbb Z}}
\def\R{{\Bbb R}}
\def\T{{\Bbb T}}
\def\C{{\Bbb C}}
\def\M{{\cal M}}

\def\Sp{{\rm Sp}}
\def\sp{{\rm sp}}

\let\cal=\mathcal
\hfuzz=\maxdimen
\begin{document}
\title[]{Symplectic Normal Form and Growth of Sobolev Norm}

\author{Zhenguo Liang}
\address{School of Mathematical Sciences and Key Lab of Mathematics for Nonlinear Science, Fudan University, Shanghai 200433, China}
\email{zgliang@fudan.edu.cn}
\thanks{Z. Liang was partially supported by NSFC grant 12071083 and partially supported by the New Cornerstone Science Foundation through the New Cornerstone Investigator Program. }

\author{Jiawen Luo}
\address{School of Mathematical Sciences and Key Lab of Mathematics for Nonlinear Science, Fudan University, Shanghai 200433, China}
\email{20110180010@fudan.edu.cn}
	
\author{Zhiyan Zhao}
\address{Universit\'e C\^ote d'Azur, CNRS, Laboratoire J. A. Dieudonn\'{e}, 06108 Nice, France}
\email{zhiyan.zhao@univ-cotedazur.fr}
\thanks{Z. Zhao was partially supported by the French government through the National Research Angency (ANR) grant for the project KEN ANR-22-CE40-0016 and partially supported by NSFC grant 12271091.}

\begin{abstract}
For a class of reducible Hamiltonian partial differential equations (PDEs) with arbitrary spatial dimension, quantified by a quadratic polynomial with time-dependent coefficients, we present a comprehensive classification of long-term solution behaviors within Sobolev space. This classification is achieved through the utilization of Metaplectic and Schr\"odinger representations. Each pattern of Sobolev norm behavior corresponds to a specific $n-$dimensional symplectic normal form, as detailed in Theorems \ref{main1} and \ref{main2}.

When applied to periodically or quasi-periodically forced $n-$dimensional quantum harmonic oscillators, we identify novel growth rates for the $\mathcal{H}^s-$norm as $t$ tends to infinity, such as $t^{(n-1)s}e^{\lambda st}$ (with $\lambda>0$) and $t^{(2n-1)s}+ \iota t^{2ns}$ (with $\iota\geq 0$). Notably, we demonstrate that stability in Sobolev space, defined as the boundedness of the Sobolev norm, is essentially a unique characteristic of one-dimensional scenarios, as outlined in Theorem \ref{thm_inverse}.

As a byproduct, we discover that the growth rate of the Sobolev norm for the quantum Hamiltonian can be directly described by that of the solution to the classical Hamiltonian which exhibits the optimal growth, as articulated in Theorem \ref{thm-ODEPDE}.

\
		

\
		
		
\end{abstract}

	\maketitle
	
	\section{Introduction of the Main Results}

This paper focuses on the Sobolev norms of solutions to Hamiltonian partial differential equations (PDEs), particularly aiming to introduce new behavioral patterns and to provide a thorough classification of growth rates for certain quantum Hamiltonians. This task is approached through a reducibility argument and a complete classification of quantum normal forms.

In view of the significance of the quadratic quantum Hamiltonian in illustrating wave propagation phenomena and the connection with the coherent state, the study concentrates on Hamiltonian PDEs quantized by a quadratic polynomial with time dependent coefficients. Specifically, the work delves into the $n-$dimensional PDE
	\begin{equation}\label{orig-equ}
		\frac{1}{2\pi {\rm i}}\partial_t \psi(t,x)= H_0(t,D, X)\psi(t,x), \qquad x=(x_1,\cdots, x_n)\in\R^n,
	\end{equation}
	where
	$\displaystyle D	=\begin{pmatrix}
		D_1\\ \vdots \\ D_n
	\end{pmatrix}:= \frac{1}{2\pi {\rm i}}\nabla$,
	$\displaystyle  X:=\begin{pmatrix}
		X_1\\  \vdots \\ X_n
	\end{pmatrix}$ with $X_j$ the multiplication by $j-$th coordinate function, i.e., $(X_jf)(x)=x_jf(x)$, and $H_0(t,D,X)$ is a time dependent quadratic polynomial in $(D,X)$, i.e.,
	\begin{eqnarray*}
		H_0(t,D,X)&=&\frac12 \left\la D,A_{12}(t)D \right\ra -\frac12  \left\la D,A_{11}(t) X \right\ra-\frac12  \left\la X, A_{11}(t)^*D \right\ra-\frac12 \left\la X,A_{21}(t) X \right\ra \\
		& &+ \,  \left\la l_1(t),D \right\ra+ \left\la l_2(t),X \right\ra+c(t),
	\end{eqnarray*}
	with the superscript ``$*$'' the transpose of matrix or vector, and $c(\cdot) \in C_b^{0}(\R,\R)$, $l_j(\cdot) \in C_b^{0}(\R,\R^n)$, $A_{ij}(\cdot)\in C_b^{0}(\R,{\rm gl}(n,\R))$ satisfying $A_{12}(t)=A_{12}(t)^*$, $A_{21}(t)=A_{21}(t)^*$. Here, ``$C_b^{0}$" means the continuous class with uniform boundedness w.r.t. $t\in \R$ (see Section \ref{sec_notations}--(\ref{Cb01})).

Define the column of operators $Z$ by $Z^*=(D^*, X^*)=(D_1,\cdots, D_n,X_1,\cdots, X_n)$. Then the quadratic quantum Hamiltonian (\ref{orig-equ}) is rewritten as
\begin{equation}\label{orig-equ-1}
\frac{1}{2\pi {\rm i}}\partial_t \psi=(\CQ_{\CA(t)}(Z)+\CL_{\ell(t)}(Z)+c(t))\psi,
\end{equation}
with the linear operator $H_0(t,D,X)$ in Eq. (\ref{orig-equ}) decomposed into
	\begin{itemize}
		\item {\it Homogeneous part} $\CQ_{\CA(t)}(Z):=-\frac12\la Z,\mathcal{A}(t)\J_nZ\ra$ with
		$$\CA(\cdot)= \left(\begin{array}{cc}
			A_{11}(\cdot) & A_{12}(\cdot) \\ A_{21}(\cdot) & -A_{11}(\cdot)^* \end{array}\right)
		\in C_b^{0}(\R, \sp(n,\R)),
		\footnote{In this context, $\sp(n,\R)$ represents the set of $2n\times 2n$ Hamiltonian matrices with real entries. Refer to Section \ref{Section_MetaRepre} for precise definitions.}
		\qquad \J_n:=\left(\begin{array}{cc}
			0 & I_n \\ -I_n & 0 \end{array}\right), $$
		\item {\it Linear part} $\CL_{\ell(t)}(Z):=\la \ell(t),Z\ra$ with $\ell(\cdot)=(l_1(\cdot)^*,l_2(\cdot)^*)^* \in C_b^{0}(\R,\R^{2n})$,
		\item {\it Scalar part} $c(t)$.
	\end{itemize}
	It is quantized by the quadratic polynomial classical hamiltonian (usually called the {\it symbol} of $H_0$)
	\begin{equation}\label{classicalHam}
		h_0(t,\xi,x) =-\frac12\la z,\mathcal{A}(t)\J_n z\ra+\la \ell(t), z\ra + c(t),\qquad \xi,x\in \R^{n},
	\end{equation}
	where $z\in \R^{2n}$ such that $z^*=(\xi^*,x^*)$, and the three terms of $h_0(t,\xi,x)$ are also called {\it homogeneous part}, {\it linear part} and {\it scalar part} respectively in the classical hamiltonian sense.

	As a typical example of the quantum Hamiltonian $H_0$ and a  well-known ``equilibrium state", the $n-$D quantum harmonic oscillator (QHO for short through the paper),
	\begin{equation}\label{nD-QHO}
		\CT=\sum_{j=1}^n\CT_j:=2\pi\left(\la D,D\ra+ \la X, X\ra\right)=-\frac{1}{2\pi} \Delta+2\pi |X|^2, \qquad \CT_j:=2\pi (D_j^2+X_j^2),
			\end{equation}
	i.e., $H_0$ with constant coefficients $\mathcal{A}(\cdot)=4\pi \J_n$, $\ell(\cdot)=0$, $c(\cdot)=0$,	as well as their perturbations, are well investigated in many recent works and it is much related to the present paper.

\subsection{Reducibility and symplectic normal forms}
	
For a specific quantum Hamiltonian $H_0(t,D,X)$ with constant coefficients, a straightforward approach exists for analyzing the propagators of Eq. (\ref{orig-equ-1}). This includes, for example, Mehler's formula for the $n-$D QHO $\CT$ defined in (\ref{nD-QHO}). The properties of the solutions are explored through direct computations. Consequently, to efficiently describe the behavior of solutions to the time-dependent equation, it is practical to remove the time-dependence from the original equation via an appropriate transformation. This process is referred to as {\it reducibility} of the quantum Hamiltonian in this paper.

A considerable amount of prior research has focused on the reducibility of time-dependent PDEs. For $1-$D QHOs, time periodic smooth perturbations were initially examined \cite{Com87, DLSV2002, EV83, HLS86, Kuk1993}. In the context of time quasi-periodic perturbations, both bounded \cite{GreTho11, LiangWZQ2022, Wang08, WLiang17} and unbounded disturbances \cite{Bam2018, Bam2017, BM2018, Liangluo19, LZZ2021, LLZ2022} have been thoroughly investigated. The study of reducibility issues has led to the comprehensive development of infinite-dimensional KAM theory for $1-$D PDEs with unbounded perturbations, as seen in works by Bambusi-Graffi \cite{BG2001}, Kuksin \cite{Kuk1997}, and Liu-Yuan \cite{LY2010}.

Regarding higher-dimensional PDEs, the concept of reducibility was pioneered by Eliasson-Kuksin \cite{EK2009, EK2010} for the quasi-periodic Schr\"odinger equation. For higher-dimensional QHOs, references \cite{GrePat16, LiangWang19, LiangWZQ2022plus} discuss bounded potentials, and \cite{BGMR2018} examines unbounded perturbations, including Hamiltonian PDEs quantized by quadratic polynomials. Through the lens of reducibility, these results offer insights into the long-term behavior of solutions (at least in terms of upper bounds of Sobolev norms) and the spectral characteristics of the Floquet operator.
Recently, the reducibility of a class of Landau Hamiltonian has been studied by Bambusi-Gr\'ebert-Maspero-Robert-Villegas Blas \cite{BGMRV2025}. Through new KAM algorithms for the classical dynamics, it is shown that the Landau Hamiltonian is reducible
to a two-dimensional harmonic oscillator and thus gives rise to bounded dynamics.

The concept of reducibility for higher-dimensional PDEs was pioneered by Eliasson-Kuksin \cite{EK2009, EK2010} in the context of the quasi-periodic Schr\"odinger equation. 
For higher-dimensional QHOs, the works \cite{GrePat16, LiangWang19, LiangWZQ2022plus} address cases with bounded potentials, while \cite{BGMR2018} extends the analysis to unbounded perturbations, including Hamiltonian PDEs quantized by quadratic polynomials. Viewed through the framework of reducibility, these studies provide insight into the long-term behavior of solutions—particularly upper bounds for Sobolev norms—and the spectral properties of the associated Floquet operator.

More recently, Bambusi-Gr\'ebert-Maspero-Robert-Villegas Blas \cite{BGMRV2025} have investigated the reducibility of a class of Landau Hamiltonians. Employing novel KAM algorithms tailored to the classical dynamics, they demonstrate that the Landau Hamiltonian is reducible to a two-dimensional harmonic oscillator, thereby establishing boundedness of the dynamics.

As deduced from a quantization argument (see Lemma \ref{reducibleprop} and \ref{reduciblepropnew}), the reducibility of the quadratic quantum Hamiltonian (\ref{orig-equ}) is closely linked to the reducibility of the corresponding finite-dimensional system. To fully comprehend the dynamics of solutions to Eq. (\ref{orig-equ-1}), it is crucial to understand the precise forms of the reduced constant equations, termed here as the {\it quantum normal forms} (QNFs). In essence, a comprehensive classification of quantum normal forms provides a complete depiction of the long-term behavior of solutions.

\medskip
	
Prior to presenting the main theorem regarding QNFs, it is essential to define reducibility for both finite-dimensional and infinite-dimensional systems precisely.
\begin{Definition}\label{def-reduc}{\bf (Reducibility in classical and quantum Hamiltonians)}
The affine system
$$z'=\CA(t)z+\ell(t),\qquad \CA(\cdot)\in C_b^0(\R, \sp(n,\R)),\qquad \ell(\cdot)\in C_b^0(\R, \R^{2n})  $$
 is called {\it reducible}, if there exist $\CB\in \sp(n,\R)$, $l\in \R^{2n}$ and $S(\cdot)\in C_b^1(\R,\Sp(n,\R))$, $p(\cdot)\in C_b^1(\R,\R^{2n})$
 \footnote{$C_b^1(\R,\bullet)$, with $\bullet=\Sp(n,\R)$ or $\R^{2n}$, is the subspace of $C^0_b(\R, \bullet)$ whose entries are all $C^1$ with uniformly bounded derivatives on $\R$. See Section \ref{sec_notations}--(\ref{Cb01}).}, such that it is conjugated to $y'=\CB y+l$ under the change of variable $z=S(t) y+p(t)$. The quadratic quantum Hamiltonian with time dependent coefficients
$\frac{1}{2\pi {\rm i}}\partial_t \psi= H(t,Z)\psi $
is called {\it reducible}, if there exist an $L^2-$unitary transformation $\psi(t)=\CU(t)\varphi(t)$, and a quadratic quantum Hamiltonian $\CR(Z)$ with constant coefficients, such that it is
conjugated to $\frac{1}{2\pi {\rm i}}\partial_t \varphi= \CR(Z)\varphi$.
\end{Definition}

\begin{remark}
It is important to note that in the definition of reducibility, the regularity required for the time-dependence in the variable transformation can differ based on the regularity of the original system. For further details on this variation, refer to Section \ref{secAppli-qp}.
\end{remark}

\begin{Theorem}\label{main1} {\bf (Reducibility and classification of QNFs)}
If the system $z'=\CA(t)z+\ell(t)$, with $\CA(\cdot)\in C_b^0(\R,\sp(n,\R))$ and $\ell(\cdot)\in C_b^0(\R, \R^{2n})$ as in Eq. (\ref{orig-equ-1}),  is reducible, then Eq. (\ref{orig-equ-1}) is also reducible, via some $L^2-$unitary transformation $\psi(t)=\CU(t)\varphi(t)$, to
$\frac{1}{2\pi {\rm i}}\partial_t \varphi(t,x)=\CR(Z)\varphi(t,x)$,
where the operator $\CR(Z)$ is a time-independent operator, and in one of following forms.
\begin{itemize}
\item [\bf QNF1]  For $n\in\N^*$, $\CR(Z)=\CR^{(n)}_1(Z):=\CQ_{\CB^{(n)}_1(\lambda)}=-\frac12\left\la Z,\CB^{(n)}_1(\lambda)\J_n Z\right\ra$, $\lambda>0$, with \begin{equation}\label{A1}\CB^{(n)}_1(\lambda):=\left(\begin{array}{cccc:cccc}
					-\lambda & -1  &    &         &  &  &  &     \\
					& \ddots  &  \ddots  &         &  &  &  &    \\
					&   &  -\lambda  &  -1       &  &  &  &    \\
					&         &   & -\lambda &  &  &  &   \\
					\hdashline
					&  &   &  & \lambda &      &        &    \\
					&  &   &  &      1    & \lambda &  &    \\
					&  &   &  &          & \ddots  & \ddots &  \\
					&  &   &  &          &   &  1  &  \lambda
				\end{array}
				\right)\in\sp(n,\R). \footnote{Throughout this paper, the $2n\times 2n$ Hamiltonian and symplectic matrices (refer to Section \ref{Section_MetaRepre} for details) are frequently expressed in the form of $n\times n$ block matrices. To streamline the presentation, we only include entries that may be non-zero, thereby simplifying the expressions.}
			\end{equation}
In particular, $\CB^{(1)}_1(\lambda)=\begin{pmatrix}
						-\lambda & 0 \\
						0 & \lambda
					\end{pmatrix}$.			
\item [\bf QNF2] For $n\in\N^*$ even, $\CR(Z)=\CR^{(n)}_2(Z)=\CQ_{\CB^{(n)}_2(\lambda_1, \lambda_2)}$,
$\lambda_1,\lambda_2>0$, with
			\begin{equation}\label{A2}
				\CB^{(n)}_2(\lambda_1, \lambda_2):=\left(\begin{array}{cccc:cccc}
					A_2 &    & & & &  &  & \\
					-I_2  & A_2 & & & & & & \\
					& \ddots & \ddots & & & & & \\
					&  & -I_2 &  A_2 & & & & \\
					\hdashline
					& & & & -A_2^* & I_2 & & \\
					& & & & & \ddots & \ddots & \\
					& & & & & & -A_2^* & I_2 \\
					& & & & & & & -A_2^*
				\end{array}
				\right)\in\sp(n,\R),
			\end{equation}
			where $A_2:=\begin{pmatrix}
				-\lambda_1 & -\lambda_2 \\
				\lambda_2 & -\lambda_1
			\end{pmatrix}$. In particular, $\CB^{(2)}_2(\lambda_1, \lambda_2)=\begin{pmatrix}
						A_2 & 0 \\
						0 & -A_2^*
					\end{pmatrix}$.

\item [\bf QNF3] For $n\in\N^*$, for $\CR(Z)=\CR^{(n)}_{3}(Z)=\CQ_{\CB^{(n)}_3(\mu,\gamma)}$,
$\gamma=\pm 1$, $\mu>0$, where
			\begin{equation}\label{A3}
				\CB^{(n)}_3(\mu,\gamma):=\gamma\left(\begin{array}{cccc:cccc}
					&   &   &  &  &  & -1 &  \mu   \\
					&   &   &  &  & \begin{sideways}$\ddots$\end{sideways} & \begin{sideways}$\ddots$\end{sideways} &    \\
					&   &   &  & -1 & \mu &  &    \\
					&   &   &  & \mu &  &  &   \\
					\hdashline
					&  &   & -\mu & &    &    &    \\
					&  & -\mu  & 1 &   &  &  &    \\
					& \begin{sideways}$\ddots$\end{sideways} & \begin{sideways}$\ddots$\end{sideways}  &  &  &  &  &  \\
					-\mu & 1 &   &  &   &   &    &
				\end{array}
				\right)\in\sp(n,\R).
			\end{equation}
In particular, $\CB^{(1)}_3(\mu,\gamma)=\gamma \begin{pmatrix}
						0 & \mu \\
						- \mu & 0
					\end{pmatrix}$.						
\item [\bf QNF4] For $n\in\N^*$, $\CR(Z)=\CR^{(n)}_4(Z)=\CQ_{\CB^{(n)}_4(\gamma)}+\CL_{l^{(n)}_4(\iota_1)}$,
$\gamma=\pm 1$, $\iota_1\in\R$,  with
					\begin{equation}\label{A4}
						\CB^{(n)}_4(\gamma)=\gamma\left(\begin{array}{cccc:cccc}
							0 &   &   &  &  &  &  &     \\
							-1 & 0  &   &  &  &  &  &    \\
							&  \ddots & \ddots  &  &  &  &  &    \\
							&   &  -1 & 0 &  &  &  &   \\
							\hdashline
							&  &  &  & 0 &  1  &    &    \\
							&  &  &  &   & \ddots & \ddots &    \\
							&  &  &  &  &  & 0 & 1 \\
							&  &  & (-1)^{n+1} &  &   &    & 0
						\end{array}
						\right)\in\sp(n,\R),
					\end{equation}
and $l^{(n)}_4(\iota_1)=\iota_1 {\bf e}_1\in\R^{2n}$, where $\{{\bf e}_1,\cdots, {\bf e}_{2n}\}$ is the canonical basis of column vectors in $\R^{2n}$. In particular, $\CB^{(1)}_4(\gamma)=\gamma\begin{pmatrix}
0 & 0 \\ 1 & 0 \end{pmatrix}$.
					
\item  [\bf QNF5] For $n\in\N^*$ with $n$ odd, $\CR(Z)=\CR^{(n)}_5(Z)=\CQ_{\CB^{(n)}_5}+\CL_{l^{(n)}_5(\iota_1,\iota_2)}$, $\iota_1,\iota_2\in\R$, with
\begin{equation}\label{A5}
						{\CB}^{(n)}_5=\left(\begin{array}{cccc:cccc}
							0 &   &   &  &  &  &  &     \\
							-1 & 0  &   &  &  &  &  &    \\
							&  \ddots & \ddots  &  &  &  &  &    \\
							&   &  -1 & 0 &  &  &  &   \\
							\hdashline
							&  &  &  & 0 &  1  &    &    \\
							&  &  &  &   & \ddots & \ddots &    \\
							&  &  &  &  &  & 0 & 1 \\
							&  &  &  &  &   &    & 0
						\end{array}
						\right)\in\sp(n,\R),
\end{equation} and $l^{(n)}_5(\iota_1,\iota_2)=\iota_1{\bf e}_1+\iota_2 {\bf e}_{2n}\in \R^{2n}$.
					In particular, $\CB^{(1)}_5=\begin{pmatrix}
						0 & 0 \\
						0 & 0
					\end{pmatrix}$ and $l^{(1)}_5(\iota_1,\iota_2)=\iota_1{\bf e}_1$.
					
\item  [\bf QNF6] (Decomposable case) For $n\in\N^*$ with $n\geq 2$, $\CR(Z)=\CR_6^{(n)}(Z)$ is composed as a direct sum of operators from the previously mentioned five cases, each operating in a lower-dimensional space. More precisely, there exists a mutually disjoint partition $\bigcup_{j\in\Lambda}{\CI_j}=\{1,\cdots, n\}$ with $\# \CI_j=m_j<n$ and $k_j\in\{1,\cdots, 5\}$ such that
\begin{equation}\label{R6}
\CR_6^{(n)}(Z)=\sum_{j\in\Lambda}\CR_{k_j}^{(m_j)}(Z_{\CI_j}),\qquad  Z_{\CI_j}:=\begin{pmatrix}
(D_l)_{l\in \CI_j} \\ (X_l)_{l\in \CI_j}\end{pmatrix}.\end{equation}
\end{itemize}
\end{Theorem}

\begin{remark}
The categorization of the homogeneous part in the quantum Hamiltonian is based on the classification of the homogeneous Hamiltonian $h_0(t,\xi,x)$ (as in (\ref{classicalHam}), with $\ell$ and $c$ omitted), which effectively corresponds to classifying the symplectic Lie algebra $\sp(n,\R)$. This classification is detailed by H\"ormander (see Cases a) to g) in Theorem 3.1 of \cite{Hor1995}).
In this paper, we reclassify the homogeneous part in an indecomposable manner, since our objective is to analyze the growth of Sobolev norms, and quantum quadratic Hamiltonians may contain linear terms. The differences between our classification and the one in Theorem 3.1 of \cite{Hor1995} are as follows:
\begin{itemize}
\item Cases a) -- c) in  \cite{Hor1995} correspond precisely to QNF1 -- QNF3 in Theorem \ref{main1}, respectively.
\item Cases d) and e) are categorized into QNF5 and QNF4 for $n=1$, respectively.
\item Cases f) and g) are classified into QNF4 and QNF5 for $n>1$, respectively.
\end{itemize}
Since the QNF classification depends on the symplectic Lie algebra classification, and each QNF is characterized by a Hamiltonian matrix in $\sp(n,\R)$, the normal forms in Theorem \ref{main1} are also referred to as (quantum) symplectic normal forms, as indicated in the title of this paper.
\end{remark}

\begin{remark}
In this paper, the familiar classical-quantum correspondence, typically established via Weyl quantization (refer to (\ref{WeylQuantization})), is slightly modified for the sake of simplicity in expression. For the symbol $h_0(t,\xi,x)$ as defined in (\ref{classicalHam}), the corresponding system of equations of motion is given by $z'=\CA(t)^*z- \J_n \ell(t)$. It is important to note that this alteration does not impact the validity or the formulation of the theorems presented.
\end{remark}
			
\subsection{Classification of growth rates (GRs) of Sobolev norms}\label{sec_gr_norms}
			
Proposed by Bourgain \cite{Bou96} and expanded upon in subsequent pioneering works (e.g., \cite{CKSTT2010, GG2016, GK2015, HPTV2015}), the behavior of solutions to Hamiltonian PDEs in Sobolev space has been a significant topic in Mathematical Physics over the past decades. Particularly, the existence of unbounded trajectories in Sobolev space is intrinsically linked to weak turbulent effects and energy cascades. In recent years, considerable advancements have been made in this area for a variety of quantum Hamiltonian systems.

As discussed in the previous subsection, examining the growth of Sobolev norms of solutions is intimately connected to the reducibility of the original quantum Hamiltonian, especially in the case of linear Hamiltonians. Following \cite{GY00}, Bambusi-Gr\'ebert-Maspero-Robert \cite{BGMR2018} analyzed a $1-$D QHO with a time-periodic linear potential, demonstrating  $t^{s}$ polynomial growth for ${\CH}^s-$norms through reducibility to a transport equation. For a $1-$D QHO with time quasi-periodic perturbations, which are quadratic polynomials of $(D,X)$, various growth rates of ${\CH}^s-$norms have been observed depending on the normal forms of reducibility \cite{LZZ2021, LLZ2022}:  $t^{s}$ polynomial growth for parabolic normal form, exponential growth for hyperbolic normal form, and  $t^{2s}$ polynomial growth when reducible to the Stark Hamiltonian.

It is important to note that almost reducibility (i.e., conjugation to a time-independent system through a sequence of transformations without demonstrating the convergence of this sequence) is also insightful for understanding the long-term behaviors of solutions. For a $1-$D QHO with a time quasi-periodic quadratic perturbation, it was shown in \cite{LZZ2024} that the growth of $\CH^s$-norms can oscillate with an optimal $o(t^s)$ upper bound if the original equation is not reducible. Bambusi and his collaborators employed a similar concept for time-dependent Schr\"odinger equations \cite{BGMR2021, BLM20, BL2022}, focusing on increasingly regularized remaining parts with each iteration, instead of the usual assumption of asymptotic smallness. They achieved a $t^\epsilon$ upper bound for any $\epsilon > 0$ on the Sobolev norms of solutions through this methodology.

Another method often used to achieve unbounded trajectories involves constructing specific perturbations that induce infinite growth. For a $1-$D QHO, Delort \cite{Del2014} (also see Maspero \cite{Mas2018}) designed a time-periodic order-zero pseudo-differential operator as the perturbation, leading to $t^{\frac{s}{2}}$ polynomial growth in ${\CH}^s-$norms. In an abstract setting, Maspero \cite{Mas2021, Mas2023} further explored time-periodic perturbations and provided conditions under which certain solutions exhibit such polynomial growth. For a $2-$D QHO, Faou-Rapha\"el \cite{FR20} embedded Arnold diffusion into an infinite-dimensional quantum system, constructing a time-decaying potential that resulted in logarithmic growth of Sobolev norms over time. Similarly, Thomann \cite{Thomann20}, based on studies of linear Lowest Landau Level equations with a time-dependent potential \cite{ST20}, designed a perturbation that projects onto Bargmann-Fock space, resulting in polynomial growth of Sobolev norms over time for some traveling wave solutions. Bourgain \cite{Bou99} and Haus-Maspero \cite{HM2020} also demonstrated logarithmic growth of Sobolev norms in various settings, including for linear Schr\"odinger equations with quasi-periodic potentials and semiclassical anharmonic oscillators with regular time-dependent potentials.

\medskip

This paper aims to extend beyond the limitations of spatial dimension and explore the long-term behavior of solutions to Eq. (\ref{orig-equ-1}) in arbitrary dimensions $n$. In the $n-$D quadratic Hamiltonian context, the {\it Sobolev space} (also referred to as {\it $\CH^s-$space} for a given $s\geq 0$) is defined as
$$
\H^s(\R^n):=\left\{u\in L^2(\R^n):\CT^{\frac{s}{2}} u \in L^2(\R^n)\right\},
$$
with $\CT$ being the $n-$D QHO defined as in (\ref{nD-QHO}), and equipped with the {\it Sobolev norm} (also known as {\it $\CH^s-$norm})
$\|u\|_{\H^s(\R^n)}:=\left(\|\CT^{\frac{s}{2}} u\|^2_{L^2}+\|u\|^2_{L^2}\right)^{\frac12}$.
Several equivalent forms of the $\CH^s-$norm will be presented in Section \ref{sec_Sobolev}.

Inspired by previous studies, the primary goal of this paper is to provide a comprehensive classification of the growth rates of Sobolev norms for solutions to Eq. (\ref{orig-equ-1}) under the assumption of reducibility. All potential growth rates of Sobolev norms for solutions to Eq. (\ref{orig-equ-1}) will be precisely delineated based on the classification of QNFs obtained in Theorem \ref{main1}.

\begin{Theorem}\label{main2} {\bf (Classification of GRs of Sobolev norms)}
Consider Eq. (\ref{orig-equ-1}) with non-vanishing initial condition $\psi(0)\in \mathcal{H}^s(\R^n)$, $s>0$.
If the system $z'=\CA(t)z+\ell(t)$ is reducible through the conjugation $z=U(t)y+v(t)$, $U\in C_b^1(\R,\Sp(n,\R))$, $v\in C_b^1(\R,\R^{2n})$, then there exists a constant $C>1$,  depending on $s$, $\psi(0)$, $\|U\|_{C_b^0}$ and $\|v\|_{C_b^0}$, and positive functions $g^{(n)}_k(s,t)$, $1\leq k\leq 6$,
 such that, corresponding to QNF1 --- QNF6 in Theorem \ref{main1}, $C^{-1} g^{(n)}_k(s,t) \leq\|\psi(t)\|_{\H^s(\R^n)}\leq C g^{(n)}_k(s,t)$ as $|t|\rightarrow \infty$.
More precisely,
\begin{itemize}
\item [\bf GR1] For $\CR(Z)=\CR^{(n)}_1(Z)$, $g^{(n)}_1(s,t)=|t|^{(n-1)s}e^{\lambda s|t|}$.	
\item [\bf GR2] For $n$ even, $\CR(Z)=\CR^{(n)}_2(Z)$, $g^{(n)}_2(s,t)=|t|^{(\frac{n}{2}-1)s}e^{\lambda_1 s|t|}$.	
\item [\bf GR3] For $\CR(Z)=\CR^{(n)}_3(Z)$, $g^{(n)}_3(s,t)=|t|^{(n-1)s}$.								
\item [\bf GR4] For $\CR(Z)=\CR^{(n)}_4(Z)$, $g^{(n)}_4(s,t)= |t|^{(2n-1)s} +  |\iota_1|^s |t|^{2n s}$.			
\item [\bf GR5] For $n$ odd, $\CR(Z)=\CR^{(n)}_5(Z)$, $g^{(n)}_5(s,t)= |t|^{(n-1)s}  +  (|\iota_1|+ |\iota_2|)^s  |t|^{ns}$.
\item [\bf GR6] (Decomposable case) For $n\geq 2$, $\CR(Z)=\CR^{(n)}_6(Z)$, $g^{(n)}_6(s,t)= \sum_{j\in\Lambda}g^{(m_j)}_{k_j}(s,t)$.
\end{itemize}
\end{Theorem}

\begin{remark} Among the growth rates mentioned, GR1 -- GR5 are uniquely characteristic of the $n-$D setting, in that they are not observed in any Hamiltonian PDE with spatial dimension lower than $n$. To the best of our knowledge, this theorem is the first to describe and comprehensively classify the growth of $\H^s(\R^n)-$norms across arbitrary spatial dimension $n$.
For greater clarity, and particularly to illustrate the decomposable case more concretely, we provide detailed explanations of the quantum normal forms and the corresponding growths of Sobolev norms for the $1-$D and $2-$D scenarios in Appendix \ref{app_2d}.
\end{remark}

Based on the classifications presented in Theorem \ref{main2}, stability of solutions to Eq. (\ref{orig-equ-1}) within the Sobolev space, which is characterized by the boundedness of the solution's Sobolev norm, can be observed exclusively under specific conditions. These conditions are either for QNF3 with $n=1$, or for QNF5 with $n=1$ and both $\iota_1$ and $\iota_2$ equal to zero. This corresponds to:
$$ \CR(Z)=\pm\frac{\mu}{2} \, (D_1^2+X_1^2),\quad \mu>0, \qquad {\rm or}  \qquad  \CR(Z)= 0,$$
or in the case of QNF6, which represents a direct sum of the aforementioned $1-$D QNFs. Consequently, it appears that stability in the Sobolev space manifests primarily as a $1-$D phenomenon.

\begin{Theorem}\label{thm_inverse} Assume that the affine system $z'=\CA(t)z+\ell(t)$ is reducible with $\CA(\cdot)\in C_b^0(\R, \sp(n,\R))$ and $\ell(\cdot)\in C_b^0(\R, \R^{2n})$ as in Eq. (\ref{orig-equ-1}). Given $s>0$, if the solution to Eq. (\ref{orig-equ-1}) satisfies
$\sup_{t\in\R} \|\psi(t)\|_{\H^s(\R^n)}<\infty$, then Eq. (\ref{orig-equ-1}) is reduced to the constant Hamiltonian
$$\CR(Z)=\sum_{j=1}^n c_j (D_j^2+X_j^2),\qquad  c_j\in \R.$$
\end{Theorem}			

Additionally, there exists a notable correlation between classical and quantum Hamiltonians regarding the long-term behaviors of their solutions. This connection facilitates the direct observation of Sobolev norm growth through computations performed on classical systems, thereby eliminating the need to specify which symplectic normal form is responsible for this growth.
\begin{Theorem}\label{thm-ODEPDE}
Assume that the affine system
\begin{equation}\label{OP-O}
z'=\CA(t)z+\ell(t),\qquad z\in\R^{2n},
\end{equation} is reducible.
Let $z^{\rm opt} (t)$ be a solution to (\ref{OP-O}) that exhibits the optimal growth as $|t|\to\infty$ in the sense that, for any solution $z_0(t)$ to (\ref{OP-O}) we have
\beq\label{fast-z*}\sup_{t\in\R}\frac{\|z_0(t)\|}{1+\|z^{\rm opt} (t)\|}<\infty.\eeq
Consequently, for the solution $\psi(t)$ to Eq. (\ref{orig-equ-1}) with $\psi(0)\in \CH^s(\R^n)\setminus\{0\}$, $s>0$,
 there exists a constant $C>1$, depending on $\psi(0)$ and $z^{\rm opt} (0)$, such that
 \begin{equation}\label{estODEPDE}
C^{-1}\|z^{\rm opt} (t)\|^s \leq \|\psi(t)\|_s \leq C\|z^{\rm opt} (t)\|^s.
\end{equation}
\end{Theorem}

\begin{remark} The asymptotic estimate presented in (\ref{estODEPDE}) is anticipated to be valid in more expansive contexts, extending beyond the limits of the reducibility assumption and the quadratic polynomial form of the Hamiltonian. Let us consider two illustrative examples in the $1-$D framework.

In the case of a $1-$D QHO perturbed by a time quasi-periodic homogeneous quadratic form of $(X,D)$, it has been demonstrated in \cite{LZZ2024} that the $\CH^s-$norm of the solution exhibits oscillatory growth with an optimal $o(t^s)-$upper bound under the condition of a {\it non-reducible} quantum Hamiltonian. Correspondingly, for the analogous classical system, a {\it non-reducible} quasi-periodic ${\rm sl}(2,\R)-$linear system, Eliasson has established the existence of oscillatory unbounded solutions with an $o(t)-$upper bound (refer to Theorem A2 and B2 in \cite{Eli1992}).

In the case of a $1-$D QHO perturbed by a specific time-periodic order$-0$ pseudo-differential operator:
  $$ \frac{1}{2\pi {\rm i}}\partial_t \psi=\left(\frac{1}{4\pi }\CT + h^W(t,D,X)\right)\psi, \qquad  h(t,\xi,x)=\epsilon \cos(2t) \eta(\xi,x) \frac{x\xi}{x^2+\xi^2},$$
where $\eta\in C^\infty(\R^2,\R_+)$ is defined such that $\eta(\xi,x)$ equals $0$ when $x^2+\xi^2\leq \frac14$ and $1$ when $x^2+\xi^2\geq \frac12$, and $h^W$ is the Weyl quantization of the classical Hamiltonian $h$ on $\R^2$ (refer to (\ref{WeylQuantization})).
Maspero \cite{Mas2023} has demonstrated that if $\epsilon>0$ is sufficiently small, then $h^W$ acts as a transporter, implying the existence of a solution $\psi(t)\in\CH^s$ for $s>0$, which satisfies $\|\psi(t)\|_s\sim t^\frac{s}{2}$ as $t\to\infty$, within the context of this paper.
Conversely, the corresponding classical system
\beq\label{exampleMaspero}
x'=\xi+\partial_\xi h(t,\xi,x),\qquad \xi'=-x-\partial_x h(t,\xi,x),\qquad x,\xi\in\R, \eeq
permits a solution for $t\geq \epsilon^{-1}$ of the form
  $$x(t)=\sqrt{\epsilon} \cos(t)\sqrt{t+\frac{\sin(4t)}4},\quad \xi(t)=-\sqrt{\epsilon} \sin(t)\sqrt{t+\frac{\sin(4t)}4}.$$
Therefore, the quantity $\sqrt{x^2(t)+\xi^2(t)}$ scales approximately as $(\epsilon t)^\frac{1}{2}$ when $t\to\infty$, representing the optimal growth rate for solutions to the system (\ref{exampleMaspero}).
\end{remark}

\subsection{Strategy of proof}\label{sec_idee}
The classical-quantum correspondence facilitated by Weyl quantization is well-established. More specifically, for any given symbol $h=h(\xi,x)$, where $\xi,x\in\R^n$ and $n\geq 1$, the Weyl operator $h^{W}$ associated with $h$ is defined as follows:
\begin{equation}\label{WeylQuantization}
\left(h^{W} u\right)(x)=\int_{\xi, y \in\R^n} e^{2\pi {\rm i}\la \xi, x-y\ra} \,  h\left(\xi, \frac{x+y}{2}\right) u(y) \, dy  \, d\xi,\qquad \forall \ u\in L^2(\R^n).
\end{equation}
Notably, if $h$ is a quadratic polynomial in $(\xi,x)$, then the corresponding Weyl operator $h^W$ becomes a quadratic polynomial in $(D,X)$ after symmetrization. Therefore, to comprehensively understand Eq. (\ref{orig-equ-1}), it is essential to examine the associated finite-dimensional Hamiltonian system. This approach was utilized in \cite{BGMR2018} to achieve reducibility and the almost preservation of Sobolev norms for the time quasi-periodic perturbed QHO, by implementing a classical KAM scheme for ``most" frequency vectors.
Similarly, in the $1-$D context, the (almost) reducibility and diverse behaviors of Sobolev norms have been explored in \cite{LZZ2021, LZZ2024, LLZ2022}. This was achieved by adapting Eliasson's KAM scheme for time quasi-periodic ${\rm sl}(2,\R)-$linear systems \cite{Eli1992}, where the reducibility transformation is not necessarily close to the identity, to the quantum Hamiltonian.

The reducibility assumption of the affine system
\begin{equation}\label{affine-y}
z'=\CA(t)z+\ell(t),\qquad \CA(\cdot)\in C_b^0(\R, \sp(n,\R)), \qquad \ell(\cdot)\in C_b^0(\R,\R^{2n}),\end{equation}
means that there exist $U(\cdot)\in C_b^1(\R, \Sp(n,\R))$ and $v(\cdot)\in C_b^1(\R,\R^{2n})$ such that, via the transformation $z= U(t) w+v(t)$, the system (\ref{affine-y}) is conjugated to the affine system with constant coefficients $w'=\CB w+l$. Similar to the strategy of reducibility in \cite{LLZ2022}, the reducibility of the quantum Hamiltonian (\ref{orig-equ-1}) is achieved in two steps.

\begin{enumerate}
\item [1)]
Through the transformation $z= U(t) y$, the system (\ref{affine-y}) is conjugated to the affine system
\begin{equation}\label{affine-z}
y'=\CB y+\tilde\ell(t),\qquad \CB\in \sp(n,\R), \qquad \tilde\ell(\cdot)\in C_b^0(\R,\R^{2n}), \ \ {\rm with} \ \ \tilde\ell(t)= U(t)^{-1}\ell(t). \end{equation}
Simultaneously, Eq. (\ref{orig-equ-1}) is conjugated to the quantum Hamiltonian
\begin{equation}\label{PDE-phi}
\frac{1}{2\pi {\rm i}} \partial_t \phi= (\CQ_{\CB} + \CL_{\tilde\ell(t)})\phi,
\end{equation}
via the $L^2-$unitary transformation $ \psi = \CM(U(t)) \phi $, where $\CM(U(\cdot))$ denotes the Metaplectic representation of $U(\cdot)$ (refer to Section \ref{Section_MetaRepre} for details).

\item [2)] By applying the translation $y=w+U(t)^{-1} v(t)$, the system (\ref{affine-z}) is transformed into
$$w'=\CB w+l,\qquad \CB\in \sp(n,\R), \qquad l\in \R^{2n},$$
where all coefficients are time-independent. In parallel, the Hamiltonian PDE (\ref{PDE-phi}) is conjugated to $ \frac{1}{2\pi {\rm i}} \partial_t \varphi= (\CQ_{\CB} +  \CL_{l})\varphi $ (up to a scalar term, typically omitted, see Remark \ref{changshubianhuan}), through the $L^2-$unitary transformation $ \phi = \rho(U(t)^{-1} v(t)) \varphi $. Here, $\rho(U(\cdot)^{-1} v(\cdot))$ represents the Schr\"odinger representation of $U(\cdot)^{-1} v(\cdot)$ (for details, see Section \ref{Section_schroRepre}).
\end{enumerate}
The strategy for proving reducibility is summarized in the following diagram.
$$
\begin{array}{rcccl}
 & {\rm Affine \ system}   & &    {\rm Hamiltonian \ PDE} &  \\
  &   &  &  &   \\
&z'=\CA(t)z+\ell(t) &\stackrel{\rm Quantized}{\longrightarrow}&  \frac{1}{2\pi {\rm i}} \partial_t \psi= (  \CQ_{\CA(t)} +  \CL_{\ell(t)} )\psi & \\
  &   &  &  &   \\
z=U(t) y &\big\downarrow &  & \big\downarrow    &  \psi =  \CM(U(t)) \phi  \\
  &  &  &  &   \\
  &y'=\CB y+\tilde\ell(t)   &\stackrel{\rm Quantized}{\longrightarrow}&  \frac{1}{2\pi {\rm i}}  \partial_t \phi= (  \CQ_{\CB} +  \CL_{\tilde\ell(t)})\phi  &  \\
    &  &  &  &   \\
y=w+ U(t)^{-1} v(t) & \big\downarrow &  & \big\downarrow    & \phi  =   \rho(U(t)^{-1}v(t))  \varphi \\
  &  &  &  & \\
&w'=\CB w+l   &\stackrel{\rm Quantized}{\longrightarrow}&   \frac{1}{2\pi {\rm i}} \partial_t \varphi= (\CQ_{\CB} +  \CL_{l} )\varphi&
  \end{array}
$$

\

The classical-quantum correspondence in this context is achieved through the Metaplectic and Schr\"odinger representations. Section \ref{sec_bounds}, the cornerstone of this article, presents detailed estimates of the Sobolev norms for both representations (see Proposition \ref{boundofmeta} and \ref{boundofschres}). The provided upper and lower bounds demonstrate the equivalence of the Sobolev norm under these representations, which allows us to concentrate on the constant QNFs established in Theorem \ref{main1} and to undertake concrete computations once the quantum Hamiltonian flow is explicitly defined. Moreover, we establish a link between the growth patterns of homogeneous classical and quantum Hamiltonian flows (refer to Corollary \ref{cor_UetA} - (ii)).
Consequently, for QNF1 to QNF3 in Theorem \ref{main1}, which are all homogeneous, only the Metaplectic representation is required and the behavior of solutions in Sobolev space can be directly inferred from the corresponding classical linear system.
For the remaining cases, by employing an inductive approach based on the Baker-Campbell-Hausdorff formula, the quantum Hamiltonian flow is still explicitly derivable (as discussed in Section \ref{sec6}).

\subsection{Index of notations}\label{sec_notations}

\begin{enumerate}
\item With the dimension $n$ fixed, the notations ``$\lesssim$" and ``$\gtrsim$" are used to denote upper and lower bounds, respectively, by a positive constant that depends solely on $s\in \R$ and is independent of any other variables. Specifically, for non-negative values $a$ and $b$, the inequality $a \lesssim b$ (or $a \gtrsim b$) indicates the existence of a constant $c_s$, which is only a function of $s$, such that $a \leq c_s b$ (or $a \geq c_s b$, respectively). The expression $a\simeq b$ signifies that both $a \lesssim b$ and $a \gtrsim b$ hold true. Furthermore, when considering the dependence on other factors, for instance, $u\in \CH^s(\R^n)\setminus\{0\}$, the notation $a \gtrsim_u b$ implies that there is a specific constant $c_{s,u}$, dependent only on $s$ and $u$, ensuring that $a \geq c_{s,u} b$ (refer to (\ref{optimallower-Meta}) in Proposition \ref{boundofmeta} for an example).

\item For $s, s'\in \R$ with $s\geq s'$, let $B({\CH}^{s},{\CH}^{s'})$ be the space of bounded operators from ${\CH}^{s}$ to ${\CH}^{s'}$ with $\|\cdot\|_{B({\CH}^{s},{\CH}^{s'})}$ the operator norm, simplified to $B({\CH}^{s})$ and  $\|\cdot\|_{B({\CH}^{s})}$ if $s'=s$.

\item Let $\Sp(n,\R)$ (resp. $\sp(n,\R)$) be the symplectic group (resp. symplectic Lie algebra) \footnote{See Section \ref{Section_MetaRepre} for the precise definitions.} of $2n\times2n$ real matrices with the matrix norm defined by
\begin{equation}\label{norm_matrix}
\|\A\|=\sum_{i,j=1}^{2n} |a_{ij}|,\quad \A=(a_{ij})_{i,j=1}^{2n}\in\Sp(n,\R) \  \  ({\rm or} \  \   \sp(n,\R)).
\end{equation}
Moreover, we define $\{{\bf e}_1,\cdots, {\bf e}_{2n}\}$ the canonical basis of column vectors in $\R^{2n}$,
the $n-$D identity matrix $I_n$ and the $2n-$D standard symplectic matrix
$$\J_n:=\left(\begin{array}{cc} 0 & I_n \\ -I_n & 0 \end{array}\right)=(-{\bf e}_{n+1}, \cdots, \underset{\underset{n^{\rm th}}{\uparrow}}{-{\bf e}_{2n}}, {\bf e}_{1}, \cdots, {\bf e}_{n})\in \Sp(n,\R),$$
as well as the $2n-$D anti-diagonal symplectic matrix
\beq\label{Janti}
\J^{\rm a}_{n} :=(-{\bf e}_{2n}, \cdots, -\underset{\underset{n^{\rm th}}{\uparrow}}{{\bf e}_{n+1}}, {\bf e}_{n}, \cdots, {\bf e}_{1})\in \Sp(n,\R).\eeq

\item Given $\A_j=\left(\begin{array}{cc}A_j & B_j \\ C_j & F_j \end{array}\right)\in \Sp(n_j,\R)$(resp. $\sp(n_j,\R)$), $j=1,2$, we define the  ``direct sum" of the two matrices by
$$\A=\A_1\bigoplus\A_2:=\left(\begin{array}{cccc}
A_1 & 0 & B_1 & 0 \\
0 & A_2 & 0 & B_2 \\
C_1 & 0 & F_1 & 0 \\
0 & C_2 & 0 & F_2
\end{array}\right)\in \Sp(n_1+n_2,\R) \  \  ({\rm or}  \  \  \sp(n_1+n_2,\R)).$$

\item\label{Cb01} Let $C^0_b(\R, \bullet)$, $\bullet=\Sp(n,\R)$, $\sp(n,\R)$ or $\R^d$, be the space of $\bullet-$valued functions whose entries are all $C^0$ and uniformly bounded functions on $\R$, equipped with the norm
$$\|U\|_{C^0_b}:=\sup_{t\in\R} \|U(t)\|,\qquad U\in C^0_b(\R, \bullet).$$
Let $C^1_b(\R, \bullet)$ be the subspace of $C^0_b(\R, \bullet)$ whose entries are all $C^1$ with uniformly bounded derivatives on $\R$, equipped with the norm
$$\|U\|_{C^1_b}:=\sup_{t\in\R} (\|U(t)\|+\|U'(t)\|),\qquad U\in C^1_b(\R, \bullet).$$

\item As in Eq. (\ref{orig-equ-1}), we define the columns of operators
\begin{equation}\label{DX-n}
D=\begin{pmatrix}D_1\\ \vdots \\ D_n\end{pmatrix}:= \frac{1}{2\pi {\rm i}}\nabla_x,
\quad  X:=\begin{pmatrix}X_1\\ \vdots \\ X_n\end{pmatrix},
\quad Z_j:=\begin{pmatrix}D_j \\ X_j\end{pmatrix},
\quad Z:=\begin{pmatrix}D \\ X\end{pmatrix} ,
\end{equation}
where $X_j$ represents the multiplication by the $j-$th coordinate function, i.e., $(X_jf)(x)=x_jf(x)$. It is important to note that the derivation in $D_j$ always refers to the $j-$th spatial coordinate if other temporal variable is involved. Naturally, for $p\in \R^n$ and $v \in \R^{2n}$, we define
$$\la p, D \ra:= \sum_{j=1}^n p_jD_j,\quad \la p, X \ra:= \sum_{j=1}^n p_jX_j,\quad \CL_{v}(Z)=\left\la v, Z \right\ra:=\sum_{j=1}^n (v_jD_j+v_{n+j}X_j),$$
and for $\CA= \left(\begin{array}{cc} A& B \\ C & -A^* \end{array}\right)\in \sp(n,\R)$,
$$\la Z, \CA \J_n Z \ra:=- \left\la D,B D \right\ra +  \left\la D,AX \right\ra+  \left\la X, A^*D \right\ra+\left\la X,CX \right\ra,$$
where the term $\left\la D,AX \right\ra:=\sum_{i,j=1}^n A_{ij}  D_i X_j $ and similar for the other terms. Let
$\CQ_{\CA}(Z):=-\frac12 \la Z, \CA \J_n Z \ra$,
whose symbol is the homogeneous quadratic polynomial
\begin{equation}\label{homoQuadratic}
\CQ_{\CA}(\xi, x)=\frac12 \left(\la \xi, B\xi\ra-2\la \xi, Ax\ra-\la x, Cx\ra \right),\qquad \xi,x\in\R^n.
\end{equation}

\item For $x\in \R^n$, let $\la x \ra:=\sqrt{1+\la x,x\ra}$. Define the operators $\la X \ra$ and $|X|$ by
$$\left(\langle X\rangle u\right)(x):= \la x \ra u(x),  \qquad \left(|X| u\right)(x):= \sqrt{\la x,x\ra} u(x),$$
and, with ${\CF}^{-1}$ the inverse of Fourier transform, define $\la D \ra$ and $|D|$ by
$$\langle D\rangle u := {\CF}^{-1} \left(\langle X\rangle\hat{u}\right),\quad |D| u := {\CF}^{-1} \left(|X| \hat{u}\right).$$
Moreover, given $\alpha\in \N^n$, we define $D^\alpha:=D_1^{\alpha_1}\cdots D_n^{\alpha_n}$ and $X^\alpha:=X_1^{\alpha_1}\cdots X_n^{\alpha_n}$.
\end{enumerate}

\subsection{Description of the remaining of paper}
			
In the preliminary Section \ref{sec_pre}, we explore the equivalent forms of Sobolev norms, along with an introduction to the Metaplectic and Schr\"odinger representations, highlighting their fundamental properties.

Section \ref{sec_bounds} presents the upper and lower bounds of Sobolev norms within these two representations. This leads to a pivotal conclusion of the paper: all $L^2-$unitary transformations, derived from these representations and conjugating time-dependent equations, maintain the growth rate of the Sobolev norms of solutions. This is articulated in Propositions \ref{boundofmeta} and \ref{boundofschres}. A comprehensive proof of Proposition \ref{boundofmeta} is included in this section, while the proof for Proposition \ref{boundofschres} is deferred to Appendix \ref{app_Schro}.

Section \ref{sec4} delves into a classical-quantum correspondence for the quadratic Hamiltonians, which facilitates the translation of conjugation in classical Hamiltonian systems to quantum Hamiltonians.

In Section \ref{sec5}, we discuss the reducibility of the affine system (\ref{affine-y}) and provide a complete classification of symplectic normal forms (see Proposition \ref{ODEreducibility}). From these findings, Theorem \ref{main1} is derived using the classical-quantum correspondence.

The examination of Sobolev norm growth for each QNF, or each variant of reduced constant quantum Hamiltonian, is conducted in Section \ref{sec6}. This analysis culminates in the proof of Theorems \ref{main2} and \ref{thm-ODEPDE}.

Finally, the application of these theorems in time-periodic and time quasi-periodic quadratic quantum Hamiltonians is explored in Sections \ref{secAppli-per} and \ref{secAppli-qp}.
	

			
\section{Preliminaries}\label{sec_pre}
			
Adhering to the conventions established by Folland \cite{Fol1989}, this section revisits the definition and characteristics of Sobolev norms. Additionally, it offers a review of the Metaplectic and Schr\"odinger representations. Key formulas, essential for subsequent discussions, are introduced here without their proofs.

\subsection{Sobolev space and Sobolev norm}\label{sec_Sobolev}
In the $n-$D setting, $n \geq 2$, with the columns of operators $D$ and $X$ defined as in (\ref{DX-n}), the $n-$D QHO in this paper,
$$\CT=2\pi\left(\la D,D\ra+ \la X, X\ra\right)=-\frac{1}{2\pi} \Delta+2\pi |X|^2,
$$
differs from the conventional harmonic operator ``$-\Delta + |X|^2$."
This distinction leads to variations in the eigenfunctions of the QHO, specifically when normalized. Precisely, the $1-$D Hermite functions
$$h_j(x)=\frac{2^\frac14}{\sqrt{j!}}\left(\frac{-1}{2\sqrt\pi}\right)^je^{\pi x^2}\frac{d^j}{dx^j}\left(e^{-2\pi x^2}\right),\quad x\in\R, \quad j\in \N,$$
which adhere to the recurrence formula
\begin{equation}\label{Hermite-recurr}
h_j(x)=\sqrt{\frac\pi{j}}(X-{\rm i}D)h_{j-1}=\sqrt{\frac{\pi^j}{j!}}(X-{\rm i}D)^jh_{0},
\end{equation}
yield the eigenfunctions of $\CT$ corresponding to the eigenvalue $2|\alpha| + n$ for $\alpha \in \N^n$, given by
$$h_\alpha(x)=h_{\alpha_1}(x_1)\cdots h_{\alpha_n}(x_n).$$
The {\it Sobolev space} (also referred to as the {\it $\CH^s$-space} for a given $s \geq 0$) is defined as
\begin{equation*}\label{Sobo_space-n}
\H^s(\R^n):=\left\{f\in L^2(\R^n):\CT^{\frac{s}{2}} f \in L^2(\R^n)\right\},\quad s\geq 0,
\end{equation*}
equipped with the graph norm, known as the {\it Sobolev norm} or {\it $\CH^s$-norm}:
\begin{equation*}\label{Sobo_norm-n}
\|f\|_{\H^s(\R^n)}:=\left(\|\CT^{\frac{s}{2}} f\|^2_{L^2(\R^n)}+\|f\|^2_{L^2(\R^n)}\right)^{\frac12},
\end{equation*}
and for $s<0$, $\H^{s}(\R^n)=(\H^{-s}(\R^n))'$ by duality. For convenience, we sometimes denote $L^2$ and $\H^s$ instead of $L^2(\R^n)$ and $\H^s(\R^n)$ respectively,
when the dimension of the base space is explicitly stated, and we abbreviate $\|\cdot\|_{\H^s}$ as $\|\cdot\|_{s}$.		
		
Define the {\it Fourier transform} ${\CF}$  on $L^2(\R^n)$ and its inverse ${\CF}^{-1}$ as follows
\begin{eqnarray*}({\CF}f)(\xi) \, = \, \hat{f}(\xi)&=& \int_{\mathbb R^n} e^{- 2\pi {\rm i}\la x, \xi\ra}  f(x) \, dx, \qquad f\in L^2(\R^n), \\
({\CF}^{-1}\hat{f})(x) \, = \, f(x)&=& \int_{\R^n}e^{2\pi {\rm i }\la x, \xi\ra} \hat{f}(\xi) \,  d\xi, \qquad \hat{f}\in L^2(\R^n).
\end{eqnarray*}
Moreover, the {\it partial Fourier transform} ${\CF}_k$, $k=1,\cdots,n$, on $L^2(\R^n)$, and its inverse are defined as
\begin{eqnarray*}
({\CF}_kf)(x_1,\cdots, x_{k-1}, \xi_k,x_{k+1},\cdots,x_n) &=& \int_{\mathbb R} e^{- 2\pi {\rm i} x_k \xi_k}  f(x) \, dx_k, \qquad f\in L^2(\R^n), \\
({\CF}_k^{-1}\hat{f})(\xi_1,\cdots, \xi_{k-1}, x_k,\xi_{k+1},\cdots,\xi_n) &=& \int_{\R}e^{2\pi {\rm i }x_k \xi_k} \hat{f}(\xi) \,  d\xi_k, \qquad \hat{f}\in L^2(\R^n).
\end{eqnarray*}
For given $s\geq 0$, it is well known that $\CF, \CF_k:\CH^s\to \CH^s$, and for $u\in\CH^s$,
\beq\label{Hs-Fourier}
\|u\|_{s}= \|\CF u\|_{s}= \|\CF_k u\|_{s}, \quad \|X_k^s u\|_{L^2} = \|D_k^s \CF_k u\|_{L^2}, \quad \|X_j^s u\|_{L^2} = \|X_j^s \CF_k u\|_{L^2},\quad j\neq k.\eeq
According to \cite{FR20}, for $u \in \CH^s$, $s \geq 0$, the equivalent forms of the $\CH^s-$norm can be expressed as follows:
\begin{equation*}\label{Sobo-equi1}
\|u\|_{s} \simeq\| \langle D\rangle^s u \|_{L^2} + \|\langle X\rangle^s u \|_{L^2}
=\| \langle X\rangle^s \hat{u} \|_{L^2} + \|\langle X\rangle^s u \|_{L^2}
\simeq\|u \|_{L^2}+ \||X|^s \hat{u} \|_{L^2} + \||X|^s u \|_{L^2},
\end{equation*}
and,	 in accordance with \cite{GrePat16}, if $s\in \N^*$, then
$$
\displaystyle \|u\|_{s}\simeq \sup_{\alpha,\beta\in \N^n \atop{|\alpha|+|\beta|\leq s}} \|X^\alpha D^\beta u\|_{L^2}.$$
 Furthermore, for $s>0$, $u \in \CH^s\setminus\{0\}$ is equivalent to one of the following arguments:
\begin{itemize}
\item [(i)] $\|u\|_{s} > 0$,
\item [(ii)] $\|\la X\ra^s u\|_{L^2}> 0$ and $\|\la X\ra^s \CF  u\|_{L^2}>  0$,
\item [(iii)] $\|X_j^s u\|_{L^2}> 0$ and $\|X_j^s\CF u\|_{L^2}> 0$ for $j=1,\cdots,n$.
\end{itemize}
			
\subsection{Schr\"odinger representation}\label{Section_schroRepre}

Let ${\bf H_n}$ represent the real $(2n+1)-$dimensional {\it Heisenberg group}, which is $\R^{2n+1}$ equipped with the group law
$$(p, q, t) (p', q', t') =\left(p+p', q+q', t+t'+\frac12(\la p,q' \ra- \la q, p'\ra)\right),\quad p, q,p',q' \in \R^{n}, \quad t,t'\in \R.$$
The map $\rho$, from ${\bf H_n}$ to the group of unitary operators on $L^2(\R^n)$, is defined as
$$\rho(p,q, t) =e^{2\pi {\rm i}(\la p,D\ra+\la q, X\ra+t{\rm Id})} = e^{2\pi {\rm i }t} e^{2\pi {\rm i}(\la p,D\ra+\la q, X\ra)}, \quad (p,q, t)\in {\bf H_n},$$
meaning, for $f\in L^2(\R^n)$,
\begin{eqnarray}\label{rhode}
\left(\rho(p,q,t) f\right)(x)=e^{2\pi {\rm i}t+2\pi {\rm i}\la q, x\ra +\pi {\rm i} \la p,q\ra} f(x+p) ,
\end{eqnarray}
$\rho$ then becomes an irreducible unitary representation of ${\bf H_n}$ on the Hilbert space $L^2(\R^n)$, known as the {\it Schr\"odinger representation}. Since the variable $t$ always acts in a straightforward manner, it is often convenient to omit it entirely. Consequently, we redefine $\rho$ from $\R^{2n}$ to $B(\CH^s)$, $s\geq 0$, as:
$$
\rho\begin{pmatrix}p \\ q \end{pmatrix} =\rho(p,q,0)=e^{2\pi {\rm i}(\la p,D\ra+\la q, X\ra)},\qquad \forall \ \begin{pmatrix} p \\ q \end{pmatrix}\in \R^{2n},$$
and (\ref{rhode}) is reformulated as
\begin{equation}\label{SchroRep-L2}
\left(\rho\begin{pmatrix} p \\ q \end{pmatrix} f\right)(x) = e^{2\pi {\rm i}\la q, x\ra +\pi {\rm i} \la p,q\ra} f(x+p), \qquad f\in L^2(\R^n).\end{equation}
A straightforward discussion shows that
\begin{eqnarray}\label{schrepmult}
\rho\begin{pmatrix} p_1 \\ q_1 \end{pmatrix}\rho\begin{pmatrix} p_2 \\ q_2\end{pmatrix}=e^{\pi{\rm i}(\la p_1,q_2\ra -\la q_1,p_2\ra)}\rho\begin{pmatrix} p_1+p_2 \\ q_1+q_2
\end{pmatrix},\qquad  p_1,q_1,p_2,q_2 \in \R^{n},
\end{eqnarray}
and it can be easily verified that
\begin{equation}\label{SchroRep-infinitesimal}
\frac{1}{2\pi{\rm i}}\frac{d}{dt}\rho(t\ell)=\left\la \ell, Z\right\ra\rho(t\ell),\qquad \forall \ \ell\in\R^{2n}.\end{equation}

\begin{Proposition} \label{diudiaogaojie}
For $\CB\in \sp(n, \R)$ and $w(\cdot)\in C_{b}^0(\R, \R^{2n})$,
\begin{equation}\label{eq-prop-schro}\left.\partial_\tau\rho\left(e^{\tau\CB}w(t)\right)\right|_{\tau=0}=\left. \partial_\tau\rho\left(w(t)+\tau\CB w(t)\right)\right|_{\tau=0},\qquad \forall \ t\in \R. \end{equation}
\end{Proposition}	
\proof For simplicity, let us express $e^{\tau\CB} w(t) = w(t)+ \tau\CB w(t)+g_\CB(\tau)w(t)$, where $|g_\CB(\tau)|=O(\tau^2)$ as $\tau\to 0$. Based on (\ref{schrepmult}), we obtain $$\rho\left(e^{\tau\CB}w(t)\right)= e^{-\pi {\rm i} \la w(t)+ \tau\CB w(t), \J_n g_\CB(\tau) w(t)\ra} \rho\left(w(t)+ \tau\CB w(t)\right) \rho(g_\CB(\tau) w(t)). $$
Consequently,
			\begin{eqnarray*}
				\left.\partial_\tau\rho\left(e^{\tau\CB}w(t)\right)\right|_{\tau=0} &=&\left.\partial_\tau\left(e^{-\pi {\rm i} \la w(t)+s\CB w(t), \J_n g_\CB(\tau) w(t)\ra}\right)\right|_{\tau=0} \rho(w(t))\\
				& &+ \, \left. e^{-\pi {\rm i} \la w(t)+\tau\CB w(t), \J_n g_\CB(\tau) w(t)\ra}\right|_{\tau=0} \left.\partial_\tau\rho\left(w(t)+\tau\CB w(t)\right)\right|_{\tau=0}\\
				& &+ \, \left. e^{-\pi {\rm i} \la w(t)+\tau\CB w(t), \J_n g_\CB(\tau) w(t)\ra}\right|_{\tau=0} \rho(w(t)) \left.\partial_\tau\rho\left(g_\CB(\tau) w(t)\right)\right|_{\tau=0} \\
				&=:&P_1+P_2+P_3,
			\end{eqnarray*}
where $g_\CB(0)=g_\CB'(0)=0$.
As a result, (\ref{eq-prop-schro}) is established since $P_1=P_3=0$, $P_2=  \left.\partial_\tau\rho\left(w(t)+\tau\CB w(t)\right)\right|_{\tau=0}$.\qed

\smallskip

The following proposition can be established through a similar proof.
\begin{Proposition} \label{diudiaogaojie1}
For $w(\cdot)\in C_{b}^1(\R, \R^{2n})$,
$\left.\partial_\tau\rho(w(t+\tau))\right|_{\tau=0}=\left.\partial_\tau\rho\left(w(t)+\tau w'(t)\right)\right|_{\tau=0}$, $t\in \R$.\end{Proposition}
			
\smallskip
		
Consider $w\in \R^{2n}$. The operator $\rho(w)$ is well-defined on $L^2(\R^n)$, and consequently, on $\CH^s(\R^n)$, $s\geq 0$, as delineated in (\ref{SchroRep-L2}). For $s > 0$, we extend the definition of $\rho(w)$ to $\CH^{-s}(\R^n)$ by employing duality:
\begin{equation}\label{dfforrho}
\la \rho(w)u, v\ra = \la u, \rho(-w) v\ra, \qquad u \in \CH^{-s}(\R^n),  \quad  v\in \CH^s(\R^n).
\end{equation}
					
\begin{Proposition}\label{daoshuforrho}
For $ w(\cdot)\in C_b^1(\mathbb R, \R^{2n})$, it follows that $\partial_t \rho(w(t)) \in B(\H^{s},\H^{s-1})$ for any $s\in\R$.
\end{Proposition}
\noindent The proof of this proposition is provided in Appendix \ref{app_Schro}.

\subsection{Metaplectic representation}\label{Section_MetaRepre}
The {\it symplectic group}, denoted by $ \Sp(n,\R)$, comprises $2n\times 2n$ real matrices that preserve the symplectic form. Specifically, a matrix $\A= \begin{pmatrix} A & B\\ C& F \end{pmatrix}$ belongs to $\Sp(n,\R)$, with $A$, $B$, $C$, and $F$ being its four $n\times n$ blocks, if and only if the conditions $A^{*} C= C^{*} A$, $B^{*}F= F^{*} B$ and $A^* F- C^* B =I_n$ are satisfied. Such a matrix $\A\in \Sp(n,\R)$ is referred to as a {\it symplectic matrix}.
								
The {\it symplectic Lie algebra}, $\sp(n, \R)$, consists of all matrices $\CA\in {\rm gl}(2n,\R)$ such that $e^{t\CA} \in \Sp(n,\R)$ for every $t\in \R$. In other words, $\CA= \begin{pmatrix} A & B\\ C& F \end{pmatrix}$ is in $\sp(n, \R)$, with $A$, $B$, $C$, and $F$ as the four $n\times n$ blocks of $\CA$, if and only if the following hold: $F=-A^*$, $B=B^{*}$, and $C=C^*$. A matrix $\CA\in \sp(n,\R)$ is known as a {\it Hamiltonian matrix}.
		
Let ${\bf Q}$ represent the space of real homogeneous quadratic polynomials on $\R^{2n}$, equipped with the {\it Poisson bracket}. Specifically, for smooth functions $Q_1$ and $Q_2$ on $\R^{2n}$, the Poisson bracket is defined as
$$\left\{Q_1(\xi, x), Q_2(\xi, x)\right\}:= \sum\limits_{j=1}^n \left( \frac{\partial Q_1}{\partial \xi_j} \frac{\partial Q_2}{\partial x_j}-\frac{\partial Q_1}{\partial x_j} \frac{\partial Q_2}{\partial \xi_j}\right). $$
Recall the definition of the homogeneous quadratic polynomial $\CQ_{\CA}(\xi, x)$ given in (\ref{homoQuadratic}).
\begin{Proposition}(Proposition 4.42 in \cite{Fol1989}) The mapping $\CA\mapsto \CQ_{\CA}(\xi, x)$ constitutes a Lie algebra isomorphism from $\sp(n,\R)$ to ${\bf Q}$.\end{Proposition}

The Metaplectic representation $\CM$ of $\Sp(n,\R)$ is defined as follows. Given $\mathbb A\in \Sp(n,\R)$, it induces an automorphism $T_{\mathbb A}$ of the Heisenberg group ${\bf H_n}$ by
$$T_{\mathbb A}(p,q, t)=\left(\mathbb A\begin{pmatrix}p  \\ q \end{pmatrix}, t\right),\qquad p,q\in \R^n, \qquad t\in\R. $$
Since the Schr\"odinger representation $\rho$, defined in (\ref{rhode}), is an irreducible unitary representation of ${\bf H_n}$ on $L^2(\R^n)$, the composite $\rho\circ T_{\mathbb A}$ is also an irreducible unitary representation of ${\bf H_n}$, satisfying
$$(\rho\circ T_{\mathbb A})(0,0, t)=e^{2\pi {\rm i}t},\qquad t\in\R. $$
By the Stone-von Neumann theorem, the two irreducible unitary representations $\rho$ and $\rho\circ T_{\mathbb A}$ are equivalent, meaning there exists a unitary operator $\M(\mathbb A)$ on $L^2(\R^n)$ such that
\begin{equation}\label{hengdeng1}
\rho\circ T_{\mathbb A}(\CX) =\CM(\mathbb A)  \rho (\CX)   \CM(\mathbb A)^{-1}, \qquad \CX\in {\bf H_n}.
\end{equation}
As demonstrated in \cite{Fol1989}, $\CM(\mathbb A)$ can be uniquely chosen, up to factors $\pm 1$, so that $\CM$ becomes a double-valued unitary representation of $\Sp(n,\R)$. Hence,
\begin{equation}\label{product-Meta}
\CM(\mathbb A \mathbb B )=\pm \,  \CM(\mathbb A) \CM(\mathbb B),\qquad \mathbb A, \mathbb B\in \Sp(n,\R),
\end{equation}
and, specifically, $\CM(\mathbb A) \CM\left(\mathbb A^{-1}\right)=\pm \, \Id$.
Therefore, $\CM$ maps $\Sp(n,\R)$ into the group of unitary operators on $L^2(\R^n)$ modulo ${\pm \, \Id}$ and is referred to as the {\it Metaplectic representation} of $\Sp(n,\R)$. In explicit formulations, the ambiguity of $\pm 1$ typically manifests as the ambiguity in the sign of a square root. For simplicity, we often omit the sign $\pm1$.

Given $\mathbb A\in \Sp(n,\R)$, the operator $\CM(\mathbb A)$ is well-defined on $L^2(\R^n)$, and consequently, on $\CH^s(\R^n)$ for $s\geq 0$, as specified in (\ref{hengdeng1}). For $s > 0$, we extend the definition of $\CM(\mathbb A)$ to $\CH^{-s}(\R^n)$ using duality:
\begin{equation}\label{dfHms}
\la \CM(\mathbb A)u, v\ra=\la u, \CM(\mathbb A^{-1}) v\ra, \qquad u \in \CH^{-s}(\R^n),  \quad  v\in \CH^s(\R^n).
\end{equation}
								
Define the infinitesimal version $d\CM$ of the Metaplectic representation as
\begin{equation}\label{defi-infinitesimal}
d\CM(\CA):=\left.\frac{d}{dt} \CM\left(e^{t\CA}\right)\right|_{t=0}, \qquad \CA\in \sp(n, \R).
\end{equation}
In this context and subsequently, the sign of $\CM(e^{t\CA})$ is chosen to ensure continuity in $t$ and equality to $\rm Id$ at $t=0$. Under this framework, the following theorem is established:
\begin{Theorem}(Theorem 4.45 in \cite{Fol1989})\label{Thm4.45}
For any $u$ in the Schwartz space $\CS(\R^n)$ and $\CA\in \sp(n, \R)$, $$ d\CM(\CA)u= 2\pi {\rm i}\CQ_{\CA}(Z) u.$$
\end{Theorem}
\begin{remark}\label{Thm4.45'} The proof of the above theorem shows that for $s\in\R$, $d\CM(\CA)=  2\pi {\rm i}\CQ_{\CA}(Z)\in B(\H^{s},\H^{s-2})$.
\end{remark}

Consider $\CA\in \sp(n, \R)$ and $u\in \H^{s}$. Let us denote $\psi(t,x) =(\CM(e^{t\CA})u)(x)$. Given that $\CM(e^{(t+\tau)\CA}) =\CM(e^{\tau\CA}) \CM(e^{t\CA})$ for any $t,\tau\in\R$, a straightforward computation, in conjunction with Theorem \ref{Thm4.45}, leads us to the following observation:
$$
\partial_t \psi(t, x) =\lim\limits_{\tau\to 0} \frac{1}{\tau}\left(\psi(t+\tau, x)-\psi(t, x)\right)
=\lim\limits_{\tau\rightarrow 0} \frac{1}{\tau} \left( \CM(e^{\tau\CA}) - {\rm Id}\right) \psi(t, x)    = 2\pi {\rm i}\CQ_{\CA}(Z) \psi(t,x),$$
which lies in $\H^{s-2}$. The aforementioned equality implies that for $\CA\in \sp(n, \R)$,
\begin{equation}\label{zizhisolutions}
e^{2\pi{\rm i}\tau \CQ_{\CA}(Z)} =\CM(e^{\tau \CA}),\qquad \forall \ \tau \in \R,
\end{equation}
thereby establishing a connection between the linear system
\begin{equation}\label{ODE-pre}
z'(t) =\CA z(t), \qquad z(0)=z_0\in \R^{2n},
\end{equation}
and its corresponding Hamiltonian PDE
\beq\label{PDE-pre}
\frac{1}{2\pi{\rm i}}\partial_t \psi= \CQ_{\CA}(Z) \psi(t,x) , \qquad \psi(0)=u\in \H^{s}(\R^n).\eeq
It is evident that $z(t)= e^{ t\CA}z_0$ solves the system (\ref{ODE-pre}), while $\psi(t)=\CM(e^{t\CA}) u$ addresses Eq. (\ref{PDE-pre}). Moreover, in light of Theorem 1.2 in \cite{MR2017}, $\psi(t)$ is globally well-posed in $\H^{s}(\R^n)$ and exhibits an upper bound for the $\H^{s}$-norm that grows exponentially with $t$.

\medskip										
										
For specific matrices $\mathbb A\in \Sp(n, \R)$, explicit formulas of the Metaplectic representation can be provided:
\begin{Theorem}(Theorem 4.51 \& 4.53 in \cite{Fol1989})\label{ThmFormMeta}
Let $\mathbb A= \begin{pmatrix}A & B\\C & F\end{pmatrix}\in \Sp(n, \R)$ and $u\in \CH^s$, $s\in \R$.
\begin{itemize}
\item [(i)] If $\det A\neq 0$, then $\displaystyle (\CM(\mathbb A) u)(x) =(\det A)^{-\frac12} \int_{\R^n} e^{2\pi {\rm i} S(\xi,x)} \hat{u}(\xi) \, d\xi$,
where $$S(\xi,x):= -\frac12\la x, CA^{-1} x\ra+\la \xi, A^{-1}x\ra+\frac12\la\xi, A^{-1}B \xi\ra . $$
\item [(ii)] If $\det B\neq 0$, then
$\displaystyle (\CM(\mathbb A) u)(x) ={\rm i}^{\frac{n}{2}}(\det B)^{-\frac12} \int_{\R^n} e^{2\pi {\rm i} S'(x,y)} u(y) \, dy$,
where $$S'(x,y)= -\frac12\la x, FB^{-1} x\ra+\la y, B^{-1}x\ra-\frac12\la y, B^{-1}A  y\ra .$$
\end{itemize}
\end{Theorem}
\begin{remark}
According to \cite{Fol1989}, the aforementioned formulas are valid when $u\in L^2(\R^n)$. Consequently, they also apply to $\CH^{s}(\R^n)$ for $s\geq 0$. Furthermore, by utilizing the definition by duality as presented in (\ref{dfHms}), it can be demonstrated that these formulas remain applicable for $u\in \CH^{-s}(\R^n)$, with $s\geq 0$.
\end{remark}

The oscillatory integrals presented in Theorem \ref{ThmFormMeta} do not provide an explicit description for arbitrary $\mathbb A\in \Sp(n, \R)$. This is due to the prevalence of $2n\times 2n$ symplectic matrices whose $n\times n$ blocks are all singular. This issue will be further explored in Section \ref{sec_bounds}.

\section{Bounds of Sobolev norms under the representations}\label{sec_bounds}
										
It is established that Eq. (\ref{PDE-pre}) is globally well-posed in $\CH^s(\R^n)$ for $s\geq0$. To achieve a more quantitative description of the $\CH^s$-norm, which is the primary objective of this paper, we provide both upper and lower bounds of Sobolev norms. These bounds are determined under the application of the Metaplectic representation for any $\mathbb A\in \Sp(n, \R)$ and the Schr\"odinger representation for any $w\in \R^{2n}$.

\subsection{Bounds of Sobolev norms under Metaplectic representation}
										
Recall from (\ref{zizhisolutions}) that for $\CB(\cdot)\in C_b^0(\R, \sp(n, \R))$ and $|\tau|\leq 1$, $e^{\tau \CB(\cdot)}\in  C_b^0(\R, \Sp(n,\R))$ and satisfies $\M(e^{\tau \CB(t)})=e^{2\pi {\rm i} \tau \CQ_{\CB(t)}}$ for every $t\in\R$. According to Lemma 2.8 in \cite{BGMR2018}, there exist constants $c_1$ and $c_2>0$, dependent on $s\geq 0$ and the uniform bounds of $\CB(\cdot)$, ensuring that
\begin{equation}\label{weakform}
c_1\|u\|_s\le\| e^{2\pi {\rm i} \tau \CQ_{\CB(t)}} u\|_s\le c_2\|u\|_s,  \qquad u\in\H^s, \qquad  t\in\R.\end{equation}
This estimate is particularly effective when the $L^2-$unitary transformation is of the exponential form, especially for near-identity transformations. However, as noted by Djokovi\'c \cite{D1980}, not all elements in $C_b^0(\R, \Sp(n,\R))$ are encompassed in (\ref{weakform}). Specifically, for a connected Lie group $G$, its Lie algebra $g$, and the exponential map $\exp: g\rightarrow G$ of $G$, the {\it index} of an element $a\in G$, denoted by ${\rm ind}(a)$, is the smallest positive integer $m$ such that $a^m$ lies in the image of $\exp$.

\begin{lemma}(Djokovi\'c \cite{D1980})\label{Dlemma} The image of ${\rm ind}$ from $\Sp(n,\R)$ to $\Z$ is $\{1, 2\}$. In other words, for any $\B\in \Sp(n,\R)$, there exists $\mathcal{B}\in \sp(n, \R)$ such that $\B^2=e^\mathcal{B}$.\end{lemma}
Considering Lemma \ref{Dlemma}, for $\mathbb A(\cdot)\in C_b^0(\R, \Sp(n,\R))$, we cannot assure the existence of $\CB(\cdot)\in C_b^0(\R, \sp(n, \R))$ such that $\mathbb A(\cdot)=e^{\CB(\cdot)}$.
To extend (\ref{weakform}), we propose the following proposition, which reveals that only uniform boundedness is required, and even continuity is not a necessity.
							
\begin{Proposition}\label{boundofmeta}
For $\A\in \Sp(n,\R)$, for any $s\in\R$, we have \begin{equation}\label{equiMeta}
\|\A\|^{-|s|}\|u\|_{s} \lesssim \|\M(\A)u\|_{s}\lesssim \|\A\|^{|s|} \|u\|_{s},\qquad \forall \  u\in\H^s.
\end{equation}
Furthermore, for $s>0$ and $u\in\H^s\setminus \{0\}$, it holds that
\begin{equation}\label{optimallower-Meta}
\|\M(\A)u\|_{s}\gtrsim_u \|\A\|^{s}.\end{equation}
\end{Proposition}
									
The estimate (\ref{equiMeta}) is crucial for the $\CH^s-$estimates of the $L^2-$unitary transformations in the reducibility argument, demonstrating the $\CH^s-$equivalence between the original time-dependent equation and the reduced equation with a constant homogeneous part. Conversely, the lower bound (\ref{optimallower-Meta}) is applied to the time propagator of the homogeneous quadratic Hamiltonian PDE via (\ref{zizhisolutions}), essentially contributing to the growth of the Sobolev norm for the solution to the reduced Hamiltonian PDE. More specifically, we arrive at the following corollary.
										
\begin{Corollary}\label{cor_UetA} Given $U\in C_b^0(\R,\Sp(n,\R))$, $\CA\in\sp(n,\R)$, we have, for $t\in\R$,
\begin{itemize}
\item [(i)] for $u\in\H^s$, $s\in\R$, $\|U\|_{C_b^0}^{-|s|}\|u\|_{s} \lesssim \|\M(U(t))u\|_{s}, \ \|\M(U(t)^{-1})u\|_{s}\lesssim \|U\|_{C_b^0}^{|s|} \|u\|_{s}$,
\item [(ii)] for $u\in\H^s\setminus \{0\}$, $s>0$, $\|e^{t\CA}\|^s\lesssim_u \|\M(e^{t\CA})u\|_s\lesssim  \|e^{t\CA}\|^s \|u\|_{s}$.
\end{itemize}
\end{Corollary}

\medskip
									
The remainder of this subsection is dedicated to the proof of Proposition \ref{boundofmeta}.
																
\subsubsection{Bounds for symplectic matrix with diagonal block form}
								
\begin{lemma}\label{towardbasicprop}
For $\A= \begin{pmatrix}A & 0\\0 & A^{*-1}\end{pmatrix}\in\Sp(n,\R)$, we have, for any $s>0$,
\begin{equation}\label{equiMetat}
\|\A\|^{-s} \|u\|_{s} \lesssim \|\M(\mathbb A)u\|_{s}\lesssim \|\A\|^{s}  \|u\|_{s},\quad \forall  \  u\in\H^s.
\end{equation}
Moreover, if $A={\rm diag}\{\zeta_0,\cdots,\zeta_{n-1}\}$, then, for $u\in\H^s\setminus \{0\}$,
\begin{equation}\label{lower-Meta}
\|\M(\mathbb A)u\|_{s}\gtrsim \sum_{j=1}^{n}\left(|\zeta_{j-1}|^s  \|X_{j}^s u\|_{L^2}+|\zeta_{j-1}|^{-s}  \|X_{j}^{s}\hat u\|_{L^2}\right).\end{equation}\end{lemma}
							\proof Since $\A= \begin{pmatrix}
								A & 0\\
								0 & A^{*-1}
							\end{pmatrix}\in\Sp(n,\R)$ implies that $A$ is non-singular, by Theorem \ref{ThmFormMeta}-(i), we have
							$$\left(\M\left(\A\right)u\right)(x)=(\det A)^{-\frac12}\int_{\R^n} e^{2\pi{\rm i} \la \xi, A^{-1}x\ra}\hat{u}(\xi)d\xi=(\det A)^{-\frac12}u\left(A^{-1}x\right),\qquad u\in\H^s,\quad s>0. $$
A straightforward computation on the expression of solution implies that
\begin{eqnarray}\||X|^s\M\left(\A\right)u\|_{L^2}&=&|\det A|^{-\frac12}\left(\int_{\R^n} |x|^{2s}\left|u\left(A^{-1}x\right)\right|^2dx\right)^{\frac12}\nonumber\\
&=&|\det A|^{-\frac12}\left(\int_{\R^n} |Ay|^{2s} |u(y)|^2 \, d\left(Ay\right) \right)^{\frac12}
								 \leq  \|A\|^s\left(\int_{\R^n} |y|^{2s} |u(y)|^2dy \right)^{\frac12},\label{x-new}
							\end{eqnarray}
							 and, if $A={\rm diag}\{\zeta_0,\cdots,\zeta_{n-1}\}$, then
								\begin{equation}\label{x-diagonal}
									\||X|^s\M\left(\A\right)u\|_{L^2}=\left(\int_{\R^n} \left(\sum_{j=1}^{n}\zeta^2_{j-1}y^2_{j}\right)^{s} |u(y)|^2 \, dy\right)^{\frac12}\gtrsim \sum_{j=1}^{n}|\zeta_{j-1}|^s  \|X_{j}^su\|_{L^2}.\end{equation}
On the other hand, since
$$\left(\CF\M\left(\A\right)u\right)(\xi)=(\det A)^{-\frac12}\int_{\R^n} e^{-2\pi{\rm i}\left\la A^*\xi, A^{-1}x\right\ra}u\left(A^{-1}x\right)dx=|\det A|^{\frac12}\hat u\left(A^*\xi\right),$$
it follows that												
\begin{equation}\label{d-new}
\left\||X|^s\CF\M\left(\A\right)u\right\|_{L^2}
=\left(\int_{\R^n} \left|A^{*-1}\eta\right|^{2s}|\hat u(\eta)|^2d\eta \right)^{\frac12}
\leq \|A^{-1}\|^s \left( \int_{\R^n}\left|\eta\right|^{2s} |\hat u(\eta)|^2d\eta \right)^{\frac12},
\end{equation}
and, if $A={\rm diag}\{\zeta_0,\cdots,\zeta_{n-1}\}$, then
\begin{equation}\label{d-diagonal}
\left\||X|^s\CF\M\left(\A\right) u\right\|_{L^2} =\left(\int_{\R^n} \left(\sum_{j=1}^{n} \zeta^{-2}_{j-1}\eta^2_{j}\right)^{s} |\hat u(\eta)|^2 \, d\eta\right)^{\frac12}\gtrsim \sum_{j=1}^{n} |\zeta_{j-1}|^{-s}  \|X_{j}^{s}\hat u\|_{L^2}.\end{equation}
Combining (\ref{x-new}) and (\ref{d-new}), we obtain
$$\left\|\M\left(\A\right) u\right\|_{s}\lesssim\left( \|A\|^s+\|A^{-1}\|^s\right) \| u\|_{s}\lesssim \|\A\|^{s} \| u\|_{s}.$$
The property (\ref{product-Meta}) implies that $\M\left(\A^{-1}\right)\M\left(\A\right)=\Id$. Then, noting that $\|\A^{-1}\|=\|\A\|$, we have
							\begin{equation}\label{uppertolower}
								\| u\|_s=\left\|\M\left(\A^{-1}\right) \M\left(\A\right) u\right\|_{s} \lesssim \|\A\|^{s} \left\| \M\left(\A\right) u\right\|_{s},
							\end{equation}
							which implies that $\left\| \M\left(\A\right) u\right\|_{s}\gtrsim \|\A\|^{-s} \| u\|_{s}$. Moreover, if $A={\rm diag}\{\zeta_0,\cdots,\zeta_{n-1}\}$, then, combining (\ref{x-diagonal}) and (\ref{d-diagonal}), the lower bound (\ref{lower-Meta}) is established.\qed

\subsubsection{Elementary symplectic (ES) matrices}
							
To extend the application of Lemma \ref{towardbasicprop} to any symplectic matrix, we introduce the concept of elementary symplectic matrices within $\Sp(n,\R)$.

A matrix $\A\in\Sp(n,\R)$ is termed {\it elementary symplectic}, abbreviated as {\bf ES} for brevity, if it falls into one of six specific categories. These categories are collectively denoted by $\A\in {\bf ES}:=\bigcup_{{\bf k}={\bf 1}}^{\bf 6} {\bf ES(k)}$.
\begin{itemize}
\item [\bf ES(1).] $\A=\CJ_k$ or $\CJ_k^{-1}$, $k=1, \cdots, n$, where
$$\CJ_k:=({\bf e}_1, \cdots, {\bf e}_{k-1},\underset{\underset{k^{\rm th}}{\uparrow}}{-{\bf e}_{n+k}}, {\bf e}_{k+1}, \cdots,  {\bf e}_n, {\bf e}_{n+1}, \cdots, {\bf e}_{n+k-1},  \underset{\underset{(n+k)^{\rm th}}{\uparrow}}{{\bf e}_{k}}, {\bf e}_{n+k+1}, \cdots, {\bf e}_{2n}).$$
Note that, for $S\in \Sp(n,\R)$ expressed by column vectors, i.e., $S= (\xi_1, \cdots, \xi_n, \eta_1, \cdots, \eta_n)$,
$$S  \CJ_k=(\xi_1, \cdots, \xi_{k-1},\underset{\underset{k^{\rm th}}{\uparrow}}{-\eta_k}, \xi_{k+1}, \cdots,  \xi_n, \eta_1, \cdots,\eta_{k-1},  \underset{\underset{(n+k)^{\rm th}}{\uparrow}}{\xi_k}, \eta_{k+1}, \cdots, \eta_n).  $$
while for $S\in \Sp(n,\R)$ expressed by row vectors, i.e., $S= (\alpha_1, \cdots, \alpha_n, \beta_1, \cdots, \beta_n)^*$,
$$ \CJ_k S=(\alpha_1, \cdots, \alpha_{k-1},\underset{\underset{k^{\rm th}}{\uparrow}}{\beta_k}, \alpha_{k+1}, \cdots, \alpha_n, \beta_1, \cdots,\beta_{k-1},  \underset{\underset{(n+k)^{\rm th}}{\uparrow}}{-\alpha_k}, \beta_{k+1}, \cdots, \beta_n)^*.$$

\item [\bf ES(2).] For $n\geq 2$, $\A=\B_{jk}$ for $1\leq j<k\leq n$, where
$$\B_{jk}:= ({\bf e}_1, \cdots,  \underset{\underset{j^{\rm th}}{\uparrow}}{{\bf e}_{k}},\cdots, \underset{\underset{k^{\rm th}}{\uparrow}}{{\bf e}_{j}}, \cdots, \underset{\underset{(n+j)^{\rm th}}{\uparrow}}{{\bf e}_{n+k}}, \cdots, \underset{\underset{(n+k)^{\rm th}}{\uparrow}}{{\bf e}_{n+j}} \cdots, {\bf e}_{2n}).$$
Note that, for $S\in \Sp(n,\R)$, $S  \B_{jk}$ means the exchanges between $j^{\rm th}$, $k^{\rm th}$ columns and between $(n+j)^{\rm th}$, $(n+k)^{\rm th}$ columns, i.e., for $S= (\xi_1, \cdots, \xi_n, \eta_1, \cdots, \eta_n)$,
$$S  \B_{jk}=(\xi_1, \cdots,\underset{\underset{j^{\rm th}}{\uparrow}}{\xi_k}, \cdots, \underset{\underset{k^{\rm th}}{\uparrow}}{\xi_j}, \cdots,  \xi_n, \eta_1, \cdots, \underset{\underset{(n+j)^{\rm th}}{\uparrow}}{\eta_k}, \cdots, \underset{\underset{(n+k)^{\rm th}}{\uparrow}}{\eta_j}, \cdots, \eta_n),  $$
while $\B_{jk}S$ means the exchanges between $j^{\rm th}$, $k^{\rm th}$ rows and between $(n+j)^{\rm th}-$row and $(n+k)^{\rm th}-$rows, i.e., for $S=(\alpha_1, \cdots, \alpha_n, \beta_1, \cdots, \beta_n)^*$,
$$\B_{jk} S=:(\alpha_1, \cdots,\underset{\underset{j^{\rm th}}{\uparrow}}{\alpha_k}, \cdots, \underset{\underset{k^{\rm th}}{\uparrow}}{\alpha_j}, \cdots,  \alpha_n, \beta_1, \cdots,\underset{\underset{(n+j)^{\rm th}}{\uparrow}}{\beta_k}, \cdots, \underset{\underset{(n+k)^{\rm th}}{\uparrow}}{\beta_j}, \cdots, \beta_n)^*. $$

\item [\bf ES(3).] For $n\geq 2$, $\A=\mathbb E_{q, j}$ for $2 \leq j \leq n$ and $|q|\leq 1$, where $$\mathbb E_{q, j}: =({\bf e}_1, {\bf e}_2, \cdots, \underset{\underset{j^{\rm th}}{\uparrow}}{{\bf e}_j+q {\bf e}_1}, \cdots, {\bf e}_n,  \underset{\underset{(n+1)^{\rm th}}{\uparrow}}{{\bf e}_{n+1}-q{\bf e}_{n+j}}, \cdots, {\bf e}_{2n}). $$
Note that, for $S= (\xi_1, \cdots, \xi_n, \eta_1, \cdots, \eta_n)\in \Sp(n,\R)$, we have
$$S  \mathbb E_{q, j}=(\xi_1, \cdots, \underset{\underset{j^{\rm th}}{\uparrow}}{\xi_j+q \xi_1}, \cdots, \xi_n,  \underset{\underset{(n+1)^{\rm th}}{\uparrow}}{\eta_1-q \eta_j}, \cdots,\eta_j, \cdots,  \eta_n).  $$

\item [\bf ES(4).] For $n\geq 2$, $\A=\mathbb F_{q, j}$ for $2 \leq j\leq n$ and $|q|\leq 1$, where
$$\mathbb F_{q, j}:=({\bf e}^*_1, \cdots, \underset{\underset{j^{\rm th}}{\uparrow}}{{\bf e}^*_j+q{\bf e}^*_{1}}, \cdots, {\bf e}^*_n, \underset{\underset{(n+1)^{\rm th}}{\uparrow}}{{\bf e}^*_{n+1}-q{\bf e}^*_{n+j}}, \cdots, {\bf e}^*_{n+j},\cdots,  {\bf e}^*_{2n})^*.$$
Note that, for $S=(\alpha^*_1, \cdots, \alpha^*_n, \beta^*_1, \cdots, \beta^*_n)^*$, we have
$$\mathbb F_{q, j} S =(\alpha^*_1,  \cdots,  \underset{\underset{j^{\rm th}}{\uparrow}}{\alpha^*_j+q\alpha^*_1}, \cdots, \alpha^*_n, \underset{\underset{(n+1)^{\rm th}}{\uparrow}}{\beta^*_1-q\beta^*_j}, \cdots, \beta^*_j, \dots, \beta^*_n)^*.  $$

\item [\bf ES(5).] For $\A=\mathbb T_{\ell}$ for $1 \leq j \leq n$ and $\ell\in [-n,n]\times  [-1,1]^{n-1}$ (or $[-1,1]$ for $n=1$), where $$\mathbb T_{\ell}: =({\bf e}_1, \dots, {\bf e}_{n}, {\bf e}_{n+1}+\sum\limits_{j=1}^n \ell_j{\bf e}_j, {\bf e}_{n+2}+\ell_2{\bf e}_{1}, \cdots, {\bf e}_{2n}+\ell_n{\bf e}_1).$$
Note that, for $S= (\xi_1, \cdots, \xi_n, \eta_1, \cdots, \eta_n)\in \Sp(n,\R)$, we have
									$$S \mathbb T_{\ell}=
									(\xi_1, \cdots, \xi_n, \eta_1+\sum\limits_{j=1}^n \ell_j\xi_j,\eta_2+\ell_2\xi_1, \cdots,\eta_n+\ell_n\xi_1).$$
																		
\item [\bf ES(6).] $\A=\mathbb S_{\ell}$ for  $\ell\in [-n,n]\times  [-1,1]^{n-1}$, where
$$\mathbb S_{\ell}: =( {\bf e}^*_1, \cdots, {\bf e}^*_{n}, {\bf e}^*_{n+1}+\sum\limits_{j=1}^n \ell_{j}{\bf e}^*_{j}, {\bf e}^*_{n+2}+\ell_2{\bf e}^*_{1}, \cdots, {\bf e}^*_{2n}+\ell_n{\bf e}^*_{1})^*.$$
Note that, for $S=(\alpha^*_1, \cdots, \alpha^*_n, \beta^*_1, \cdots, \beta^*_n)^*$, we have
									$$ \mathbb S_{\ell} S=(\alpha^*_1, \cdots, \alpha^*_n, \beta^*_{1}+\sum\limits_{j=1}^n\ell_j\alpha^*_{j}, \beta^*_{2}+\ell_2\alpha^*_{1}, \cdots, \beta^*_n+\ell_n\alpha^*_{1})^*.$$
								\end{itemize}
For ${\bf k}={\bf 1},\cdots, {\bf 6}$, it is easy to verify that $\A^{-1}\in {\bf ES(k)}$ if $\A\in {\bf ES(k)}$.

\begin{lemma}\label{chudengbianhua0}For $\A\in {\bf ES(k)}$, ${\bf k}={\bf 1},\cdots, {\bf 6}$, for $u \in \CH^s\setminus \{0\}$, $s>0$, we have \begin{equation}\label{equi-SoboNorm}
\left\|\M(\A)u \right\|_{s}\simeq \left\|u\right\|_{s}, \qquad  \|X^s_j\M(\A)u\|_{L^2}, \  \|X^s_j\CF\M(\A)u\|_{L^2} \gtrsim_u 1, \quad j=1,\cdots, n.\end{equation}
\end{lemma}								
\proof Let us show (\ref{equi-SoboNorm}) for every type of {\bf ES} matrix.

For $\A=\CJ_k$ or $\CJ_k^{-1}\in {\bf ES(1)}$, in view of (\ref{Hs-Fourier}), it is sufficient to show $\M(\CJ_k)={\rm i}^\frac12\mathcal{F}_k^{-1}$. Without loss of generality, let us suppose $k = 1$.
For $J_1:=\left(\begin{array}{cccc:cccc}
	 &   &        &   & 1 &   &        &    \\
	&  &        &   &   & 0 &        &    \\
	&   &  &   &   &   & \ddots &    \\
	&   &        &   &   &   &        & 0  \\
	\hdashline
	-1  &   &        &   &  &   &        &    \\
	& 0 &        &   &   &  &        &    \\
	&   & \ddots &   &   &   &  &    \\
	&   &        & 0 &   &   &        &
\end{array}
\right)\in \sp(n,\R)$, noting that
$$e^{t J_1}=\left(\begin{array}{cccc:cccc}
	\cos(t)  &   &        &   & \sin(t) &   &        &    \\
	& 1 &        &   &   & 0 &        &    \\
	&   & \ddots &   &   &   & \ddots &    \\
	&   &        & 1 &   &   &        & 0  \\
	\hdashline
	-\sin(t)  &   &        &   & \cos(t) &   &        &    \\
	& 0 &        &   &   & 1 &        &    \\
	&   & \ddots &   &   &   & \ddots &    \\
	&   &        & 0 &   &   &        & 1
\end{array}
\right)\in\Sp(n,\R),$$
and especially, $e^{\frac\pi2 J_1}=\CJ_1$, we have, by Theorem \ref{Thm4.45}, that
$$d\M(J_1)u=2\pi{\rm i}\CQ_{J_1}(Z)u=\pi{\rm i}(D_1^2+X_1^2)u,\qquad u\in\mathcal{S}(\R^n).$$
Then, for $h_\alpha(x)=h_{\alpha_1}(x_1)\cdots h_{\alpha_n}(x_n)$, with $\alpha=(\alpha_1,...,\alpha_n)\in\N^n$ and $h_k$ the $k-$th Hermite function,
$$d\M(J_1)h_\alpha=\pi{\rm i}(D_1^2+X_1^2)h_{\alpha_1}(x_1)...h_{\alpha_n}(x_n)=\frac{\rm i}2(2\alpha_1+1)h_\alpha.$$
Since $v_{\alpha}(t,\cdot):=\left(\M(e^{t J_1}\right)h_\alpha)(\cdot)$ solves the equation $\partial_t v(t,\cdot)=d\M(J_1)v(t,\cdot)$ with $v(0,\cdot)=h_\alpha(\cdot)$,
we have that $\M(e^{t J_1})h_\alpha=e^{td\M( J_1)}h_\alpha=e^{\frac{{\rm i}t}2(2\alpha_1+1)}h_\alpha$, and hence, for $t=\frac\pi2$, for $\alpha\in\N^n$, $\M\left(\CJ_1\right)h_\alpha={\rm i}^{\alpha_1+\frac12}h_\alpha$.
On the other hand, with $\mathcal{F}h_0=h_0$ and by the recurrence formula (\ref{Hermite-recurr}),
$$\mathcal{F}h_k=\sqrt{\frac{\pi^{k}}{k!}}\mathcal{F}(X-{\rm i}D)^k h_0=\sqrt{\frac{\pi^{k}}{k!}}(-D-{\rm i}X)^k\mathcal{F} h_0=\sqrt{\frac{\pi^{k}}{k!}}(-{\rm i})^{k}(-{\rm i}D+X)^k h_0=(-{\rm i})^{k}h_k. $$
Thus, for every $\alpha\in\N^n$, $\mathcal{F}_1^{-1}h_\alpha={\rm i}^{\alpha_1}h_\alpha={\rm i}^{-\frac12}\M(\CJ_1)h_\alpha$. Since $\{h_\alpha\}$ is a complete orthonormal set, this establishes $\M(\CJ_1)={\rm i}^{-\frac12}\mathcal{F}_1^{-1}$.

For $\A=\B_{jk} \in {\bf ES(2)}$, $1\leq j<k\leq n$, we have that
$$ (\M(\A)u )(x)=u(x_1, \cdots,\underset{\underset{j^{\rm th}}{\uparrow}}{x_{k}}, \cdots, \underset{\underset{k^{\rm th}}{\uparrow}}{x_{j}}, \cdots,  x_n),\quad
(\CF\M(\A)u )(\xi)=\hat  u(\xi_1, \cdots,\underset{\underset{j^{\rm th}}{\uparrow}}{\xi_{k}}, \cdots, \underset{\underset{k^{\rm th}}{\uparrow}}{\xi_{j}}, \cdots,  \xi_n).$$
Hence, $\|X^s_l\M(\A) u \|_{L^2}=\|X^s_l u \|_{L^2}$,  $\|X^s_l\CF\M(\A) u \|_{L^2}=\|X^s_l\CF u \|_{L^2}$ for $l\neq k,j$, and
$$\|X^s_{j}\M(\A) u \|_{L^2}=\|X^s_{k} u \|_{L^2},\quad \|X^s_{k}\M(\A) u \|_{L^2}=\|X^s_{j} u \|_{L^2},$$
$$\|X^s_{j}\CF\M(\A) u \|_{L^2}=\|X^s_{k}\CF u \|_{L^2},\quad \|X^s_{k}\CF\M(\A) u \|_{L^2}=\|X^s_{j}\CF u \|_{L^2},$$
which implies (\ref{equi-SoboNorm}).
							
For $\A= \mathbb E_{q, j}\in {\bf ES(3)}$, we see that $\|\A\|\simeq 1$ since $-1\leq q\leq 1$. Noting that $\A$ is of diagonal block form, in view of Lemma \ref{towardbasicprop}, we have $ \left\|\M(\A)u \right\|_{s}\simeq \left\|u\right\|_{s}$ for $u\in\CH^s$.
A direct computation shows that
$$(\M(\A) u )(x)= u(x_{1}-q x_{j},x_2,\cdots,x_n),\quad
(\CF\M(\A) u )(\xi)= \hat u(\xi_1,\cdots, \underset{\underset{j^{\rm th}}{\uparrow}}{\xi_{j}-q \xi_{1}},\cdots,\xi_n),\quad u\in\CH^s.$$
Hence, $\|X^s_l\M(\A) u \|_{L^2}=\|X^s_l u \|_{L^2}$ for $l\neq 1$,  $\|X^s_l\CF\M(\A) u \|_{L^2}=\|X^s_l\CF u \|_{L^2}$ for $l\neq j$. Moreover,
$$\|X^s_{1}\M(\A) u \|_{L^2}=\left(\int_{\R^n} |x_{1}|^{2s}| u(x_{1}-q x_{j},\cdots,x_n)|^2 dx\right)^{\frac12}
=\left(\int_{\R^n} |y_{1}+q y_{j}|^{2s}| u(y)|^2 dy\right)^{\frac12}  =: F_{ u}(q). $$
Since the function $F_{u}$ is non-vanishing and continuous on $[-1,1]$ provided that $ u\in \CH^s\setminus \{0\}$, we see that $F_{u}(q)\gtrsim_u 1$. Similarly, $\|X^s_{j}\CF\M(\A) u \|_{L^2}\gtrsim_ u 1$, which completes the proof of (\ref{equi-SoboNorm}) for $\A \in {\bf ES(3)}$. The proof is similar for $\A \in {\bf ES(4)}$.
							
For $\A=\mathbb T_{\ell}\in {\bf ES(5)}$, according to Theorem \ref{ThmFormMeta}-(i), we have
$$(\M(\A) u )(x)= \int_{\R^n} e^{2\pi {\rm i} (\la \xi,x\ra+\frac12\la\xi, B_\ell \xi \ra)} \hat{ u}(\xi) \, d\xi,\quad
B_\ell :=\begin{pmatrix}
\ell_1 & \ell_2 &\cdots&\ell_n \\
\ell_2 & 0 & & \\
\vdots &  &\ddots & \\
\ell_n  &  &  & 0\end{pmatrix}.$$
Then $(\CF\M(\A) u)(\xi)= e^{\pi {\rm i} \la\xi, B_\ell \xi \ra} \hat{ u}(\xi)$. For $j=1,\cdots, n$, we have $\|X_j^s\CF\M(\A) u\|_{L^2}= \|X_j^s \CF u\|_{L^2}$,
$$\|X_j^s\M(\A) u\|_{L^2}=\left(\int_{\R^n} |x_j|^{2s}\left|\int_{\R^n} e^{2\pi {\rm i} (\la \xi,x\ra+\frac12\la\xi, B_\ell \xi \ra)} \hat{ u}(\xi) \, d\xi \right|^2 dx\right)^{\frac12}=:G_{j, u}(\ell).$$
Since the function $G_{j, u}$ is non-vanishing and continuous on $[-n,n] \times [-1,1]^{n-1}$ provided that $ u\in \CH^s\setminus \{0\}$, we see that $G_{j, u}(\ell)\gtrsim_ u 1$. Moreover, for $s\in \N^*$, we have
$$ \|X_j^s\M(\A) u\|_{L^2}=\|D_j^s \CF \M(\A) u\|_{L^2}\lesssim  \sum_{k=0}^s \sup_{|\beta|=s-k}\|\la X\ra^k D^\beta \hat{u}\|_{L^2} \lesssim\|u\|_{s}.
							$$
							%
							By interpolation and (\ref{uppertolower}), we obtain (\ref{equi-SoboNorm}) for $\A \in {\bf ES(5)}$. The proof is similar for ${\bf ES(6)}$.\qed

\subsubsection{Proof of (\ref{equiMeta}) for $s>0$.}
The subsequent proposition demonstrates that Lemma \ref{towardbasicprop} can always be applied, up to a finite number of {\bf ES} transformations. By combining this with Lemma \ref{chudengbianhua0}, we can establish the estimates (\ref{equiMeta}) and (\ref{optimallower-Meta}) presented in Proposition \ref{boundofmeta} for the case where $s>0$.
\begin{Proposition}\label{prop-SSdiagonal} {(\bf ES-digonalization of symplectic matrix)} For given $\A\in \Sp(n, \R)$, there exist $\zeta_0,\zeta_1,\cdots,\zeta_{n-1}\in \R$ satisfying
\begin{equation}\label{esti-chi_j}
|\zeta_0|\simeq \|\A\|, \qquad 1 \leq|\zeta_j|\leq  4^{j} \left(\frac{n!}{(n-j)!}\right)^2 \|\A\|,\quad j=0,\cdots,n-1,
\end{equation}
such that, through at most $(2n^2+6n)$ {\bf ES} transformations, $L_1,\cdots L_{k_L},R_1,\cdots R_{k_R}\in {\bf ES}$, \begin{equation}\label{diagonal}
 \left(\prod_{l=k_L}^{1} L_l\right) \A \left(\prod_{j=1}^{k_R} R_j\right)={\bf \Lambda}:={\rm diag}\left\{\zeta_0,\cdots,\zeta_{n-1},\zeta^{-1}_0,\cdots,\zeta^{-1}_{n-1}\right\}.
								\end{equation}
							\end{Proposition}
																			
\proof Let us denote the quantities in (\ref{esti-chi_j}) as
$$M_j:=4^{j} \left(\frac{n!}{(n-j)!}\right)^2 \|\A\|,\qquad j=0,\cdots , n.$$
We shall  ``diagonalize" $\A$ through {\bf ES} transformations, utilizing an induction argument based on the {\bf ES} block-diagonalization process.
Throughout the proof, the superscript ``$j$" (without parentheses), for $j=0,1,\cdots,n$, corresponds to a specific stage in the induction process (but is not an index of power). Furthermore, to streamline the notations, we define the {\bf ES} transformations as follows:
\begin{itemize}
\item $\mathbb A \xrightarrow{R} \mathbb B$ indicates that $\A R= \B$ for some $R\in{\bf ES}$,
\item $\A \xrightarrow[L]{}\B$ signifies that $L\A= \B$ for some $L\in{\bf ES}$.
\end{itemize}

\smallskip
							
\noindent$\bullet$ {\bf First ES block-diagonalisation for $\A$}.
							
\smallskip	
							
Initially, let us ``move" the largest entry of $\A$ to the $(1,1)$-entry through exchanges of columns and rows. For $\A=\A^0=\left(a^0_{ij}\right)_{i,j\in \{1, \cdots, 2n\}}\in\Sp(n,\R)$, let $\chi_0$ be the matrix entry such that $|\chi_0|=\max_{i,j} \left|a^0_{ij}\right|$. In light of the matrix norm definition (\ref{norm_matrix}), it follows that $\frac{\|\A\|}{4n^2} \leq|\chi_0|\leq \|\A\|$. By multiplying $\mathbb J_n$ or $-\mathbb J_n$ (consisting of the product of $n$ matrices in {\bf ES(1)}) to each side of $\A^0$ at most once, which entails exchanges among the four $n\times n$ blocks of $\A^0$ with a possible change of sign, we can obtain a new symplectic matrix, still denoted by $\A^0$, with the block form
$\mathbb A^0= \begin{pmatrix}
A^0_{11} & A^0_{12}\\
A^0_{21} & A^0_{22}\end{pmatrix}$,
ensuring that $\chi_0$ is one of the entries in the upper-left block $A^0_{11}$. Assuming that $\chi_0$ is the $(i_0,j_0)$-entry of $\A^0$ for some $i_0,j_0\in \{1, \cdots, n\} $, we apply the {\bf ES} transformations $\A^0 \xrightarrow{ \B_{1j_0} } \ \xrightarrow[ \mathbb B_{1i_0} ]{} \A_1^{0}$, resulting in $\A_1^0 \in \Sp(n,\R)$ with $\chi_0$ as its $(1,1)$-entry.

Let us now focus on the first row and first column of $\A_1^0$ and symbolically represent $\A_1^0$ as follows:
$$\mathbb A_1^0:  =\left(\begin{array}{cccc:cccc}
\chi_0   & b_2  &  \cdots      & b_n  & g_1 & g_2  & \cdots       &  g_n   \\
a_2 &  &     &    &   &  &       &    \\
\vdots &   &   * &   &   &   & * &    \\
a_n &   &        &  &   &   &       &   \\
\hdashline
a_{n+1} &   &        &   &  &   &         &    \\
\vdots &  &      *   &   &   &  &     *    &    \\
a_{2n} &   &         &  &   &   &      &
\end{array}\right),$$
with ``$*$" denoting other entries which are temporarily irrelevant. We then apply $n-1$ column transformations to $\mathbb A_1^0$:
$$\mathbb A_1^0 \xrightarrow{ \mathbb E_{-b_2/\chi_0,2} } \cdots \xrightarrow{ \mathbb E_{-b_n/\chi_0, n} } \left(\begin{array}{cccc:cccc}
\chi_0   & 0 &  \cdots      & 0 & g_1+\chi_0^{-1}\sum\limits_{j=2}^{n}{b_jg_j} & g_2  & \cdots       &  g_n   \\
a_2 &  &     &    &   &  &       &    \\
\vdots &   &   * &   &   &  * &  &    \\
a_n &   &        &  &   &   &       &   \\
\hdashline
a_{n+1} &   &        &   &  &   &         &    \\
\vdots &  &      *   &   &   & * &         &    \\
a_{2n} &   &         &  &   &   &      &
\end{array}\right)=:\A_2^0,$$
resulting in $\A_2^0\in \Sp(n,\R)$ with the $(1,n+1)$-matrix entry $g':= g_1+\chi_0^{-1}\sum\limits_{j=2}^{n}{b_jg_j}$ satisfying $|g'|\leq n\|\A\|$. A subsequent column transformation:
$$\A_2^0 \xrightarrow{ \mathbb T_{\ell}^{(n)} } \A_3^0:  =\left(\begin{array}{cccc:cccc}
\chi_0   & 0 &  \cdots      & 0 & 0 & 0  & \cdots       &  0  \\
a_2 &  &     &    &   &  &       &    \\
\vdots &   &   * &   &   &   & * &    \\
a_n &   &        &  &   &   &       &   \\
\hdashline
a_{n+1} &   &        &   &  &   &         &    \\
\vdots &  &      *   &   &   &  &     *    &    \\
a_{2n} &   &         &  &   &   &      & \end{array}\right) \quad {\rm with}\quad  \ell:=-\chi_0^{-1}\left(g',g_2,\cdots, g_n\right)$$
yields $\A_3^0$, where every entry falls within $\left[-2n\|\A\|,2n\|\A\|\right]$. In the first column of $\A_3^0$, we have:
$$|\chi_0|, \ |a_2|, \ \cdots, \ |a_{2n}|\leq \|\A\|. $$
Then, using $n-1$ row transformations on ${\mathbb A}^0_3$, we obtain ${\mathbb A}^0_4$:
$${\mathbb A}^0_3 \xrightarrow[\mathbb F_{-a_2/\chi_0, 2}]{}  \cdots \xrightarrow[\mathbb F_{-a_n/\chi_0, n}]{}\left(\begin{array}{cccc:cccc}
\chi_0   & 0 &  \cdots      & 0 & 0 & 0  & \cdots       &  0  \\
0 &  &     &    &   &  &        &    \\
\vdots &   &   * &   &   &   & * &    \\
0 &   &        &  &   &   &         &   \\
\hdashline
a'  &   &          &   &  &   &         &    \\
a_{n+2} &  &         &   &   &  &         &    \\
\vdots &   &    * &   &   &   &  *  &    \\
a_{2n} &   &        &  &   &   &        &
\end{array}\right)  =:{\mathbb A}^0_4$$
with $a':=a_{n+1}+\chi_0^{-1}\sum\limits_{j=2}^n a_ja_{n+j}$. All entries of ${\mathbb A}^0_4$ are within $\left[-4n^2\|\A\|,4n^2\|\A\|\right]$. A further row transformation: $${\mathbb A}^0_4 \xrightarrow[ \mathbb S_{\tau}]{} \left(\begin{array}{cccc:cccc}
\chi_0   & 0 &  \cdots      & 0 & 0 & 0  & \cdots       &  0 \\
0 &  &     &    &   &  &        &    \\
\vdots &   &   * &   &   &   & * &    \\
0 &   &        &  &   &   &         &   \\
\hdashline
0  &   &          &   &  &   &         &    \\
0 &  &         &   &   &  &         &    \\
\vdots &   &  * &   &   &   &  *  &    \\
0 &   &        &  &   &   &        &
\end{array}\right) =:{\mathbb A}_5^0  \quad {\rm with}\quad  \tau:=-\chi_0^{-1}\left(a' , a_{n+2},\cdots, a_{2n}\right), $$
leads to ${\mathbb A}_5^0$, where $\chi_0^{-1}$ is the sole non-vanishing entry in both the $(n+1)^{\rm th}$ row and column, i.e.,
\begin{equation}\label{A50}
{\mathbb A}_5^0= \left(\begin{array}{cccc:cccc}
\chi_0   & 0 &  \cdots      & 0 & 0 & 0  & \cdots       &  0 \\
0 &  &     &    & 0   &  &        &    \\
\vdots &   &   * &   & \vdots  &   & * &    \\
0 &   &        &  & 0   &   &         &   \\
\hdashline
0   & 0 &  \cdots      & 0 &\chi_0^{-1}     & 0  & \cdots       &  0    \\
0  &   & &   &0   &  &         &    \\
\vdots &   &    * &   & \vdots   &   &  *  &    \\
0 &   &        &  &  0 &   &        &
\end{array}\right).\end{equation}
Indeed, by writing ${\mathbb A}_5^0\in \Sp(n, \R)$ in the forms of column vectors and row vectors, i.e.,
$${\mathbb A}_5^0=:(\xi_1, \cdots, \xi_n, \eta_1, \cdots, \eta_n)=:(\alpha_1, \cdots, \alpha_n, \beta_1, \cdots, \beta_n)^*,$$				
by the facts that $\la \xi_1, \mathbb J \xi_j\ra=0$, $\la \xi_1, \mathbb J \eta_j\ra=\delta_{ij}$, and $({\mathbb A}_5^0)^{*}\in \Sp(n, \R)$, one obtains $\beta_1=\chi_0^{-1}{\bf e}_{n+1}$ and $\eta_{1}=\chi_0^{-1}{\bf e}_{n+1}$.
							
Now, from (\ref{A50}), we observe that $\A_5^0$ is already in a block-diagonal form. If $|\chi_0|<1$, then by multiplying $\CJ_1$ and $-\CJ_1$ to each side of $\A_5^0$ respectively, we can facilitate the exchanges between $\chi_0$ and $\chi_0^{-1}$, thereby positioning $\chi_0^{-1}$ as the $(1,1)$-entry. Consequently, it can be seen that, through a maximum of $(4n+4)$ {\bf ES} transformations, $\A=\A^0$ is converted into the first block-diagonal form:
$$\A^0_6= \A_1\bigoplus\A^1\quad  {\rm with} \quad  \A_1=\left(\begin{array}{cc}
								\zeta_0   & 0 \\
								0 &   \zeta_0^{-1}
							\end{array}\right), \quad \A^1\in \Sp(n-1, \R),$$				
where $\zeta_0=\chi_0$ or $\chi_0^{-1}$ and satisfies $\max\left\{\frac{\|\A\|}{4n^2},1 \right\}\leq |\zeta_0|\leq 4n^2\|\A\|=M_1$. Furthermore, every entry of $\A^1$, corresponding to the implicit entries denoted by ``$*$" in (\ref{A50}), falls within $\left[-4n^2\|\A\|,4n^2\|\A\|\right]=[-M_1,M_1]$.		
									
\smallskip
							
\noindent$\bullet$ {\bf Successive ES block-diagonalisations}.
							
\smallskip
											
In pursuit of the {\bf ES} diagonalization (\ref{diagonal}), we repeatedly apply the {\bf ES} block-diagonalization process to reduce the dimension of the non-diagonal block.
				
Suppose that $\A$ has been transformed into $\A_{n-m}\bigoplus \A^{n-m}$, for $1\leq m\leq n-1$, where			
$$\A_{n-m}:={\rm diag}\left\{\zeta_0,\cdots,\zeta_{n-m-1}, \zeta_0^{-1},\cdots,\zeta_{n-m-1}^{-1}\right\},$$
with $\zeta_j$ satisfying $1 \leq|\zeta_j| \leq M_{j+1}$ for $j=0,\cdots,n-m-1$. Here, $\A^{n-m}\in \Sp(m,\R)$ has entries in the range of $[-M_{n-m},M_{n-m}]$. Following the previously outlined {\bf ES} block-diagonalization process, through at most $(4m+4)$ {\bf ES} transformations, $\A^{n-m}$ is converted to
$\left(\begin{array}{cc}
\zeta_{n-m}   & 0 \\0 &   \zeta_{n-m}^{-1}
\end{array}\right) \bigoplus\A^{n-m+1}$,
where $\zeta_{n-m}$ is either the largest matrix entry of $\A^{n-m}$ or its inverse, such that
$1 \leq |\zeta_{n-m}|\leq 4m^2M_{n-m}=M_{n-m+1}$. The matrix $\A^{n-m+1}\in \Sp(m-1,\R) $ has every entry within $\left[-M_{n-m+1},M_{n-m+1}\right]$. Subsequently, $\A$ is transformed into
$$\A_{n-m+1}\bigoplus \A^{n-m+1} \quad {\rm with} \quad \A_{n-m+1}:={\rm diag}\left\{\zeta_0,\cdots,\zeta_{n-m}, \zeta_0^{-1},\cdots, \zeta_{n-m}^{-1}\right\}.$$
Notably, in the final state when $m=1$, $\A$ is transformed into ${\bf \Lambda}$ as specified in (\ref{diagonal}).
Furthermore, there are at most $\sum_{m=1}^n (4m+4)=(2n^2+6n)$ {\bf ES} transformations between $\A$ and ${\bf \Lambda}$.\qed
							
\medskip
							
By synthesizing Proposition \ref{prop-SSdiagonal}, Lemma \ref{towardbasicprop} (specifically, (\ref{equiMetat})), and Lemma \ref{chudengbianhua0}, we establish the bounds (\ref{equiMeta}) in Proposition \ref{boundofmeta} for the case where $s>0$.

Moreover, (\ref{diagonal}) in Proposition \ref{prop-SSdiagonal} suggests that
\begin{equation}\label{SS-diagonal}\A=  \left(\prod^{k_L}_{l=1} L_l^{-1}\right) {\bf \Lambda} \left(\prod^{1}_{j=k_R} R_j^{-1}\right).\end{equation}
Applying (\ref{product-Meta}) and Lemma \ref{chudengbianhua0}, we deduce that
$$\|\CM(\A)u\|_s=\left\|\prod^{k_L}_{l=1} \CM\left( L_l^{-1}\right) \CM({\bf \Lambda})\prod^{1}_{j=k_R}  \CM\left(R_j^{-1}\right) u\right\|_s
\simeq\left\|\CM({\bf \Lambda})\prod^{1}_{j=k_R}  \CM\left(R_j^{-1}\right)u\right\|_s. \ \footnote{Here, the notation `` $\prod$ " represents the compositions of $L^2$-unitary transformations, while in (\ref{SS-diagonal}) it denotes the product of matrices.}$$
Let $\tilde u:=\prod_j \CM\left(R_j^{-1}\right)u$. For $u\in \CH^s\setminus\{0\}$, Lemma \ref{chudengbianhua0} ensures that $\tilde u\in \CH^s\setminus\{0\}$ and $\|X_1^s\tilde u\|_{L^2}\gtrsim_u 1$. Therefore, in conjunction with (\ref{lower-Meta}) in Lemma \ref{towardbasicprop}, it follows that
$$\|\CM(\A)u\|_s\simeq\|\CM({\bf \Lambda})\tilde u\|_s\gtrsim|\zeta_0|^s\|X_1^s\tilde u\|_{L^2}\gtrsim_u \|\A\|^s,$$
which leads to the derivation of the lower bound (\ref{optimallower-Meta}).						
					
\subsubsection{Proof of (\ref{equiMeta}) for $s<0$.}
Fix $s<0$. Given $\mathbb A\in \Sp(n,\R)$, as per (\ref{dfHms}), we have:
\begin{equation}\label{metasinR3}
\left\langle \CM(\mathbb A)u, v\right\rangle =\left\langle u, \CM\left(\mathbb A^{-1}\right)v\right\rangle, \qquad \forall \ u\in \CH^s, \ v\in \CH^{-s}.\end{equation}
Considering that (\ref{equiMeta}) is valid for $\A^{-1}\in \Sp(n,\R)$ and $-s>0$, it follows that:
\begin{equation}\label{metasinR4}
\|\CM(\mathbb A) u\|_{s} =\sup_{\|v\|_{-s}\leq 1} |\langle \CM(\mathbb A) u, v\rangle|\leq\sup_{\|v\|_{-s}\leq 1} \|u\|_{s}\left\|\CM\left(\mathbb A^{-1}\right) v\right\|_{-s}\lesssim \|\A\|^{-s}\|u\|_{s}.
\end{equation}
Furthermore, for $u\in \CH^{s}$ and $v\in \CH^{-s}$, equality (\ref{metasinR3}) implies:
$$\left\la \CM\left(\mathbb A^{-1}\right)\CM(\mathbb A)u, v\right\ra=\left\la \CM(\mathbb A)u, \CM(\mathbb A)v\right\ra=\left\la u, \CM\left(\mathbb A^{-1}\right) \CM(\mathbb A) v\right\ra, $$
which equals to $\la u, v\ra$, considering the fact that $\CM\left(\mathbb A^{-1}\right) \CM(\mathbb A) =\Id$ on $\CH^{-s}$ for $-s>0$.
Thus, for $u\in \CH^{s}$ with $s<0$, we have $\left\|\CM\left(\mathbb A^{-1}\right)\CM(\mathbb A) u\right\|_{s} =\|u\|_{s}$. As demonstrated in (\ref{metasinR4}),
\begin{equation}\label{metasinR5}\|u\|_{s} = \left\|\CM\left(\mathbb A^{-1}\right)\CM(\mathbb A) u\right\|_{s}  \lesssim \|\A\|^{-s} \|\CM(\mathbb A) u\|_{s}.\end{equation}
By combining (\ref{metasinR4}) and (\ref{metasinR5}), we establish (\ref{equiMeta}) for $s<0$.

\subsection{Bounds of Sobolev norms under Schr\"odinger representation}

\begin{Proposition}\label{boundofschres}
 For $ w \in  \R^{2n}$, for $u\in\H^s$, $s\in\R$, $\|w\|^{-|s|} \|u\|_s\lesssim \|\rho(w) u\|_s\lesssim \|w\|^{|s|} \|u\|_s$.
\end{Proposition}

The proof of this proposition is more straightforward than that of Proposition \ref{boundofmeta}. We defer it to Appendix \ref{app_Schro} for brevity.

\begin{Corollary}\label{bianhuanguji2}
Given $w(\cdot)\in C_b^0(\R,\R^{2n})$, for $s\in\R$ and $t\in \R$, we have
$$\|w\|_{C_b^0}^{-|s|}\|u\|_{s} \lesssim \|\rho(w(t))u\|_{s}\lesssim \|w\|_{C_b^0}^{|s|} \|u\|_{s}, \qquad  \forall \  u\in\H^s.$$
\end{Corollary}

\section{Classical-quantum correspondence for quadratic Hamiltonians}\label{sec4}
						
In this section, we will elucidate the relationship between reducibility in classical quadratic Hamiltonian systems and the corresponding Hamiltonian PDEs.

Initially, we introduce the following preparatory lemmas.
\begin{lemma}\label{prelemma} Given $a>0$, for $\mathbb A(\cdot)\in C^1(]-a, a[,\Sp(n,\R))$ satisfying
\begin{equation}\label{jibenjiashe1}
\mathbb A(0)=I_{2n}, \qquad \left.\frac{d}{d\tau}\mathbb A(\tau)\right|_{\tau=0}=0,
\end{equation}
we have that, for $s\in\R$,
\begin{eqnarray}\label{guanjian1}
\left.\partial_\tau\M(\mathbb A(\tau)) u\right|_{\tau=0}=0,\qquad \forall \ u\in \CH^s.
\end{eqnarray}\end{lemma}
\proof
If $\mathbb A(\cdot)=\begin{pmatrix}
A(\cdot) & 0 \\ 0 & A(\cdot)^{*-1}
\end{pmatrix}$ with $A(\cdot)\in C^{1}(]-a, a[, {\rm GL}(n,\R))$ satisfying
$$
A(0)=A(0)^{-1}=I_n,\qquad \left. \frac{d}{d\tau}A(\tau)\right|_{\tau=0}= \left. \frac{d}{d\tau}A(\tau)^{-1}\right|_{\tau=0}=0,$$
then, for $|\tau|$ sufficiently small, we have
$\left(\M(\A(\tau))u\right)(x)=(\det A(\tau))^{-\frac12}u\left(A(\tau)^{-1}x\right)$, and
\begin{eqnarray*}
\left. \partial_\tau \left((\det A(\tau))^{-\frac12}u\left(A(\tau)^{-1}x\right)\right)\right|_{\tau=0}
&=&\left.\frac{d}{d\tau}(\det A(\tau))^{-\frac12}\right|_{\tau=0} u(x)+\left.\partial_\tau u\left(A(\tau)^{-1}x\right)\right|_{\tau=0}\\
&=&2\pi {\rm i}\left\la\left.\frac{d}{d\tau}A(\tau)^{-1}\right|_{\tau=0}x,(D u)(x)\right\ra \ =\ 0.
\end{eqnarray*}
If $\mathbb A(\cdot)=\begin{pmatrix}
I_n & 0 \\ C(\cdot) & I_n
\end{pmatrix}\in\Sp(n,\R)$ with $C(\cdot)=C(\cdot)^*\in C^{1}((-a, a), {\rm gl}(n,\R))$ satisfying
$$C(0)=0,\qquad  \left.\frac{d}{d\tau} C(\tau)\right|_{\tau=0}=0,$$
then, for $|\tau|$ sufficiently small, we have
$\left(\M(\A(\tau))u\right)(x)=e^{-\pi{\rm i}\la x, C(\tau)x\ra}u(x)$, we obtain (\ref{guanjian1}) since
$$\left.\partial_\tau\left(e^{-\pi{\rm i}\la x, C(\tau)x\ra}u(x)\right)\right|_{\tau=0}=-\pi{\rm i}\left\la x, \left.\frac{d}{d\tau}C(\tau)\right|_{\tau=0}x\right\ra u(x)=0.$$

For the general $\mathbb A(\cdot)=\begin{pmatrix}
A(\cdot) & B(\cdot) \\ C(\cdot) &F(\cdot)
\end{pmatrix}\in C^1(]-a, a[,\Sp(n,\R))$ satisfying (\ref{jibenjiashe1}), it follows that
\begin{eqnarray}
A(0)=F(0)=I_n, \qquad B(0)=C(0)=0,\label{mathbbAdaoshuwei0}\\
\left.\frac{d}{d\tau}A(\tau)\right|_{\tau=0}=\left.\frac{d}{d\tau}B(\tau)\right|_{\tau=0}=\left.\frac{d}{d\tau}C(\tau)\right|_{\tau=0}=\left.\frac{d}{d\tau}F(\tau)\right|_{\tau=0}=0.\label{mathbbAdaoshuwei1}
\end{eqnarray}
Considering $A(0)=I_n$, $\mathbb A(\tau)$ can be decomposed for sufficiently small $|\tau|$ as
						$$\mathbb A(\tau)=\begin{pmatrix}
							I & 0 \\
							C(\tau)A(\tau)^{-1} & I
						\end{pmatrix}\begin{pmatrix}
							A(\tau)  & 0 \\
							0 & A(\tau)^{*-1}
						\end{pmatrix}\J_n\begin{pmatrix}
							I & 0 \\
							-A(\tau)^{-1}B(\tau) & I
						\end{pmatrix}\J_n^{-1},$$
and from the above two elementary cases, it follows from the conditions (\ref{mathbbAdaoshuwei0}) and (\ref{mathbbAdaoshuwei1}) that
\begin{eqnarray*}\left.\partial_\tau\M(\mathbb A(\tau))u\right|_{\tau=0}&=&\left.\partial_\tau\M\begin{pmatrix}
								I & 0 \\
								C(\tau)A^{-1}(\tau) & I
							\end{pmatrix}u\right|_{\tau=0}
							+\left.\partial_\tau\M\begin{pmatrix}
								A(\tau) & 0 \\
								0 & A(\tau)^{*-1}
							\end{pmatrix} u \right|_{\tau=0} \\
							& &+ \, \left.\M(\J_n)\partial_\tau\M\begin{pmatrix}
								I & 0 \\
								-A(\tau)^{-1}B(\tau) & I
							\end{pmatrix} \M(\J_n^{-1})u \right|_{\tau=0}
							\ = \  0.	 \qed
						\end{eqnarray*}

\begin{lemma}\label{generalofthm4.45}
Given $Y(\cdot)\in C^1(\R,\Sp(n,\R))$ \footnote{It is easy to verify that $Y^{-1}(\cdot) Y'(\cdot)\in C^0(\R,\sp(n,\R))$ for $Y(\cdot)\in C^1(\R,\Sp(n,\R))$.} and $u\in\H^{s}$, $s\in\R$, we have
\begin{eqnarray}\label{generalofThm4.45}
\frac{1}{2\pi{\rm i}}\partial_t \M(Y(t)) u=\M(Y(t)) \CQ_{Y^{-1}(t) Y'(t)}(Z)u\in \H^{s-2},\quad \forall \ t\in\R.
\end{eqnarray}
\end{lemma}
\begin{remark} Lemma \ref{generalofthm4.45} is perceived as an extension of Theorem \ref{Thm4.45}, wherein the time-parameterized infinitesimal variant of the Metaplectic representation on $\H^{s}(\R^n)$ is considered, applicable for any $s\in \R$.\end{remark}
\proof Recalling (\ref{defi-infinitesimal}) and Remark \ref{Thm4.45'}, we have for $t\in\R$,
$$\CQ_{Y^{-1}(t) Y'(t) }(Z) u= \frac{1}{2\pi{\rm i}} \left.\partial_{\tau}\M\left(e^{\tau Y^{-1}(t) Y'(t) }\right)\right|_{\tau =0}u\in \H^{s-2},\qquad u\in \H^{s},$$
then (\ref{generalofThm4.45}) is equivalent to
\begin{equation}\label{Y1tau}
\left.\partial_{\tau}\M(Y^{-1}(t)Y(t+\tau))\right|_{\tau=0}u=\left.\partial_{\tau}\M\left(e^{\tau Y^{-1}(t) Y' (t)}\right)\right|_{\tau =0}u.\end{equation}
Noting that, for $t,\tau\in\R$,
\begin{eqnarray*}
& &\M\left(Y^{-1}(t)Y(t+\tau)\right)-\M\left(e^{\tau Y^{-1}(t) Y' (t)}\right)\\
& =&\left(\M\left(Y^{-1}(t)Y(t+\tau)e^{-\tau Y^{-1}(t)Y' (t)}\right)-\Id \right) \M\left(e^{\tau Y^{-1}(t) Y' (t)}\right),
\end{eqnarray*}
for obtaining (\ref{Y1tau}), it is sufficient to show
\begin{equation}\label{2ndEqui}
\left.\partial_{\tau}\M\left(Y^{-1}(t)Y(t+\tau)e^{-\tau Y^{-1}(t)Y' (t)}\right)\right|_{\tau=0} u=0, \qquad   \forall  \ u\in \CH^s .
\end{equation}
With any fixed $t\in\R$, define $\mathbb A_t(\tau):= Y^{-1}(t)Y(t+\tau)e^{-\tau Y^{-1}(t)Y' (t)}$, which satisfies
$$\A_t(\cdot) \in C^1(\R,\Sp(n,\R)),\qquad  \A_t(0)=I_{2n},\qquad \left.\frac{d}{d\tau}\mathbb A_t(\tau)\right|_{\tau=0}=0.$$
Applying Lemma \ref{prelemma} to $\A_t(\cdot)$, we obtain (\ref{2ndEqui}).\qed

\medskip

The subsequent proposition can be viewed as analogous to Proposition 2.7 in \cite{BGMR2018}.
\begin{lemma}\label{reducibleprop} Given two Hamiltonian systems
\begin{equation}\label{ODEtype1prop4.4}
z_j'= \CA_j(t) z_j + \ell_{j}(t), \quad \CA_j(\cdot)\in C_{b}^0(\mathbb R, \sp(n, \R)), \quad \ell_j(\cdot)\in C_{b}^0(\R, \R^{2n}), \qquad j=1,2,\end{equation}
if there is $S(\cdot)\in C_b^1(\R, \Sp(n,\R))$ such that $z_1(t)=S(t) z_2(t)$, then, for the two Hamiltonian PDEs
\begin{equation}\label{Heqtype1prop4.4}\frac{1}{2\pi {\rm i}} \partial_t \psi_j=\left(\CQ_{\CA_j(t)}(Z)+\CL_{\ell_{j}(t)}(Z)\right)\psi_j,\qquad j=1,2, \end{equation}
we have the conjugation under the $L^2-$unitary transformation $\psi_1(t)= \M(S(t))  \psi_2(t)$.
\end{lemma}
\begin{remark}\label{rmk-SoboNorm-1}
Regarding the two Hamiltonian PDEs in (\ref{Heqtype1prop4.4}), if the initial conditions $\psi_1(0)$ and $\psi_2(0)$ are both in $\CH^s$, for $s\geq0$, then the solutions are globally well-posed in $\CH^s$. Furthermore, as $S(\cdot)\in C_b^0(\R, \Sp(n,\R))$, in accordance with Corollary \ref{cor_UetA} - (i), the norms $\|\psi_1(t)\|_s$ and $\|\psi_2(t)\|_s$ demonstrate parallel growth rates with respect to time $t$ as $|t|\to\infty$. This is evident from the relation:
$$\|S\|_{C_b^0}^{-s}   \|\psi_2(t)\|_s \lesssim \|\psi_1(t)\|_s\lesssim  \|S\|_{C_b^0}^{s}  \|\psi_2(t)\|_s.$$
\end{remark}

\proof The conjugation $z_1(t)=S(t) z_2(t)$ between two Hamiltonian systems in (\ref{ODEtype1prop4.4}) means that
\begin{equation}\label{conj_ode}
S'(t)=\CA_1( t)S(t) - S(t)\CA_2(t),\qquad \ell_1(t)=S(t) \ell_2(t).
\end{equation}
Since $S(\cdot)\in C_b^1(\R,\Sp(n,\R))$, from Proposition \ref{boundofmeta} and Lemma \ref{generalofthm4.45}, we have
$$\M(S(t))\in  B(\H^{s}),\qquad \frac{d}{dt}\M( S(t))\in  B(\H^{s},\H^{s-2}),\qquad t\in\R. $$
To prove the conjugation $\psi_1(t)= \M( S(t))\psi_2(t)$ between two equations in (\ref{Heqtype1prop4.4}), it is sufficient to prove, in $ B(\CH^s, \CH^{s-2})$, that
\begin{eqnarray}
\frac{1}{2\pi {\rm i}}  \M\left( S(t)\right)^{-1}\frac{d}{dt}\M( S(t))&=&\M\left( S(t)\right)^{-1} \CQ_{\CA_1(t)}(Z) \M( S(t))-\CQ_{\CA_2(t)}(Z),\label{ndterms}\\
\CL_{\ell_2(t)}(Z)&=&\M( S(t))^{-1}\CL_{\ell_1(t)}(Z)\M( S(t)). \label{stterms}
\end{eqnarray}
According to Theorem \ref{Thm4.45} and (\ref{conj_ode}), we have that
\begin{eqnarray*}
\M\left( S(t)\right)^{-1} \CQ_{\CA_1(t)}(Z) \M( S(t))
&=&\frac{1}{2\pi {\rm i}} \M\left( S(t)^{-1}\right) \left. \partial_{\tau}\CM\left(e^{\tau\CA_1(t)}\right)\right|_{\tau=0}\CM\left( S(t)\right)\\
&=&\frac{1}{2\pi {\rm i}}  \left.\partial_{\tau}\CM\left(e^{\tau S(t)^{-1}\CA_1(t) S(t)}\right)\right|_{\tau=0}\\
&=&\frac{1}{2\pi {\rm i}}  \left.\partial_{\tau}\CM\left(e^{\tau\left( S(t)^{-1} S'(t)+\CA_2(t)\right)}\right)\right|_{\tau=0}\\
&=& \CQ_{ S(t)^{-1} S'(t)+\CA_2(t)}(Z)
\ = \ \CQ_{ S(t)^{-1} S'(t)}(Z)+\CQ_{\CA_2(t)}(Z),
\end{eqnarray*}
which means that the right hand side of (\ref{ndterms}) equals to $\CQ_{ S(t)^{-1} S'(t)}(Z)$. Combining with  Lemma \ref{generalofthm4.45}, we obtain (\ref{ndterms}). As for (\ref{stterms}), notice that the property (\ref{hengdeng1}) of the Schr\"odinger representation, as well as the second equality in (\ref{conj_ode}), we have
\begin{eqnarray*}
\M( S(t))^{-1}\CL_{\ell_1(t)}(Z)\M( S(t))
&=&\frac{1}{2\pi{\rm i}}\CM\left( S(t)^{-1}\right)\left.\partial_{\tau}\rho\left(\tau  \ell_1(t)\right)\right|_{\tau=0}\CM( S(t))        \nonumber       \\
&=&\frac{1}{2\pi{\rm i}} \left.\partial_{\tau}\left(\CM\left( S(t)^{-1}\right)\rho\left(\tau \ell_1(t)\right)\CM( S(t))\right)\right|_{\tau=0}             \nonumber  \\
&=&\frac{1}{2\pi{\rm i}} \left.\partial_{\tau}\rho\left(\tau S(t)^{-1} \ell_1(t)\right)\right|_{\tau=0} \ = \ \CL_{\ell_2(t)}(Z).\qed
\end{eqnarray*}

\begin{lemma}\label{reduciblepropnew}
Given two Hamiltonian systems
\begin{equation}\label{ODEtype2prop4.6}
z_j'= \CA z_j+\ell_{j}(t),\quad \CA\in \sp(n, \R), \quad \ell_j(\cdot)\in C_{b}^0(\R, \R^{2n}),  \qquad j=1,2,
\end{equation}
if there is $w(\cdot)\in C_b^1(\mathbb R, \R^{2n})$ such that $z_1(t)=z_2(t) + w(t)$,
then the Hamiltonian PDE
\begin{equation}\label{H_1eqprop4.6}
\frac{1}{2\pi {\rm i}} \partial_t \psi_1=\left(\CQ_{\CA}(Z)+\CL_{\ell_1(t)}(Z)\right)\psi_1
\end{equation}
is conjugated, via the $L^2-$unitary transformation $\psi_1= \rho(w(t)) \psi_2$, to
\begin{equation}\label{H_2eqprop4.6}
\frac{1}{2\pi {\rm i}} \partial_t \psi_2=\left(\CQ_{\CA}(Z)+\CL_{\ell_2(t)}(Z)+\CC(t)\right)\psi_2,\qquad \CC(t):= -\frac 12\la w(t), \J_n(\ell_1(t)+\ell_2(t))\ra.
\end{equation}\end{lemma}
						
\begin{remark}\label{rmk-SoboNorm-2}
In parallel with Remark \ref{rmk-SoboNorm-1}, the condition that $w(\cdot)\in C_b^0(\mathbb R, \R^{2n})$ enables us to employ Corollary \ref{bianhuanguji2}. From this, we can infer that the growth rates of the solutions to Eq. (\ref{H_1eqprop4.6}) and Eq. (\ref{H_2eqprop4.6}) are comparable.
\end{remark}
					
\begin{remark}\label{changshubianhuan}
Eq. (\ref{H_2eqprop4.6}) is transformable into the equation
$$\frac{1}{2\pi {\rm i}} \partial_t\tilde\psi_2 = \left(\CQ_{\CA}(Z)+\CL_{\ell_2(t)}(Z)\right)\tilde\psi_2$$
via the scalar transformation $\psi_2 = e^{2\pi {\rm i}\int_0^{t}\CC(\tau)d\tau} \tilde\psi_2$. Consequently, the scalar part is often omitted when considering the norm of the solution to Eq. (\ref{H_2eqprop4.6}).
\end{remark}

\proof The conjugation $z_1(t) = z_2(t) + w(t)$ between the two systems in (\ref{ODEtype2prop4.6}) implies that
\begin{equation}\label{1conj_ODE}
\ell_2(t)=\CB w(t)-w'(t)+\ell_1(t),\qquad \forall \ t\in\R.
\end{equation}
For $w(\cdot)\in C_b^1(\R,\R^{2n})$, Propositions \ref{daoshuforrho} and \ref{boundofschres} lead to
$$ \rho(w(t))\in  B(\H^{s}),\quad \partial_{t}\rho(w(t)) \in  B(\H^{s},\H^{s-1}),\qquad t\in\R. $$
To establish the conjugation $\psi_1(t) = \rho(w(t))\psi_2(t)$ between Eq. (\ref{H_1eqprop4.6}) and Eq. (\ref{H_2eqprop4.6}), it suffices to demonstrate, within $B(\CH^s, \CH^{s-2})$, that
$$\rho(w(t))^{-1}\left(\CQ_{\CB}(Z) \rho(w(t))-\frac{1}{2\pi{\rm i}}\frac{d}{dt}\rho(w(t))+\CL_{\ell_1(t)}(Z)\rho(w(t))\right)=\CQ_{\CB}(Z)+\CL_{\ell_2(t)}(Z)+\CC(t).$$
Since, through (\ref{hengdeng1}), we have
\begin{eqnarray*}
\rho(w(t))^{-1}\CQ_{\CB}(Z)\rho(w(t))&=&\frac{1}{2\pi{\rm i}} \left.\partial_\tau\left(\rho(w(t))^{-1}\M\left(e^{\tau\CB}\right)\rho(w(t))\M\left(e^{-\tau\CB}\right)\M\left(e^{\tau\CB}\right)\right)\right|_{\tau=0}\\
&=&\frac{1}{2\pi{\rm i}} \left.\partial_\tau \left(\rho(w(t))^{-1}\rho\left(e^{\tau\CB} w(t)\right)\M\left(e^{\tau\CB}\right)\right)\right|_{\tau=0}\\
&=&\frac{1}{2\pi{\rm i}}\left.\left(\rho(w(t))^{-1}\rho\left(e^{\tau\CB}w(t)\right)\right)\right|_{\tau=0}+\frac{1}{2\pi{\rm i}} \left.\frac{d}{d\tau} \M\left(e^{\tau\CB}\right)\right|_{\tau=0}\\
&=&\frac{1}{2\pi{\rm i}}\left. \rho(w(t))^{-1} \partial_\tau\rho(e^{\tau\CB}w(t))\right|_{\tau=0}+ \CQ_\CB(Z),
\end{eqnarray*}
it remains to show that
\begin{equation}\label{1orderterms}
\rho(w(t))\left(\CL_{\ell_2(t)}(Z)+\CC(t)\right)=\frac{1}{2\pi {\rm i}}\left.\partial_\tau\rho\left(e^{\tau\CB}w(t)\right)\right|_{\tau=0}-\frac{1}{2\pi {\rm i}}\frac{d}{dt}\rho(w(t))+\CL_{\ell_1(t)}(Z)\rho(w(t)).\end{equation}
Indeed, following Proposition \ref{diudiaogaojie} and \ref{diudiaogaojie1}, the terms in the right hand side of (\ref{1orderterms}) satisfy that
\begin{eqnarray*}
\CL_{\ell_1(t)}(Z)\rho(w(t))&=&\frac{1}{2\pi {\rm i}}\left.\partial_\tau\left( \rho(\tau\ell_1(t))\rho(w(t))\right)\right|_{\tau=0},\\
\frac{d}{dt}\rho(w(t)) \ = \  \left.\partial_\tau \rho(w(t+\tau))\right|_{\tau=0}&=&\left.\partial_\tau \rho\left(w(t)+\tau w'(t)\right)\right|_{\tau=0},\\
\left.\partial_\tau\rho\left(e^{\tau\CB}w(t)\right)\right|_{\tau=0}&=&\left.\partial_\tau\rho\left(w(t)+\tau\CB w(t)\right)\right|_{\tau=0},\end{eqnarray*}
then, through (\ref{schrepmult}) and (\ref{1conj_ODE}), the right hand side of (\ref{1orderterms}) equals to
\begin{eqnarray*}
& &\frac{1}{2\pi {\rm i}}\left(\left. \partial_\tau \left(e^{2\pi{\rm i}\tau \la w(t), \J_n\ell_1(t)\ra} \rho(w(t)) \rho(\tau\ell_1(t))\right) \right|_{\tau=0}-\left. \partial_\tau \left( e^{-\pi{\rm i}\tau  \la w(t),  \J_n w'(t)\ra  }\rho(w(t))\rho\left(\tau w'(t)\right)\right)\right|_{\tau=0}\right.\\
& & \  \  \  \  \  \   \  \   \  \   +   \left.  \left. \partial_\tau \left( e^{-\pi{\rm i}\tau\la w(t),  \J_n \CB w(t)\ra }\rho(w(t))\rho\left(\tau\CB w(t)\right)\right) \right|_{\tau=0}\right)\\
&=&\frac{1}{2}\left( 2\la w(t), \J_n\ell_1(t)\ra + \la w(t),  \J_n w'(t)\ra-\la w(t),  \J_n \CB w(t)\ra   \right) \rho(w(t))\\
& &+  \,  \frac{1}{2\pi {\rm i}} \rho(w(t))  \left. \partial_\tau\left(\rho\left(\tau\CB w(t)\right)-\rho\left(\tau w'(t)\right)+\rho(\tau\ell_1(t))\right)\right|_{\tau=0} \\
&=&-\frac{1}{2}\left\la w,  \J_n(\CB w-w'(t)+2\ell_1(t))\right\ra   \rho(w(t)) +\rho(w(t))  \left\la \CB w(t)-w'(t)+ \ell_1(t), Z \right\ra \\
&=&-\frac{1}{2}\left\la w(t),  \J_n (\ell_1(t)+\ell_2(t))\right\ra\rho(w(t))+\rho(w(t))\left\la \ell_2(t), Z\right\ra\\
&=& \rho(w(t))\left(\CL_{\ell_2(t)}(Z)-\frac{1}{2}\left\la w(t),  \J_n(\ell_1(t)+\ell_2(t))\right\ra \right).
\end{eqnarray*}
Hence Lemma \ref{reduciblepropnew} is shown. \qed

\medskip

Combining Lemma \ref{reducibleprop}, Lemma \ref{reduciblepropnew} and Remark \ref{changshubianhuan}, we have
\begin{Proposition}\label{MetaandSch}
If the affine system $z'= \CA(t) z + \ell(t)$ with $\CA(\cdot)\in C_{b}^0(\mathbb R, \sp(n, \R))$ and $\ell(\cdot)\in C_{b}^0(\R, \R^{2n})$
is reducible, i.e., it is conjugated to $y'= \CB y + l$, $\CB\in \sp(n, \R)$, $l\in \R^{2n}$, through
some transformation
$$z(t)=U(t) y(t)+ v(t),\qquad  U(\cdot)\in  C_b^1(\R, \Sp(n,\R)),\qquad v(\cdot)\in C_b^1(\mathbb R, \R^{2n}),$$
then, for the two Hamiltonian PDEs
$$
\frac{1}{2\pi {\rm i}} \partial_t \psi_1=\left(\CQ_{\CA(t)}(Z)+\CL_{\ell(t)}(Z)\right)\psi_1,\qquad \frac{1}{2\pi {\rm i}} \partial_t \psi_2=\left(\CQ_{\CB}(Z)+\CL_{l}(Z)\right)\psi_2,
$$
we have $\psi_1(t)= e^{2\pi {\rm i}\int_0^{t}\CC(\tau)d\tau} \rho(v(t)) \M( U(t))  \psi_2(t)$ with for some $\CC(\cdot)\in C_b^0(\R, \R)$.\end{Proposition}
\begin{remark} In light of the diagram presented in Section \ref{sec_idee}, the aforementioned $L^2$-unitary transformation can be expressed, modulo a scalar coefficient, as $\psi_1(t) = \M(U(t))\rho(U(t)^{-1}v(t))\psi_2(t)$. In fact, according to (\ref{hengdeng1}), it follows that
$\M( U(t))\rho(U(t)^{-1} v(t)) =  \rho(v(t)) \M( U(t))$.
\end{remark}
																		
\section{Quantum symplectic normal forms -- Proof of Theorem \ref{main1}}\label{sec5}
						
Let us return to the time-dependent quadratic quantum Hamiltonian as presented in (\ref{orig-equ-1}). Considering Remark \ref{changshubianhuan}, we set aside the scalar component and concentrate on the equation
\begin{equation}\label{orig-equ1}
\frac{1}{2\pi {\rm i}}\partial_t \psi=(\CQ_{\mathcal{A}(t)}(Z)+\CL_{\ell(t)}(Z))\psi=\left(-\frac12\la Z,\mathcal{A}(t)\mathbb JZ\ra+\la \ell(t), Z\ra\right) \psi, \quad x \in\R^n,
\end{equation}
where $\mathcal{A}(\cdot)\in C_b^0(\R, \sp(n,\R))$ and $\ell(\cdot)\in C_b^0(\R, \R^{2n})$.
Our objective is to establish the reducibility framework for Eq. (\ref{orig-equ1}) through the classical-quantum correspondence delineated in Section \ref{sec4}, thereby providing a comprehensive categorization of quantum symplectic normal forms.	
Under the assumption of reducibility for the linear system
\begin{equation}\label{nonHomoODE}z'=\mathcal{A}(t)z+ \ell(t), \end{equation}
we identify $\tilde\CB\in \sp(n,\R)$, $\tilde l\in \R^{2n}$, and transformations $S(\cdot)\in C_b^1(\R,\Sp(n,\R))$, $r(\cdot)\in C_b^1(\R,\R^{2n})$, such that the system (\ref{nonHomoODE}) is conjugated to $y'=\tilde\CB y+\tilde l$ under the change of variable $z=S(t) y+r(t)$. In light of Proposition \ref{MetaandSch}, our focus is directed towards the classification of the reduced classical system $y'=\tilde\CB y+\tilde l$.

It was shown by H\"ormander \cite{Hor1995} that any Hamiltonian matrix can be conjugated to the direct sum of the normal forms
as delineated in (\ref{A1}) -- (\ref{A5}). We will now present this theorem in the context of Hamiltonian matrices.
\begin{lemma}\cite{Hor1995}\label{lem-Hormander}
For any $\tilde\CB\in\sp(n,\R)$, there exist $T\in\Sp(n,\R)$ and a mutually disjoint partition $\bigcup_{j\in\Lambda}{\CI_j}=\{1,\cdots, n\}$ with $\#   \CI_j=m_j>0$ and $k_j\in\{1,\cdots, 5\}$, such that
\beq\label{hormander}
T^{-1}\tilde\CB T= \bigoplus_{j\in\Lambda}  \CB_{k_j}^{(m_j)}.
\eeq
If $\# \Lambda=1$, then $T^{-1}\tilde\CB T = \CB_{k}^{(n)}$ with $\CB_{k}^{(n)}$, $k=1,\cdots,5$, introduced in (\ref{A1}) -- (\ref{A5}). Specifically,
\begin{itemize}
\item [\it Case 1.] If $\tilde\CB$ has two real eigenvalues $\pm \lambda$, $\lambda>0$, each with multiplicity $n$, then $T^{-1}\tilde\CB T = \CB_{1}^{(n)} = \CB_{1}^{(n)}(\lambda)$.
\item [\it Case 2.]  For $n$ even, if $\tilde\CB$ has four complex eigenvalues $\pm\lambda_1\pm{\rm i} \lambda_2$, $\lambda_1,\lambda_2>0$, each with multiplicity $\frac{n}{2}$, then $T^{-1}\tilde\CB T = \CB_{2}^{(n)} = \CB_{2}^{(n)}(\lambda_1,\lambda_2)$.
\item [\it Case 3.] If $\tilde\CB$ has two purely imaginary eigenvalues $\pm{\rm i} \mu$, $\mu>0$, each with multiplicity $n$, then $T^{-1}\tilde\CB T = \CB_{3}^{(n)} = \CB_{3}^{(n)}(\mu,\gamma)$.
\item [\it Case 4.] If $\tilde\CB^{2n} = 0$ and $\tilde\CB^{2n-1} \neq 0$, then $T^{-1}\tilde\CB T = \CB_{4}^{(n)} = \CB_{4}^{(n)}(\gamma)$.
\item [\it Case 5.] For $n$ odd, if $\tilde\CB^{n} = 0$ and, for $n \geq 3$, ${\rm Rank}(\tilde\CB^{n-1}) = 2$, then $T^{-1}\tilde\CB T = \CB^{(n)}_{5}$.
\end{itemize}
\end{lemma}
\begin{remark}\label{rmk-B5even}
One may naturally question why $\CB^{(n)}_{5}$, which remains definable as in (\ref{A5}), is not included in the classification for the case when $n$ is even. Indeed, as elucidated in \cite{Hor1995}, such a matrix can be conjugated to $\CB^{(\frac{n}{2})}_{4}(1)\bigoplus\CB^{(\frac{n}{2})}_{4}(-1)$. This result categorizes it into the decomposable case.
\end{remark}

\medskip

Given that the full dimension $n$ is arbitrarily specified, our attention turns to the simplified system
\begin{equation}\label{sys-Bl-const}
y'=\CB y+\tilde l,\qquad\CB=\CB_{k}^{(n)},\quad  k=1,\cdots,5.
\end{equation}
For $k=1,2,3$, where $\CB=\CB^{(n)}_k$ is invertible, the system (\ref{sys-Bl-const}) can be conjugated to $w'=\CB w$ via the translation $y=w-\CB^{-1} \tilde l$.
For $k=4,5$, with $\CB=\CB^{(n)}_k$ being nilpotent, the systems
$$\CB^{(n)}_4 f_4=\tilde l- \la \tilde l, {\bf e}_1\ra {\bf e}_1,\qquad \CB^{(n)}_5 f_5=\tilde l- \la \tilde l, {\bf e}_1\ra {\bf e}_1- \la \tilde l, {\bf e}_{2n}\ra {\bf e}_{2n}, $$
admit solutions $f_4$ and $f_5$, respectively. Therefore, for $k=4,5$, the system (\ref{sys-Bl-const}) is conjugated to $w'=\CB^{(n)}_k+l^{(n)}_k$ using the translation $y=w-f_k$. Here, $l^{(n)}_k$ is defined as $l^{(n)}_k:= \left\{\begin{array}{ll}
\la \tilde l, {\bf e}_1\ra {\bf e}_1,&k=4\\
\la \tilde l, {\bf e}_1\ra {\bf e}_1+ \la \tilde l, {\bf e}_{2n}\ra {\bf e}_{2n},&k=5
\end{array}\right.$. Moreover, for the special case of $n=1$ and $k=5$, where $\CB^{(1)}_5=\begin{pmatrix}
0&0\\
0&0
\end{pmatrix}$, the reduced system $w'=l^{(1)}_5=\begin{pmatrix}
\iota_1\\
\iota_2
\end{pmatrix}$, given that $(\iota_1,\iota_2)\neq (0,0)$,  can be further conjugated to $\tilde w'=\begin{pmatrix}
\sqrt{\iota_1^2+\iota_2^2}\\
0
\end{pmatrix}$ under the rotation $w=\frac{1}{\sqrt{\iota_1^2+\iota_2^2}}\begin{pmatrix}
\iota_1&-\iota_2\\
\iota_2&\iota_1
\end{pmatrix}\tilde w$.

\medskip

To summarize, the following proposition elucidates Theorem \ref{main1} in light of Proposition \ref{MetaandSch}.
\begin{Proposition}\label{ODEreducibility} {\bf (Classical normal forms)} If the system (\ref{nonHomoODE}) is reducible, then
there exists $U(\cdot)\in C_b^{1}(\R, \Sp(n,\R))$ and $v(\cdot)\in C_b^{1}(\R, \R^{2n})$ such that, via the transformation $z(t)=U(t) w(t)+v(t)$, the system (\ref{nonHomoODE}) is conjugated to
\begin{equation}\label{indecom-constant}
w'=\CB^{(n)}_k w+l^{(n)}_k, \qquad
l^{(n)}_k= \left\{\begin{array}{ll}
0,& k=1,2,3\\[2mm]
\iota_1 {\bf e}_1,&k=4\\[2mm]
\iota_1 {\bf e}_{1}+\iota_2 {\bf e}_{2n},&k=5
\end{array}\right. ,\qquad \iota_1, \iota_2\in\R,\end{equation}
where, for $n=1$, $l^{(1)}_5=\iota_1 {\bf e}_{1}$, or to $w'=\CB^{(n)}_6 w+l^{(n)}_6$, which is a direct product of systems in the above $5$ cases with lower dimensions, i.e., there exists a mutually disjoint partition $\bigcup_{j\in\Lambda}{\CI_j}=\{1,\cdots, n\}$ with $\# \CI_j=m_j<n$ and $k_j\in\{1,\cdots, 5\}$ such that $\CB_{6}^{(n)}=\bigoplus_{j\in\Lambda}  \CB_{k_j}^{(m_j)}$, $l_{6}^{(n)}=\bigoplus_{j\in\Lambda}  l_{k_j}^{(m_j)}$, and (\ref{nonHomoODE}) is conjugated to the direct product of
\begin{equation*}\label{thesecondstep1}
w'_{\CI_j}=  \CB_{k_j}^{(m_j)} w_{\CI_j} + l_{k_j}^{(m_j)}, \qquad w_{\CI_j}:=\begin{pmatrix}
(w_d)_{d\in \CI_j}\\ (w_{n+d})_{d\in \CI_j}
\end{pmatrix}\in \R^{2m_j}.
\end{equation*}
\end{Proposition}

\medskip

\noindent{\it Proof of Theorem \ref{main1}.}
Proposition \ref{ODEreducibility} demonstrates that the non-homogeneous system (\ref{nonHomoODE}) can be reduced to symplectic normal forms in the indecomposable cases. Specifically, the system (\ref{nonHomoODE}) is conjugated to the system (\ref{indecom-constant}) via the transformation $z(t)=U(t)w(t)+v(t)$. Following Proposition \ref{MetaandSch}, the Hamiltonian PDE (\ref{orig-equ-1}) is then conjugated to the quantum normal form
\begin{equation}\label{newSchequ}
\frac{1}{2\pi {\rm i}} \partial_t \varphi=\CR(Z)\varphi = (\CQ_{\CB}(Z) + \CL_{l}(Z))\varphi
\end{equation}
through an $L^2-$unitary transformation, up to a scalar multiplication that does not alter the norm,
\begin{equation}\label{new-unitrans}
\psi(t)=\CU(t) \varphi(t),\qquad \CU(t): = \rho(v(t)) \CM(U(t)).
\end{equation}
The five specific forms of $\CB$ and $l$ in (\ref{indecom-constant}) yield the quantum Hamiltonian $\CR(Z)$, as delineated in QNF1 -- QNF5 of Theorem \ref{main1}.

In the decomposable case of $\CB$, symplectic normal forms of $\#\Lambda$ lower-dimensional independent systems can be deduced using Proposition \ref{ODEreducibility}. According to Proposition \ref{MetaandSch}, Eq. (\ref{orig-equ-1}) is also conjugated to Eq. (\ref{newSchequ}), which effectively becomes the direct sum of $\#\Lambda$ independent quantum Hamiltonians,
$$\CR(Z)=\sum_{j\in\Lambda}\CR_{k_j}^{(m_j)}(Z_{\CI_j})=\sum_{j\in\Lambda}
 \left(\CQ_{\CB_{k_j}^{(m_j)}}(Z_{\CI_j})+  \CL_{l_{k_j}^{(m_j)}}(Z_{\CI_j})\right),\quad  Z_{\CI_j}:=\begin{pmatrix}(D_d)_{d\in\CI_j} \\ (X_d)_{d\in\CI_j}\end{pmatrix}.$$
This is consistent with QNF6 in Theorem \ref{main1}.


\section{Growth rate of Sobolev norms -- Proof of Theorem \ref{main2} \& \ref{thm-ODEPDE}}\label{sec6}

Considering Remark \ref{rmk-SoboNorm-1} and \ref{rmk-SoboNorm-2}, all $L^2-$unitary transformations derived through the Metaplectic or Schr\"odinger representation, as expressed in (\ref{new-unitrans}), preserve the growth rate of Sobolev norms. It is therefore essential to concentrate on the QNFs $\CR_k^{(n)}(Z)$, $k=1,\cdots,6$, and evaluate $\|e^{2\pi{\rm i}t\CR^{(n)}_k}u\|_s$ for $u\in \CH^s(\R^n)\setminus\{0\}$.

For $k=1,2,3$, the QNF is homogeneous, meaning $\CR^{(n)}_k=\CQ_{\CB^{(n)}_k}$. Through straightforward calculations, we find that
for $\CB^{(n)}_1=\CB^{(n)}_1(\lambda)$ with $\lambda>0$,
\begin{equation}\label{etB1}
e^{t\CB^{(n)}_1(\lambda)}=\left(\begin{array}{cccc:cccc}
e^{-t\lambda} & (-t )e^{-t\lambda} &  \dots & \frac{(-t)^{n-1}}{(n-1)!} e^{-t\lambda}   &   &  &   &    \\
&   \ddots & \ddots  & \vdots  &   &    &   &     \\
&  &  e^{-t\lambda}  & (-t )e^{-t\lambda}  &  &  &   &      \\
&   &  &  e^{-t\lambda} &  &  &  &    \\
\hdashline
&  &   &  & e^{t\lambda}  &    &  &  \\
&  &  &  & t e^{t\lambda} &  e^{t\lambda}  &  &   \\
&  &  & & \vdots  &  \ddots & \ddots &      \\
&  &  &  &  \frac{t^{n-1}}{(n-1)!} e^{t\lambda}   &   \dots &   t e^{t\lambda}  &    e^{t\lambda} \end{array}\right),
\end{equation}
for $\CB^{(n)}_2=\CB^{(n)}_2(\lambda_1,\lambda_2)$ with $\lambda_1,\lambda_2>0$,
\begin{equation}\label{etB2}
e^{t\CB^{(n)}_2(\lambda_1,\lambda_2)}=\left(\begin{array}{cccc:cccc}
e^{tA_2} &    & & & &  &  & \\
-te^{tA_2}  & e^{tA_2} & & & & & & \\
\vdots & \ddots & \ddots & & & & & \\
\frac{(-t)^{\frac{n}2-1}}{(\frac{n}2-1)!}e^{tA_2} & \dots & -te^{tA_2} &  e^{tA_2} & & & & \\[2mm]
\hdashline
& & & & e^{-tA_2^*} & te^{-tA_2^*} & \dots & \frac{t^{\frac{n}2-1}}{(\frac{n}2-1)!}e^{-tA_2^*}\\
& & & & & \ddots & \ddots & \vdots \\
& & & & & & e^{-tA_2^*} & te^{-tA_2^*} \\
& & & & & & & e^{-tA_2^*}\end{array}\right),
\end{equation}
where $e^{tA_2}=\begin{pmatrix}
e^{-\lambda_1t}\cos(\lambda_2t) & -e^{-\lambda_1t}\sin(\lambda_2t)\\
e^{-\lambda_1t}\sin(\lambda_2t) & e^{-\lambda_1t}\cos(\lambda_2t)
\end{pmatrix}$, and for $\CB_{3}^{(n)}=\CB_{3}^{(n)}(\mu,\gamma)$, $\gamma=\pm 1$, $\mu>0$,
\begin{equation}\label{etB3}
e^{t\CB_{3}^{(n)}(\mu,\gamma)}=\left(\begin{array}{cccc:cccc}
a_0(t) & a_1(t)  &  \cdots & a_{n-1}(t)   &  b_{n-1}(t)  &  \cdots   &  b_1(t)  & b_0(t)   \\
	 &   \ddots & \ddots  & \vdots  &   \vdots  &   \begin{sideways}$\ddots$\end{sideways} & \begin{sideways}$\ddots$\end{sideways}    &     \\
        &  &  a_0(t)  & a_1(t)  &  b_1(t)  & b_0(t) &   &      \\
	 &   &  &  a_0(t) &  b_0(t)   &  &  &    \\
	\hdashline
	 &  &   &    -b_0(t) & a_0(t)   &    &  &  \\
	 &  & -b_0(t) &    -b_1(t) & a_1(t) &  a_0(t)   &  &   \\
	 &   \begin{sideways}$\ddots$\end{sideways} &   \begin{sideways}$\ddots$\end{sideways} &  \vdots & \vdots  &  \ddots & \ddots &      \\
	   -b_0(t)&    -b_1(t) &   \cdots  &  - b_{n-1}(t) &   a_{n-1}(t)    &   \cdots &  a_1(t)  &   a_0(t)
	\end{array}
\right),
\end{equation}
where $a_j(t):=\frac{(\gamma t)^j}{j!} \cos(\gamma \mu  t+\frac{j \pi}{2})$, $b_j(t):=\frac{(\gamma t)^j}{j!} \sin(\gamma \mu t-\frac{j \pi}{2})$. Hence, as $|t|\to\infty$,
$$\|e^{t\CB_1^{(n)}(\lambda)}\| \simeq  |t|^{n-1} e^{\lambda |t|},\qquad \|e^{t\CB_2^{(n)}(\lambda_1,\lambda_2)}\| \simeq |t|^{\frac{n}{2}-1}e^{\lambda_1 |t|},\qquad \|e^{t\CB_3^{(n)}(\mu,\gamma)}\| \simeq  |t|^{n-1}.$$
Leveraging Corollary \ref{cor_UetA} - (ii), for $s> 0$, the norm $\|e^{2\pi{\rm i}t\CR^{(n)}_k}u\|_s$ aligns with the growth of $\|e^{t\CB_k^{(n)}}\|^s=g^{(n)}_k(s,t)$, $k=1,2,3$, thus substantiating GR1 -- GR3 in Theorem \ref{main2}.

\medskip

\begin{Proposition}\label{prop6.2} Given $u\in \CH^s\setminus\{0\}$, $s>0$, we have $\|e^{2\pi{\rm i}t\CR^{(n)}_k}u\|_s \simeq_u g^{(n)}_k(s,t)$ for $k=4,5$.
\end{Proposition}

\proof For QNF4 $\CR^{(n)}_4=\CQ_{\CB^{(n)}_4}+\CL_{l_4^{(n)}}$, let us focus on $\gamma=1$. A direct computation yields that
\begin{eqnarray}
e^{t\CB^{(n)}_4(1)}&=&\left(\begin{array}{cccc:cccc}
1 & & &  &   &  &  &      \\
-t & 1 & &  &   &  &  &      \\
\vdots & \ddots  & \ddots & &  &  &  &      \\
\frac{(-t)^{n-1}}{(n-1)!} &  \dots & -t & 1 &  &  &  &   \\
\hdashline
\frac{t^{2n-1}}{(2n-1)!} & \cdots &  (-1)^{n}\frac{t^{n+1}}{(n+1)!}  & (-1)^{n+1}\frac{t^{n}}{n!} & 1 & t   &  \dots  & \frac{t^{n-1}}{(n-1)!}   \\
\vdots &  & \vdots  & \vdots &   & \ddots & \ddots & \vdots   \\
\frac{t^{n+1}}{(n+1)!} & \cdots  &   (-1)^{n}\frac{t^3}{3!} & (-1)^{n+1}\frac{t^2}{2!} &  &  & 1 & t \\
\frac{t^{n}}{n!} & \cdots & (-1)^{n}\frac{t^2}{2!}  & (-1)^{n+1}t &   &   &    &  1
\end{array}\right)\nonumber\\
&=:&\begin{pmatrix}
A(t) & 0 \\C(t) & A(t)^{*-1}\end{pmatrix}\in \Sp(n,\R).\label{etB4}
\end{eqnarray}
For $t\in \R$, define $b_4(t)\in \R^{2n}$ and $c_4(t)\in\R$ by
 $$b_4(t)^*=\iota_1\left(t,\cdots,\frac{t^n}{n!},-\frac{(-t)^{2n}}{(2n)!},\cdots,-\frac{(-t)^{n+1}}{(n+1)!}\right)=:(p_4(t)^*,q_4(t)^*),\qquad c_4(t):=-\frac{\iota_1^2}2\cdot\frac{t^{2n+1}}{(2n+1)!}.$$
Then, a straightforward calculation verifies that
\begin{eqnarray}
\left(e^{2\pi{\rm i}t\CR^{(n)}_4}u\right)(x)&=&e^{2\pi{\rm i}c_4(t)}\left(\M\left(e^{t\CB^{(n)}_4(1)}\right)\rho(b_4(t)) \, u\right)(x)\nonumber\\
&=&e^{2\pi{\rm i}c_4(t)}e^{-\pi{\rm i}\la x,C(t)A(t)^{-1}x\ra}e^{2\pi{\rm i}\la q_4(t), A(t)^{-1}x\ra+\pi{\rm i}\la p_4(t),q_4(t)\ra}u(A(t)^{-1}x+p_4(t)). \ \label{solCase4}
\end{eqnarray}
For $1\leq j\leq n$, we have
\begin{eqnarray}\left\|X_j^s e^{2\pi{\rm i}t\CR^{(n)}_4}u\right\|_{L^2}&=&\left(\int_{\R^n}  |x_j|^{2s}|u(A(t)^{-1}x+p_4(t))|^2dx\right)^{\frac12}\nonumber\\
&=&\left(\int_{\R^n} |(A(t)y-A(t)p_4(t))_j|^{2s}|u(y)|^2dy\right)^{\frac12}.\label{XsR4}
 \end{eqnarray}
Then, by noting that $\|A(t)\|\simeq  |t|^{n-1}$ and $(A(t)p_4(t))^*
=\iota_1\left(t,-\frac{t^2}{2}, \cdots,\frac{(-1)^{n+1}}{n!}t^n\right)$, we obtain
\begin{equation}\label{A4es1}
\sum_{j=1}^n\left\|X_j^s e^{2\pi{\rm i}t\CR^{(n)}_4}u\right\|_{L^2}\simeq |t|^{(n-1)s}\|X_1^s u\|_{L^2}+|\iota_1|^s|t|^{ns}\|u\|_{L^2}, \qquad |t|\to \infty.
\end{equation}
On the other hand, for $1\leq j\leq n$ and $s\in\N^*$,
\begin{eqnarray*}
& &\left|\left(D_j^s e^{2\pi{\rm i}t\CR^{(n)}_4}u\right)(x)\right|\\
&=&\left|\frac{1}{(2\pi {\rm i})^s}\partial^s_{x_j} \left(e^{-\pi{\rm i}\la x,C(t)A(t)^{-1}x\ra}e^{2\pi{\rm i}\la q_4(t), A(t)^{-1}x\ra}u(A(t)^{-1}x+p_4(t))\right)\right|\\
&=&\frac{1}{(2\pi)^s}\left|\sum_{k=0}^s C_s^k \partial^{s-k}_{x_j} \left(e^{-2\pi{\rm i}\left\la \frac12 C(t)A(t)^{-1} x-A(t)^{*-1}q_4(t),x\right\ra}\right) \partial^{k}_{x_j} \left(u(A(t)^{-1}x+p_4(t))\right)\right|.
\end{eqnarray*}
By computations, we have, as $|t|\to\infty$, $\|A(t)^{-1}\|\simeq  |t|^{n-1}$, $\|C(t)A(t)^{-1}\|\simeq  |t|^{2n-1}$ and
$$(A(t)^{*-1}q_4(t))_1= -\iota_1t^{2n} \frac{(-1)^{n+1}}{2(n!)^2} ,\qquad |(A(t)^{*-1}q_4(t))_j|\simeq  |\iota_1||t|^{2n+1-j},\quad 2\leq j\leq n,$$
$$(C(t)p_4(t))_1= \iota_1t^{2n} \left(\frac{(-1)^{n+1}}{2(n!)^2} +\frac{1}{(2n)!}\right),\qquad |(C(t)p_4(t))_j|\simeq  |\iota_1||t|^{2n+1-j},\quad 2\leq j\leq n.$$
Hence, as $|t|\to\infty$,
\begin{eqnarray}
\left\|D_j^s e^{2\pi{\rm i}t\CR^{(n)}_4}u\right\|_{L^2}&\lesssim&\left(\int_{\R^n}  |(C(t)A(t)^{-1} x -A(t)^{*-1}q_4(t))_j|^{2s}|u(A(t)^{-1}x+p_4(t))|^2dx\right)^{\frac12}  \nonumber\\
&=&\left(\int_{\R^n}  |(C(t)y- C(t)p_4(t) -A(t)^{*-1}q_4(t))_j|^{2s}|u(y)|^2dy\right)^{\frac12}  \nonumber\\
&\lesssim& |t|^{(2n-2)s \|u\|_{s}} + |\iota_1|^s|t|^{(2n-1)s}\|u\|_{L^2},\quad 2\leq j\leq n, \label{A4es2}\\
\left\|D_1^s e^{2\pi{\rm i}t\CR^{(n)}_4}u\right\|_{L^2}&\simeq&\left(\int_{\R^n}  |(C(t)A(t)^{-1} x -A(t)^{*-1}q_4(t))_1|^{2s}|u(A(t)^{-1}x+p_4(t))|^2dx\right)^{\frac12} \nonumber\\
&=&\left(\int_{\R^n}  |(C(t)y- C(t)p_4(t) -A(t)^{*-1}q_4(t))_1|^{2s}|u(y)|^2dy\right)^{\frac12}  \nonumber\\
&\simeq& (1+|t|^{(2n-1)s})\|X_1^su\|_{L^2} + |\iota_1|^s |t|^{2ns}\left|\frac{(-1)^{n+1}}{(n!)^2} +\frac{1}{(2n)!}\right|^s \|u\|_{L^2}  \nonumber\\
&\simeq& |t|^{(2n-1)s}\|X_1^su\|_{L^2} +   |\iota_1|^s |t|^{2ns}\|u\|_{L^2}. \label{L411}
\end{eqnarray}
Combining (\ref{A4es1}) -- (\ref{L411}) and the interpolation, we have $\|e^{2\pi{\rm i}t\CR^{(n)}_4}u\|_s \simeq_u |t|^{(2n-1)s}+ |\iota_1|^s |t|^{2ns}$ as $|t|\to\infty$ for $s>0$. It is similar for $\gamma=-1$.

\medskip

For QNF5 $\CR^{(n)}_5=\CQ_{\CB^{(n)}_5}+\CL_{l_5^{(n)}}$ with $n$ odd, we have
\beq\label{etB5}
e^{t\CB_5^{(n)}}=\left(\begin{array}{cccc:cccc}
1 & & &  &   &  &  &      \\
-t & 1 & &  &   &  &  &      \\
\vdots & \ddots  & \ddots & &  &  &  &      \\
\frac{(-t)^{n-1}}{(n-1)!} &  \dots & -t & 1 &  &  &  &   \\[1mm]
\hdashline
 &  &   &  & 1 & t   &  \dots  & \frac{t^{n-1}}{(n-1)!}   \\
 &  &  &  &   & \ddots & \ddots & \vdots   \\
 &  &   &  &  &  & 1 & t \\
 &  &   &  &   &   &    &  1\end{array} \right)
 =\begin{pmatrix} A(t) & 0 \\ 0 & A(t)^{*-1} \end{pmatrix}.\eeq
Define $b_5(t)\in \R^{2n}$ for $t\in \R$ by
$$b_5(t)^*=\left(\iota_1t,\cdots,\iota_1\frac{t^{n}}{n!},-\iota_2\frac{(-t)^{n}}{n!},\cdots,-\iota_2(-t)\right)\ := \ (p_5(t)^*,q_5(t)^*), $$
A straightforward calculation shows that
\begin{eqnarray}
\left(e^{2\pi{\rm i}t\CR^{(n)}_5}u\right)(x)&=&\left(\M\left(e^{t\CB_5^{(n)}}\right)\rho(b_5(t))u\right)(x)\nonumber\\
&=&e^{2\pi{\rm i}\la q_5(t), A(t)^{-1}x\ra+\pi{\rm i}\la p_5(t),q_5(t)\ra}u(A(t)^{-1}x+p_5(t)).\label{solCase5}
\end{eqnarray}
By expressing $\left\|X_j^se^{2\pi{\rm i}t\CR^{(n)}_5}u\right\|_{L^2}$ for $1\leq j\leq n$ as (\ref{XsR4}), and noting that, as $|t|\to \infty$,
$$\|A(t)\|\simeq |t|^{n-1},\qquad (A(t)p_5(t))_n=\frac{\iota_1}{n!} t^n,\qquad |(A(t)p_5(t))_j|\lesssim |\iota_1|  |t|^{j},\quad 1\leq j \leq n-1,$$
we obtain
\begin{eqnarray}
\sum_{j=1}^n\left\|X_j^se^{2\pi{\rm i}t\CR^{(n)}_5}u\right\|_{L^2}&=&\sum_{j=1}^n\left(\int_{\R^n}|(A(t)y-A(t)p_5(t))_j|^{2s}|u(y)|^2dy\right)^\frac12\nonumber\\
&\simeq&|t|^{(n-1)s}\|X_1^s u\|_{L^2}+|\iota_1|^s|t|^{ns}\|u\|_{L^2}. \label{A5feiqici1}
\end{eqnarray}
On the other hand, for $1\leq j\leq n$ and $s\in\N^*$,
\begin{eqnarray*}
\left|\left(D_j^s e^{2\pi{\rm i}t\CR^{(n)}_5}u\right)(x)\right|
&=&\left|\frac{1}{(2\pi {\rm i})^s}\partial^s_{x_j} \left(e^{2\pi{\rm i}\la q_5(t), A(t)^{-1}x\ra}u(A(t)^{-1}x+p_5(t))\right)\right|\\
&=&\frac{1}{(2\pi)^s}\left|\sum_{k=0}^s C_s^k \partial^{s-k}_{x_j} \left(e^{2\pi{\rm i}\left\la A(t)^{*-1}q_5(t),x\right\ra}\right) \partial^{k}_{x_j} \left(u(A(t)^{-1}x+p_5(t))\right)\right|,
\end{eqnarray*}
The fact that  $\|A(t)^{-1}\|\simeq  |t|^{n-1}$ as $|t|\to\infty$ implies
$$\left|\partial^{k}_{x_j} \left(u(A(t)^{-1}x+p_5(t))\right)\right|\lesssim   |t|^{(n-1)k} \left|(D_j^{k}u)\left(A(t)^{-1}x+p_5(t)\right)\right|.$$
By computations, we have, as $|t|\to \infty$,
$$(A(t)^{*-1}q_5(t))_1= \iota_2 \frac{t^{n}}{n!} ,\qquad |(A(t)^{*-1}q_5(t))_j|\simeq  |\iota_2||t|^{n+j-1},\quad 2\leq j\leq n. $$
Hence, as $|t|\to\infty$,
\begin{eqnarray}
 \left\|D_j^s e^{2\pi{\rm i}t\CR^{(n)}_5}u\right\|_{L^2}&\lesssim&\sum_{k=0}^s |\iota_2|^{s-k} |t|^{(n+j-1)(s-k)}  |t|^{(n-1)k} \left(\int_{\R^n}  |(D_j^{k}u) (A(t)^{-1}x+p_5(t))|^2dx\right)^{\frac12}  \nonumber\\
&\lesssim& |t|^{(n-1)s}\sum_{k=0}^s |\iota_2|^{s-k} \| D^k_j u\|_{L^2},\quad 2\leq j\leq n, \label{A5es2}\\
\left\|D_1^s e^{2\pi{\rm i}t\CR^{(n)}_5}u\right\|_{L^2}&\simeq&|(A(t)^{*-1}q_5(t))_1|^{s}\left(\int_{\R^n}  |u(A(t)^{-1}x+p_5(t))|^2dx\right)^{\frac12}
\ = \ \frac{ |\iota_2|^s }{n!^s} |t|^{ns}\|u\|_{L^2}. \  \  \label{L511}
\end{eqnarray}
Combining (\ref{A5feiqici1}) -- (\ref{L511}) and the interpolation, we have, for $s>0$,
$$\|e^{2\pi{\rm i}t\CR^{(n)}_5}u\|_s \simeq_u |t|^{(n-1)s}+ (|\iota_1|+|\iota_2|)^s |t|^{ns},\quad |t|\to\infty.\qed$$
	
\begin{remark} The explicit solutions expressed in (\ref{solCase4}) and (\ref{solCase5}) are derived by formally applying the Baker-Campbell-Hausdorff formula. Further details on this process will be provided in Appendix \ref{Sec_BCH}.\end{remark}
				
\medskip

For QNF6, we have the following proposition, which completes the proof of Theorem \ref{main2}.

\begin{Proposition}
For $\CR_6^{(n)}(Z)=\sum_{j\in\Lambda}\CR_{k_j}^{(m_j)}(Z_{\CI_j})$ as (\ref{R6}), for given $u\in \CH^s\setminus\{0\}$, $s>0$,
\begin{equation}\label{L6guji}
\|e^{2\pi{\rm i}t\CR_6^{(n)}(Z)}u\|_{s}\simeq_u \sum_{j\in\Lambda}g^{(m_j)}_{k_j}(s,t).
\end{equation}
\end{Proposition}
\proof Note that, for the mutually disjoint partition $\bigcup_{j\in\Lambda}{\CI_j}=\{1,\cdots, n\}$, $\CR_{k_j}(Z_{\CI_j})$ commutes with $\CR_{k_i}(Z_{\CI_i})$ for any $i\ne j$. For $\CT_l:=2\pi(D_l^2+X_l^2)$, $1\leq l \leq n$, $\CR_{k_j}^{(m_j)}(Z_{\CI_j})$ also commutes with $\CT_l$ for $l\notin \CI_j$. Since $\CR_{k_j}^{(m_j)}(Z_{\CI_j})$ are all self-adjoint, we have
\begin{eqnarray*}
\|e^{2\pi{\rm i}t\CR_6^{(n)}(Z)}u\|_{s}&\simeq&\|u\|_{L^2(\R^n)}+\sum_{l=1}^n\|\CT_l^{\frac s2}e^{2\pi{\rm i}t\CR_6^{(n)}(Z)}u\|_{L^2(\R^n)}\nonumber\\
&=&\|u\|_{L^2(\R^n)}+\sum_{j\in\Lambda}\sum_{l\in\CI_j}\|\CT_l^{\frac s2}e^{2\pi{\rm i}t\CR_{k_j}^{(m_j)}(Z_{\CI_j})}u\|_{L^2(\R^n)}\nonumber\\
&\simeq&\|u\|_{L^2(\R^n)}+\sum_{j\in\Lambda}\|e^{2\pi{\rm i}t\CR_{k_j}^{(m_j)}(Z_{\CI_j})}u\|_{\CH^s(\R^n)}.\label{L61}
\end{eqnarray*}
For $\CI=\CI_j$, $j\in\Lambda$, with $\#\CI_j=m_j<n$ and $k=k_j$, which gives the decomposition
$$\{1,\cdots, n\}=\CI\cup\CI^c,\qquad x=(x_\CI,x_{\CI^c})\in \R^m\bigoplus \R^{n-m}, \quad m=m_j $$
let us show that, for $u\in \CH^s(\R^n)\setminus\{0\}$,
$$\|e^{2\pi{\rm i}t\CR_k^{(m)}(Z_{\CI})}u\|_{\CH^s(\R^n)}\simeq_u g^{(m)}_k(s,t),\qquad |t|\to\infty.$$
In view of Fubini's theorem, we have, for any $u\in \CH^s(\R^n)\setminus \{0\}$, $u(\cdot, x_{\CI^c})\in \CH^s(\R^m)\setminus \{0\}$ for a.e. $x_{\CI^c}\in \R^{n-m}$, $u(x_{\CI},\cdot)\in \CH^s(\R^{n-m})\setminus \{0\}$ for a.e. $x_{\CI}\in \R^{m}$, and
$$\|u\|_{\CH^s(\R^n)}\simeq  \left(\int_{\R^{m}}\|u(x_{\CI}, \cdot)\|^2_{\CH^s(\R^{n-m})}dx_{\CI}\right)^{\frac12} + \left(\int_{\R^{n-m}}\|u(\cdot, x_{\CI^c})\|^2_{\CH^s(\R^m)}dx_{\CI^c}\right)^{\frac12}.  $$
From a straightforward computation we have, as $|t|\to\infty$,
$$
\|e^{2\pi{\rm i}t\CR_k^{(m)}(Z_{\CI})}u\|^2_{s}\simeq\sum\limits_{j\in \CI^c}\|\CT_j^{\frac{s}{2}} u\|^2_{L^2} +  \int_{\R^{n-m}}\|e^{2\pi{\rm i}t\CR_k^{(m)}(Z_{\CI})}u(\cdot, x_{\CI^c})\|^2_{\CH^{s}(\R^{m})} dx_{\CI^c}\nonumber\\
\simeq_u g_k^{(m)}(s,t). \qed$$	

\medskip

\noindent{\it Proof of Theorem \ref{thm-ODEPDE}.} Since the system (\ref{OP-O}) is reducible, there exist $U(\cdot)\in C^1_b(\R,\Sp(n,\R))$ and $\ell(\cdot)\in C^1_b(\R,\R^{2n})$ such that, through the transformation $z=U(t)w+v(t)$, (\ref{OP-O}) is conjugated to
\beq\label{OP-reduO}
w'=\CB w+ l,\qquad \CB=\CB_k^{(n)},\qquad l_k^{(n)}\in \R^{2n}, \qquad k=1,\cdots,6,\eeq
as in Proposition \ref{ODEreducibility}. Let $w^{\rm opt} (t)$ be the solution to (\ref{OP-reduO}) which presents the optimal growth in the sense that, for any solution $w_0(t)$ to (\ref{OP-reduO}),
\beq\label{fast-w*}\sup_{t\in\R}\frac{\|w_0(t)\|}{1+\|w^{\rm opt} (t)\|}<\infty.\eeq
For $k=1,2,3$, we have $w^{\rm opt} (t)=e^{t\CB_k^{(n)}} w^{\rm opt} (0)$. In view of (\ref{etB1}) -- (\ref{etB3}), $w^{\rm opt} (t)$ can be determined as
$$w^{\rm opt} (0)=\left\{\begin{array}{lc}
{\bf e}_n+{\bf e}_{n+1}, & \CB=\CB_1^{(n)}\\[1mm]
{\bf e}_{1}+{\bf e}_{2n},& \CB=\CB_2^{(n)}\\[1mm]
{\bf e}_n+{\bf e}_{n+1}&\CB=\CB_3^{(n)}
\end{array}\right. ,\qquad \|w^{\rm opt} (t)\|\simeq \left\{\begin{array}{lc}
|t|^{n-1}e^{\lambda |t|}, & \CB=\CB_1^{(n)}\\[1mm]
|t|^{\frac{n}{2}-1}e^{\lambda_1 |t|},& \CB=\CB_2^{(n)}\\[1mm]
|t|^{n-1},&\CB=\CB_3^{(n)}
\end{array}\right. ,\quad |t|\to\infty,
$$
which implies that, for $k=1,2,3$,
$$\|w^{\rm opt} (t)\|\simeq\left(g^{(n)}_k(s,t)\right)^{\frac1s},\qquad  |t|\to\infty.$$
For $k=4$, by a direct computation on the system $w'=\CB_4^{(n)} w+\iota_1 {\bf e}_1$, we have that
$$w(t)=e^{t\CB_4^{(n)}} w(0)+\iota_1\left(t, -\frac{t^2}{2},\cdots,(-1)^{n+1}\frac{t^n}{n!},\frac{t^{2n}}{(2n)!},\cdots,\frac{t^{n+1}}{(n+1)!}\right)^*.$$
Then, in view of (\ref{etB4}), $w^{\rm opt} (t)$ can be determined as
$$w^{\rm opt} (0)= {\bf e}_{1}, \qquad \|w^{\rm opt} (t)\|\simeq |t|^{2n-1}+ |\iota_1||t|^{2n}\simeq \left(g^{(n)}_4(s,t)\right)^{\frac1s},\quad \quad |t|\to\infty.$$
For $k=5$, a direct computation on the system $w'=\CB_5^{(n)} w+\iota_1 {\bf e}_1+\iota_{2} {\bf e}_{2n}$ yields that
$$w(t)=e^{t\CB_5^{(n)}} w (0)+\iota_1\left(t, -\frac{t^2}{2}, \cdots,(-1)^{n+1}\frac{t^n}{n!},0,\cdots,0\right)^*+\iota_2\left(0,\cdots,0,  (-1)^{n+1} \frac{t^n}{n!},\cdots, -\frac{t^2}{2}, t\right)^*.$$
Then, in view of (\ref{etB5}), $w^{\rm opt} (t)$ can be determined as
$$w^{\rm opt} (0)= {\bf e}_{1}, \qquad \|w^{\rm opt} (t)\|\simeq |t|^{n-1}+ (|\iota_1|+|\iota_2|)|t|^{n}\simeq \left(g^{(n)}_5(s,t)\right)^{\frac1s},\quad \quad |t|\to\infty.$$
For $k=6$, we have the partition $\bigcup_{j\in\Lambda}{\CI_j}=\{1,\cdots, n\}$ with $\# \CI_j=m_j<n$ and $k_j\in\{1,\cdots, 5\}$. Since the full dimension $n$ is arbitrarily given, the solution $w^{\rm opt}_{\CI_j}$ to the system $w_{\CI_j}'=\CB_{k_j}^{(m_j)} w_{\CI_j}+l_{k_j}^{(m_j)}$ which exhibits the optimal growth in the sense as (\ref{fast-w*}) can be determined as above, and hence $\|w^{\rm opt}_{\CI_j}(t)\|\simeq\left(g_{k_j}^{(m_j)}(s,t)\right)^{\frac1s}$, as $|t|\to\infty$. Then, we have
$$\|w^{\rm opt}\|=\left(\sum_{j\in\Lambda}\|w^{\rm opt}_{\CI_j}\|^2\right)^\frac12\simeq \left(\sum_{j\in\Lambda} g_{k_j}^{(m_j)}(s,t)\right)^{\frac1s}=\left(g^{(n)}_6(s,t)\right)^{\frac1s},\qquad |t|\to\infty.$$
Therefore, according to Theorem \ref{main2}, there exists a constant $C_1>1$, depending on $s$, $\psi(0)$, $\|U\|_{C_b^0}$, $\|v\|_{C_b^0}$ and $w^{\rm opt} (0)$, such that the solution $\psi$ to Eq. (\ref{orig-equ-1}) satisfies
$$C_1^{-1}  \|w^{\rm opt}(t)\|^s \leq\|\psi(t)\|_{s}\leq C_1 \|w^{\rm opt}(t)\|^s, \qquad |t|\rightarrow \infty.$$
Moreover, with the transformation $z=U(t)w+v(t)$, we see that there exists a constant $c_1>1$, depending on $\|U\|_{C^0_b}$ and $\|v\|_{C^0_b}$, such that $c_1^{-1} \|w(t)\| \leq \|z(t)\|\leq c_1 \|w(t)\|$, $\forall \  t\in \R$. Hence, combining with (\ref{fast-w*}), $z^{\rm opt}(t)=U(t)w^{\rm opt}(t)+v(t)$ is a solution to (\ref{OP-O}) which exhibits the optimal growth in the sense of (\ref{fast-z*}), since for any solution $z_0(t)=U(t)w_0(t)+v(t)$ to (\ref{OP-O}),
$$\frac{\|z_0(t)\|}{1+\|z^{\rm opt}(t)\|}\leq    \frac{c_1\|w_0(t)\|}{1+c_1^{-1} \|w^{\rm opt} (t)\|},\qquad \forall \  t\in \R.$$
Then we have $c_1^{-s}C_1^{-1}  \|z^{\rm opt}(t)\|^s \leq\|\psi(t)\|_{s}\leq C_1c_1^s \|z^{\rm opt}(t)\|^s$ as $|t|\rightarrow \infty$, which implies (\ref{estODEPDE}).\qed

\section{Application in time periodic situation}\label{secAppli-per}

\subsection{Quadratic quantum Hamiltonian with time periodic coefficients}

The quadratic quantum Hamiltonian (\ref{orig-equ-1}) with time $\omega-$periodic coefficients, where $\omega>0$, $\CA(\cdot)\in C^0(\T_\omega, \sp(n,\R))$, $\T_\omega:=\R/\omega\Z$, and $\ell(\cdot)\in C^0(\T_\omega, \R^{2n})$, can be interpreted as a periodically forced $n-$D QHO with a time-periodic driving force that may not necessarily be weak. This scenario is commonly encountered in quantum mechanics. The effectiveness of Theorem \ref{main1} -- \ref{thm-ODEPDE} is contingent upon the reducibility of the associated classical $\omega-$periodic system described by $z'=\CA(t)z+\ell(t)$. The reducibility of such periodic linear systems, given by
\begin{equation}\label{linearODE}
z'(t)=\mathcal{A}(t)z(t),\qquad \CA(\cdot)\in C^0(\T_\omega, \sp(n,\R)),
\end{equation}
is affirmed by classical Floquet theory, which provides a foundational framework for analyzing the behavior of solutions to periodic differential equations.

\begin{lemma}\label{SpFloq} (Floquet-Lyapunov Theorem)
There is $\tilde\CB\in \sp(n,\R)$ such that the system (\ref{linearODE}) is conjugated to
$\tilde z'=\tilde\CB\tilde z$ under the transformation $z(t)=S(t)\tilde z(t)$ for some $S(\cdot)\in C^{1}(\T_{2\omega}, \Sp(n,\R))$.
\end{lemma}

\proof Designate $X(t)$ as the fundamental matrix for the linear system (\ref{linearODE}) with the initial condition $X(0)=I_{2n}$. According to Theorem 2.1.3 in \cite{MeyOff2017}, $X(t)$ belongs to $\Sp(n,\R)$ for any $t\in \R$, and it satisfies the property $X(t+\omega)=X(t) X(\omega)$. Consequently, we find that $X(2\omega)=X(\omega)^2$ is also an element of $\Sp(n,\R)$. By invoking Lemma \ref{Dlemma}, we can assert the existence of $\tilde\CB\in \sp(n,\R)$ such that $X(2\omega)=e^{2\omega\tilde\CB}$. We then define $S(t):=X(t)e^{-t\tilde\CB}$ for all $t\in\R$, which yields
$$S'(t)=X'(t)e^{-t\tilde\CB}+X(t)\left(e^{-t\tilde\CB}\right)'
=\mathcal{A}(t)X(t)e^{-t\tilde\CB}-X(t)e^{-t\tilde\CB}\tilde\CB =  \mathcal{A}(t)S(t)-S(t)\tilde\CB. $$
Through direct calculation, it can be shown that the system (\ref{linearODE}) is reducible to $\tilde z'=\tilde\CB\tilde z$ under the symplectic change of variable $z(t)=S(t)\tilde z(t)$. Furthermore, for any $t\in\R$, we have
$$S(t+2\omega)=X(t+2\omega)e^{-(t+2\omega)\tilde\CB} =X(t)X(\omega)^2 e^{-2\omega \tilde\CB}e^{-t \tilde\CB}=S(t),$$
indicating that $S(t)$ is $2\omega$-periodic. Given that $\CA(\cdot)\in C^0(\T_\omega, \sp(n,\R))$, it follows that $X(\cdot)$ belongs to $C^1(\R, \Sp(n,\R))$, thereby ensuring that $S(\cdot)$ is in $C^1(\T_{2\omega}, \Sp(n,\R))$.\qed

\medskip

In conjunction with Lemma \ref{lem-Hormander}, it is evident that, by employing a $2\omega$-periodic transformation $z=U(t)y$, where $U(\cdot) \in C^1(\T_{2\omega}, \Sp(n,\R))$, the system $z'=\CA(t)z+\ell(t)$ can be transformed into
\begin{equation}\label{peri-Btill}
y'=\CB y+\tilde\ell(t),\qquad \tilde\ell(\cdot)\in C^0(\T_{2\omega}, \R^{2n}),
\end{equation}
where $\CB=\CB^{(n)}_k$, for $k=1,\cdots,6$, as specified in Proposition \ref{ODEreducibility}. This transformation effectively reduces the system to one of the symplectic normal forms, facilitating a more straightforward analysis of the system's dynamics and the behavior of its solutions.

To further eliminate the time dependence in $\tilde\ell(t)$, let us recall the following theorem.
\begin{Theorem}(Theorems 1.5.2 in \cite{Ad1995})\label{Thm1.5.2}
Let a linear homogeneous $\omega-$periodic system, $\omega>0$,
\begin{equation}\label{HomoSysAdrianova}
x'=A(t)x,\quad x\in\C^d, \quad  A(t+\omega)=A(t),\quad  A(\cdot)\in C^0(\R,{\rm gl} (d))
\end{equation}
have $k\le d$ linearly independent $\omega-$periodic solutions $\varphi_1(t)$, $\cdots$, $\varphi_k(t)$.
\begin{enumerate}
\item The system $z'=-A^*(t)z$ also has $k$ linearly independent $\omega-$periodic solutions $\psi_1(t)$, $\cdots$, $\psi_k(t)$.
\item The corresponding non-homogeneous system
\begin{equation}\label{nonHomoSysAdrianova}
y'=A(t)y+f(t), \qquad f(t+\omega)=f(t),  \quad   f(\cdot)\in C^0(\R, \R^d),
\end{equation} has $\omega-$periodic solutions if and only if the orthogonality conditions
$$
\int_0^{\omega}\la \psi_j(t),f(t)\ra dt=0,\qquad j=1,\cdots,k,
$$
are satisfied and in this case the $\omega-$periodic solutions form a $k-$parametric family.
\end{enumerate}\end{Theorem}
\begin{remark}\label{REThm1.5.1}
Theorem \ref{Thm1.5.2} also holds for $k=0$. More precisely, if the linear homogeneous $\omega-$periodic system (\ref{HomoSysAdrianova}) has no non-trivial $\omega-$periodic solutions, then any non-homogeneous system (\ref{nonHomoSysAdrianova}) has a unique $\omega-$periodic solution (see Theorems 1.5.1 in \cite{Ad1995}).\end{remark}

According to Theorem \ref{Thm1.5.2} and Remark \ref{REThm1.5.1}, to reduce the system (\ref{peri-Btill}) to constant, it is crucial to determine if there are $2\omega-$periodic solutions for the homogeneous system $y'=\CB y$. Since the full dimension $n$ is arbitrarily given, we focus on the indecomposable cases $k=1,\cdots,5$.
\begin{itemize}
\item [a)] For $\CB=\CB_{k}^{(n)}$, $k=1,2$, or $\CB=\CB_{3}^{(n)}(\mu,\gamma)$ with $\mu\not\in\frac{\pi}{\omega}\N$, in view of (\ref{etB1}) -- (\ref{etB3}), there is no $2\omega-$periodic solution for the homogeneous system $y'=\CB y$. Hence, by Remark \ref{REThm1.5.1}, there exists a unique $2\omega-$periodic solution $v$ to (\ref{peri-Btill}) such that, under the translation $y=w+v(t)$, the system (\ref{peri-Btill}) is conjugated to $w'=\CB w$.

\item [b)] For $\CB=\CB_{4}^{(n)}$, the homogeneous system $y'=-\CB^* y$ admits the $2\omega-$periodic solution $y(t)=y(0)\in \Ker(\CB_{4}^{(n)*})$. According to Theorem \ref{Thm1.5.2}, the system (\ref{peri-Btill}) has $2\omega-$periodic solutions if and only if $$\int_0^{2\omega} \la \tilde\ell(\tau),   \xi \ra d\tau=0,\qquad \forall \  \xi\in \Ker(\CB_{4}^{(n)*}).$$
Note that $\Ker(\CB_{4}^{(n)*})={\rm Span}\{{\bf e}_1\}$. The system
$$y'=\CB y+\tilde\ell(t)-\iota_1 {\bf e}_1,\qquad \iota_1:=\frac{1}{2\omega}\int_0^{2\omega}\la \tilde\ell(\tau),{\bf e}_1\ra d\tau$$
has a $2\omega-$periodic solution $v(t)$. Hence, under the translation $y=w+v(t)$, the system (\ref{peri-Btill}) is conjugated to $w'=\CB w+\iota_1 {\bf e}_1$.

\item [c)] For $\CB=\CB_{5}^{(n)}$, similar to the above case, since $\Ker(\CB_{5}^{(n)*})={\rm Span}\{{\bf e}_1, {\bf e}_{2n}\}$, the system (\ref{peri-Btill}) is conjugated to $w'=\CB w+\iota_1 {\bf e}_1+\iota_2 {\bf e}_{2n}$ with
$$
\iota_1:=\frac{1}{2\omega}\int_0^{2\omega}\la \tilde\ell(\tau),{\bf e}_1\ra d\tau,\qquad  \iota_2:=\frac{1}{2\omega}\int_0^{2\omega}\la \tilde\ell(\tau),{\bf e}_{2n}\ra d\tau,$$
through the translation $y=w+v(t)$ with $v(t)$ a $2\omega-$periodic solution to the system
$y'=\CB y+\tilde\ell(t)-\iota_1 {\bf e}_1-\iota_2 {\bf e}_{2n}$.
\end{itemize}

\begin{remark}\label{RE-resonant}
 Comparing with the above argument a), we notice that $\CB=\CB_{3}^{(n)}(\mu,\gamma)$ with $\mu\in\frac{\pi}{\omega}\N$ does not appear in the classification of symplectic normal forms in this time $\omega-$periodic situation.
Indeed, for $\mu\in\frac{\pi}{\omega}\N$, the system $y'=\CB_{3}^{(n)}(\mu,\gamma) y$ is conjugated to $\tilde y'=\CB_{3}^{(n)}(0,\gamma) \tilde y$
through the time $2\omega-$periodic transformation $y(t)=R_{\mu,\gamma}(t)\tilde y(t)$ with $R_{\mu,\gamma}(t):=\cos(\mu t) I_{2n}+\gamma\sin(\mu t)\J^{\rm a}_{n}$ (recalling (\ref{Janti})).
Since $\CB_{3}^{(n)}(0,\gamma)^n=0$ and, for $n\geq 2$, ${\rm Rank}\left(\CB_{3}^{(n)}(0,\gamma)^{n-1}\right)=2$, we have $\CB_{3}^{(n)}(0,\gamma)\sim \CB_{5}^{(n)}$, which is indeed classified as QNF5 if $n$ is odd, or into the decomposable case if $n$ is even (recalling Remark \ref{rmk-B5even}).
\end{remark}

\begin{Theorem} For the quadratic Hamiltonian PDE (\ref{orig-equ-1}), if $\CA(\cdot)\in C^0(\T_\omega,\sp(n,\R))$ and $\ell(\cdot)\in C^0(\T_\omega,\R^{2n})$ for arbitrary $\omega>0$, then it is reducible to one of QNFs in Theorem \ref{main1}. Moreover, for $\psi(0)\in \CH^s(\R^n)\setminus\{0\}$, $s>0$, the $\CH^s-$norm of the solution $\|\psi(t)\|_s$ exhibits one of GRs in Theorem \ref{main2} corresponding to the reduced QNF.
\end{Theorem}

\subsection{Example of dichotomy of GRs}\label{sec-dicho}

In the realm of time-periodic PDEs, a particular case of interest is the $1-$D Hamiltonian PDE given by:
\begin{equation}\label{qpQHO-omega}
\frac{1}{2\pi{\rm i}}\partial_t \psi=\left(\frac{1}{4\pi }\CT+a\sin(\alpha t)X\right)\psi, \qquad \psi\in \CH^s(\R^1)\setminus\{0\}, \quad s>0,
\end{equation}
where $a$ and $\alpha$ are nonzero real constants. This equation represents a $1-$D QHO subjected to a time-periodic external force characterized by amplitude $a$ and frequency $\alpha$. The growth of the Sobolev norm of the solution to this equation has been thoroughly investigated in the appendix of \cite{BGMR2018}.
\begin{Theorem}\label{thm-BGMR}[Appendix of \cite{BGMR2018}] For the solution $\psi$ to Eq. (\ref{qpQHO-omega}),
\begin{itemize}
\item if $\alpha\neq\pm 1$, then $\sup_t\|\psi(t)\|_s<\infty$,
\item if $\alpha=\pm 1$, then there exists a constant $C>1$ such that $C^{-1}|t|^s<\|\psi(t)\|_s<C|t|^s$ as $|t|\to\infty$.
\end{itemize}
\end{Theorem}

It has been shown in the proof that the above two GRs correspond to QNF3 and QNF5 in $1-$D situation respectively. We give a simplified proof by applying Theorem \ref{thm-ODEPDE}. Since the homogeneous part $\frac{1}{4\pi }\CT=\CQ_{\J_1}=\CQ_{\CB^{(1)}_3(1,1)}$ is constant, Eq. (\ref{qpQHO-omega}) can be seen as Eq. (\ref{orig-equ-1}) in $1-$D situation with $\frac{2\pi}{\alpha}-$periodic coefficients and the corresponding classical system is
$z'=\begin{pmatrix}
0&1\\
-1&0
\end{pmatrix} z+ \begin{pmatrix}
0\\
a\sin(\alpha t)
\end{pmatrix}
$.

If $\alpha\neq\pm1$, then the solution to the above system is
$$z(t)= \begin{pmatrix}
 \cos(t)& \sin(t) \\
 -\sin(t) & \cos(t)
\end{pmatrix} z(0) +\frac{a}{1-\alpha^2}\begin{pmatrix}
\sin(\alpha t)-\alpha\sin(t)\\
\alpha\cos(\alpha t)-\alpha\cos(t)
\end{pmatrix} .$$
Since $\sup_t\|z(t)\|<\infty$ for any $z(0)\in\R^2$, according to Theorem \ref{thm-ODEPDE}, we have that $\sup_t\|\psi(t)\|_s<\infty$.

If $\alpha=\pm1$, then the solution to the above system is
$$z(t)=\begin{pmatrix}
 \cos(t)& \sin(t) \\
 -\sin(t) & \cos(t)
\end{pmatrix} z(0) \pm \frac{a}{2}\begin{pmatrix}
-t\cos( t)+\sin(t)\\
t\sin(t)
\end{pmatrix}. $$
For any $z(0)\in\R^2$, we have $\|z(t)\|\simeq |a| |t|$. Then the $t^s$ growth rate is obtained through Theorem \ref{thm-ODEPDE}.

\subsection{QHO with time periodic quadratic perturbation}

 Consider the time periodic $n-$D QHO described by:
\begin{equation}\label{pQHO-pertu}
\frac{1}{2\pi{\rm i}}\partial_t \psi=\left(\frac{1}{4\pi }\CT-\frac12\left\la Z,  F(t) \J_n Z \right\ra\right)\psi, \qquad \psi\in \CH^s(\R^n)\setminus\{0\}, \quad s>0,
\end{equation}
with $F(\cdot)\in C^0(\T_\omega,\sp(n,\R))$, $\omega>0$. We aim to demonstrate the stability of the $n-$D QHO under a certain non-resonance condition.
\begin{Theorem}\label{thm-pQHO-pertu}
For the solution $\psi(t)$ to Eq. (\ref{pQHO-pertu}),  the following holds.
\begin{itemize}
\item If $\omega\notin\pi\N$ and $\|F\|_{C^0}$ is sufficiently small (depending on $\omega$), then $\sup_t\|\psi(t)\|_s<\infty$.
\item If $\omega\in\pi\N$, then for each GR $g^{(n)}_k(s,t)$, $k=1,\cdots, 5$, as specified in Theorem \ref{main2} with $\iota_1 = \iota_2 = 0$ for $k=4,5$, and for $g^{(n)}_6(s,t) = \sum_{j\in\Lambda}g^{(m_j)}_{k_j}(s,t)$ - the sum of the aforementioned growth rates in lower dimensions with $\sum_{j\in\Lambda} m_j = n$ as described in Theorem \ref{main2} - there exists $F(\cdot)\in C^0(\T_\omega,\sp(n,\R))$ such that $\|\psi(t)\|_s\simeq g^{(n)}_k(s,t)$.
\end{itemize}
\end{Theorem}
\begin{remark} The stability of the $n-$D QHO under time quasi-periodic quadratic perturbations was established in \cite{BGMR2018} under certain non-resonance conditions. Contrary to the periodic case, the quasi-periodic time dependence typically requires higher smoothness assumptions. This aspect will be further discussed in Section \ref{secAppli-qp}.
\end{remark}
The classical system corresponding to Eq. (\ref{pQHO-pertu}) is
\begin{equation}\label{ODE-pertu}
z'=(\J_n+ F(t)) z= (\underbrace{\J_1\bigoplus\cdots \bigoplus\J_1}_{n \ {\rm times}} +F(t)) z.
\end{equation}
The proof of Theorem \ref{thm-pQHO-pertu} relies on the following two claims.

\begin{claim}\label{claim-nonres}
If $\omega\notin\pi\N$, then, with the decomposition
$$\omega=\tilde\omega+2\pi n_\omega,\qquad \tilde\omega\in]-\pi,0[ \, \cup \, ]0,\pi[, \qquad   n_\omega\in\N,$$
there exists $S(\cdot)\in C^1\left(\T_{2\omega},\Sp(n,\R)\right)$ and $\CB\in\sp(n,\R)$ satisfying
$\left\|\CB-\frac{\tilde\omega}{\omega} \J_n\right\|\lesssim_\omega \|F\|_{C^0}=:\epsilon,$
such that, via the transformation $z=S(t)\tilde z$, the linear system (\ref{ODE-pertu}) is conjugated to $\tilde{z}'=\CB \tilde{z}$.
\end{claim}

\proof The proof of Claim \ref{claim-nonres} is indeed a refined one of Lemma \ref{SpFloq} for the concrete time periodic linear system (\ref{ODE-pertu}).
Denote $X(t)$ the fundamental matrix of the system (\ref{ODE-pertu}) with $X(0)=I_{2n}$. By Gr\"onwall's inequality, we have $\|X(t)-e^{t\J_n}\|\lesssim_{\omega} \|F\|_{C^0}=\epsilon$ for $ t\in \left[0, \omega\right]$. Then, by Gershgorin circle Theorem, the eigenvalues $\lambda_j$ of $X(\omega)$ satisfy $|\lambda_j-e^{{\rm i}\omega}|\lesssim_\omega\epsilon$ for $1\le j\le n$ and $|\lambda_j-e^{-{\rm i}\omega}|\lesssim_\omega\epsilon$ for $n+1\le j\le 2n$, suitably reordering the eigenvalues of $X(\omega)$.

If $0<\omega<\pi$, then we have $\tilde\omega=\omega$. According to Lemma \ref{principlelog}, we see that
$\ln (e^{\omega\J_n})=\omega\J_n$ and, by direct computations,
\begin{eqnarray*}
\ln (X(\omega))&=&\ln (e^{\omega\J_n})+\left.\frac{d}{d\tau}\ln \left(e^{\omega\J_n}+\tau(X(\omega)-e^{\omega\J_n})\right)\right|_{\tau=0}\\
& & + \, \int_0^1(1-\tau) \frac{d^2}{d\tau^2}\ln\left(e^{\omega\J_n}+\tau(X(\omega)-e^{\omega\J_n})\right)d\tau\\
&=&{\omega\J_n}+\int_0^1\left(se^{\omega\J_n}+(1-s)I_{2n}\right)^{-1}(X(\omega)-e^{\omega\J_n})\left(se^{\omega\J_n}+(1-s)I_{2n}\right)^{-1}ds+O(\epsilon^2).
\end{eqnarray*}
By Lemma \ref{logofsymp} and Remark \ref{prilog}, $\CB:=\omega^{-1}\ln(X(\omega))\in\sp(n,\R)$, and we have
$$\|\CB-\J_n\|\le\frac{1}{\omega}\sup_{s\in[0,1]}\left\|\left(se^{\omega\J_n}+(1-s)I_{2n}\right)^{-1}\right\|^2\|X(\omega)-e^{\omega\J_n}\|+O(\epsilon^2)\lesssim_{\omega}\epsilon.$$
Let $S(t):=X(t)e^{-t\CB}$, then $S(\cdot)\in C^1(\T_{\omega},\Sp(n,\R))\subset C^1(\T_{2\omega},\Sp(n,\R))$ is the required transformation.

If $\omega> \pi$, then
$\ln(e^{\omega\J_n})=\tilde{\omega}\J_n.$
A similar discussion shows that there is $\CB\in\sp(n,\R)$ with $ \left\|\CB-\frac{\tilde\omega}{\omega}\J_n\right\|\lesssim_\omega \epsilon$, such that $X(\omega )=e^{\omega \CB}$, and $S(t)=X(t)e^{-t\CB}$ is the required.
\qed

\begin{claim}\label{claim-res}
If $\omega\in \pi\N$, then for any $Q\in \sp(n,\R)$, there exist $F(\cdot)\in C^0(\T_{\omega}, \sp(n,\R))$
and $S(\cdot)\in C^1\left(\T_{2\omega},\Sp(n,\R)\right)$ such that, via the transformation $z=S(t)\tilde z$, the system (\ref{ODE-pertu}) is conjugated to $\tilde{z}'=Q \tilde{z}$.
\end{claim}
\proof
For any $Q\in \sp(n,\R)$, the claim is shown by defining $F$ and $S$ as
$$F( t):=\begin{pmatrix}
	\cos(t) \, I_n & \sin(t) \, I_n\\[1mm]
	-\sin(t) \, I_n & \cos(t) \, I_n
\end{pmatrix} \, Q  \,  \begin{pmatrix}
\cos(t) \, I_n & -\sin(t) \, I_n\\[1mm]
\sin(t) \, I_n & \cos(t) \, I_n
\end{pmatrix},\qquad S(t):=\begin{pmatrix}
\cos(t) \, I_n & \sin(t) \, I_n\\
-\sin(t) \, I_n & \cos(t) \, I_n
\end{pmatrix}.\qed$$

\medskip

\noindent{\it Proof of Theorem \ref{thm-pQHO-pertu}.} If $\omega\not\in \pi\N$, then the system (\ref{ODE-pertu}) is conjugated to $\tilde{z}'=\CB \tilde{z}$ with $\|-\J_n\CB-\frac{\tilde\omega}{\omega} I_{2n}\|\simeq\|\CB-\frac{\tilde\omega}{\omega} \J_n\|\lesssim_\omega \epsilon$.
Suppose that $\tilde\omega>0$, then $-\J_n\CB$ has $2n$ positive eigenvalues near $\frac{\tilde\omega}{\omega}$, which implies that the signature of $-\J_n\CB$ is $(2n,0)$
\footnote{For a symmetric $d\times d$ real matrix $A$, its signature is denoted by $(d_+,d_-)$, where $d_+$ (resp. $d_-$) is the number of strictly positive (resp. negative) eigenvalues of $A$.}, and for any $S\in\Sp(n,\R)$, the signature of $-\J_nS^{-1}\CB S$ is still $(2n,0)$.
Hence, $\CB$ is conjugated to the symplectic normal form
$$\CB^{(1)}_3(\mu_1,1)\bigoplus\cdots \bigoplus\CB^{(1)}_3(\mu_n,1)= \bigoplus_{j=1}^n \mu_j\J_1,\qquad \left|\mu_j-\frac{\tilde\omega}{\omega}\right|\lesssim_\omega \epsilon.$$
Similar, if $\tilde\omega<0$, then the signature of $-\J_n\CB$ is $(0, 2n)$ and $\CB$ is conjugated to
$$\CB^{(1)}_3(\mu'_1,1)\bigoplus\cdots \bigoplus\CB^{(1)}_3(\mu'_n,1)= \bigoplus_{j=1}^n \mu'_j\J_1,\qquad \left|\mu'_j+\frac{\tilde\omega}{\omega}\right| \lesssim_\omega \epsilon.$$
Thus, $\sup_t\|\psi(t)\|_s<\infty$ for the solution $\psi(t)$ of Eq. (\ref{pQHO-pertu}) by Theorem \ref{thm-ODEPDE}. \\
\indent If $\omega\in \pi\N$, then the system (\ref{ODE-pertu}) is conjugated to $\tilde{z}'=Q \tilde{z}$.
Since $Q\in \sp(n,\R)$ is arbitrarily chosen, and there is no linear term in the symplectic normal form,
$\|\psi(t)\|_s$ could present every growth rate corresponding to the homogeneous QNF.\qed

\section{Application in time quasi-periodic situation}\label{secAppli-qp}

In the context of time quasi-periodic quadratic Hamiltonian PDEs, the scenario becomes more complex due to the need for higher regularity and the fact that the linear system is not always reducible, even when time analyticity is assumed.

Let $\alpha\in \R^d$ satisfy the Diophantine condition: there exist $\CK>0$ and $\tau>d-1$ such that
\beq\label{Dio}\inf_{j\in\Z}|\la K,\alpha \ra- j| \geq\frac{\CK}{|K|^\tau},\qquad \forall \  K\in\Z^d\setminus \{0\}.
\eeq
A complex number $\beta\in \C$ is called {\it Diophantine w.r.t.} $\alpha$ if there exist $\kappa$, $\sigma\geq 0$ such that
$$
|2\beta-\la K, \alpha\ra|\geq \frac{\kappa}{1+|K|^\sigma},\qquad \forall \  K\in\Z^d.$$
It is easy to verify that $\beta$ is trivially Diophantine w.r.t. $\alpha$ with $\sigma=0$ if $\Im\beta\neq 0$. The real number $\beta\in \R$ is called {\it rational w.r.t.} $\alpha$ if there exists $K_*\in\Z^d$ such that $2\beta=\la K_*, \alpha\ra$.

Consider the time quasi-periodic quadratic Hamiltonian PDE
\begin{equation}\label{nDqp-nonhomo}
\frac{1}{2\pi {\rm i}}\partial_t \psi=\left(\CQ_{\CA(t\alpha)}(Z) +\CL_{\ell(t\alpha)}(Z)\right)\psi,
\end{equation}
with $\alpha\in \R^d$ satisfy the Diophantine condition (\ref{Dio}), and $\CA(\cdot)\in C^r(\T^d,\sp(n,\R))$, $\ell(\cdot)\in C^r(\T^d,\R^{2n})$ with $r\in [0,\infty]$ or $r=\omega$ (analytic class, interpreted as ``$>\infty$").
We assume that the related $\sp(n,\R)-$linear system $z'=\CA(t\alpha)z$ is reducible in the sense that, it can be conjugated to $\tilde z'=\CB\tilde z$, $\CB\in \sp(n,\R)$, through the conjugation $z=S(\frac{t\alpha}2)\tilde z$, where $S(\cdot)\in C^{r+1}(\T^d,\Sp(n,\R))$ (or $\in C^{r}(\T^d,\Sp(n,\R))$) if $r\in [0,\infty[$ (or $r=\infty,\omega$).

\begin{Theorem}\label{thm-qp}
For the time quasi-periodic quadratic Hamiltonian PDE (\ref{nDqp-nonhomo}), assume that the quasi-periodic linear system $z'=\CA(t \alpha) z$ is reducible to $\tilde z'=\CB\tilde z$ in the above sense, and the imaginary part of any eigenvalue of $\CB$ is Diophantine or rational w.r.t. $\alpha$.
If the time-dependence is sufficiently smooth, i.e. $r$ is sufficiently large or $r=\infty,\omega$, then Eq. (\ref{nDqp-nonhomo}) is reducible to one of QNFs in Theorem \ref{main1}. Moreover, for $\psi(0)\in \CH^s(\R^n)\setminus\{0\}$, $s>0$, the $\CH^s-$norm of the solution $\|\psi(t)\|_s$ exhibits one of GRs in Theorem \ref{main2} corresponding to the reduced QNF.
\end{Theorem}

\begin{remark} The theorem referenced provides a framework for understanding the GRs of Sobolev norms in the context of time quasi-periodic quadratic Hamiltonian PDEs such as \cite{BGMR2018} and \cite{LLZ2022}.
In $1-$D scenario, the imaginary part of an eigenvalue of a Hamiltonian matrix $\CB \in \sp(1,\R)$ corresponds to the fibered rotation number associated with the quasi-periodic $\sp(1,\R)$-linear system, as detailed by Eliasson in \cite{Eli1992}. The assumption that it is Diophantine (or rational) w.r.t. $\alpha$ provides the non-resonance (or complete resonance) condition.
The distinction between these arithmetic conditions - Diophantine (implying non-resonance) versus rational (implying complete resonance) - leads to a dichotomy in the behavior of growth rates of Sobolev norms for solutions to the PDE. This dichotomy is highlighted in Theorem \ref{thm-BGMR} and further discussed in Theorem 3.3 and Remark 3.4 of \cite{BGMR2018}.

\end{remark}

\proof To obtain the reducibility of Eq. (\ref{nDqp-nonhomo}), it is sufficient to focus on the affine system $z'=\CA(t\alpha)z+\ell(\alpha t)$, which is transformed into
\beq\label{nonhomo-qp}
\tilde z'=\CB \tilde z + \tilde\ell(t\alpha), \qquad \tilde\ell(\cdot)=S\left(\frac{\cdot}2\right)^{-1} \ell(\cdot)  \in C^r(2\T^d,\R^{2n}),\eeq
through the conjugation $z=S(\frac{t\alpha}2)\tilde z$ in the assumption.
According to Lemma \ref{lem-Hormander}, we can assume that $\CB=\CB^{(n)}_k$, $k=1,\cdots,5$, as in (\ref{A1}) -- (\ref{A5}), or $\CB=\bigoplus_{j\in\Lambda}  \CB_{k_j}^{(m_j)}$ as in (\ref{hormander}).

If the time-dependence is sufficiently smooth, then, for the Fourier expansion
$$\tilde\ell(\theta) =\sum_{K\in \Z^d} \tilde\ell_K e^{\frac{\rm i}2 \la K, \theta\ra},\qquad \tilde\ell_K\in \C^{2n},$$
$|\tilde\ell_K|$ decays sufficiently fast as $|K|\to\infty$. By direct computations, we see that the elimination of time-dependence in $\tilde\ell(t\alpha)$ relies on the invertibility of the matrix $\frac{\rm i}2\la K, \alpha\ra I_{2n}-\CB$, which admits the solvability of equation
\beq\label{eq_coeff}
\left(\frac{\rm i}2\la K, \alpha\ra I_{2n} -\CB\right) v_K=\tilde\ell_K.\eeq

For $\CB=\CB^{(n)}_k$, $k=1,2$, $\frac{\rm i}2\la K,\alpha\ra-\CB$ is invertible for all $K\in \Z^d$, and, for $a_{\pm}:=\frac{\rm i}2\la K,\alpha\ra\pm\lambda$,
$$\left(\frac{\rm i}2\la K,\alpha\ra I_{2n}-\CB^{(n)}_1(\lambda)\right)^{-1} =\left(\begin{array}{cccc:cccc}
a_{+}^{-1} & -a_{+}^{-2}  & \cdots   & (-1)^{n+1}a_{+}^{-n}       &  &  &  &     \\[1mm]
		 & \ddots  &  \ddots  &   \vdots      &  &  &  &    \\[1mm]
	      &   & a_{+}^{-1}   &  -a_{+}^{-2}       &  &  &  &    \\[1mm]
		  &         &   & a_{+}^{-1}  &  &  &  &   \\[1mm]
		\hdashline
		 &  &   &  & a_{-}^{-1}   &      &        &    \\[1mm]
		 &  &   &  &      a_{-}^{-2}    & a_{-}^{-1}  &  &    \\[1mm]
		 &  &   &  &          \vdots   & \ddots  & \ddots &  \\[1mm]
		 &  &   &  &         a_{-}^{-n}   &  \cdots  & a_{-}^{-2}  & a_{-}^{-1}
	\end{array}
	\right),$$
and for $\CC_{-}:=\frac{\rm i}2\la K,\alpha\ra I_{2} - A_2$, $\CC_{+}:=\frac{\rm i}2\la K,\alpha\ra I_{2}+ A_2^*$, with $A_2=\begin{pmatrix}
-\lambda_1&-\lambda_2\\
\lambda_2& -\lambda_1
\end{pmatrix}$,
$$\left(\frac{\rm i}2\la K,\alpha\ra I_{2n}-\CB^{(n)}_2(\lambda_1,\lambda_2)\right)^{-1}=\left(\begin{array}{cccc:cccc}
		\CC_{-}^{-1} &    & & & &  &  & \\
		-\CC_{-}^{-2}  & \CC_{-}^{-1}& & & & & & \\
		\vdots & \ddots & \ddots & & & & & \\
		(-1)^{\frac{n}{2}+1}\CC_{-}^{-\frac{n}{2}} & \dots & -\CC_{-}^{-2}& \CC_{-}^{-1} & & & & \\
		\hdashline
		& & & & \CC_{+}^{-1} &  \CC_{+}^{-2} & \dots & \CC_{+}^{-\frac{n}{2}}\\
		& & & & & \ddots & \ddots & \vdots \\
		& & & & & &  \CC_{+}^{-1} &  \CC_{+}^{-2} \\
		& & & & & & &  \CC_{+}^{-1}
	\end{array}\right),$$
which implies that
$$\left\|\left(\frac{\rm i}2\la K,\alpha\ra I_{2n}-\CB^{(n)}_1(\lambda)\right)^{-1} \right\|\lesssim 1+\lambda^{-n},\qquad \left\|\left(\frac{\rm i}2\la K,\alpha\ra I_{2n}-\CB^{(n)}_2(\lambda_1,\lambda_2)\right)^{-1} \right\|\lesssim 1+\lambda_1^{-\frac{n}{2}}.$$
According to Proposition 3.3.12 in \cite{Gra2014}, if $r\geq d+2$ (including $r=\infty,\omega$), then, with coefficients $v_K$ the solutions to Eq. (\ref{eq_coeff}) for all $K\in \Z^d$, we obtain the translation
\beq\label{translation-qp}
\tilde z=w+v(t \alpha),\qquad v(\cdot)=\sum_{K\in \Z^d} \left(\frac{\rm i}2\la K,\alpha\ra I_{2n}-\CB\right)^{-1}\tilde\ell_K e^{\frac{\rm i}2  \la K, \cdot\ra}\in C^1(2\T^d, \R^{2n}),\eeq
which conjugates the system (\ref{nonhomo-qp}) to $w'=\CB w$, $\CB=\CB^{(n)}_k$, $k=1,2$.

For $\CB=\CB^{(n)}_3(\mu,\gamma)$ with $\mu>0$ Diophantine w.r.t. $\alpha$, i.e., there exist $\kappa,\sigma>0$ such that
$$|\la K, \alpha\ra\pm 2\mu |\geq \frac{\kappa}{1+|K|^{\sigma}},\qquad \forall \ K\in \Z^d,$$
we have that
$\frac{\rm i}2\la K,\alpha\ra I_{2n} -\CB$ is invertible for every $K\in \Z^d$.
Let $S_{\mu,\gamma}$ be the invertible $2n\times 2n$ matrix which realizes the Jordan decomposition of $\CB^{(n)}_3(\mu,\gamma)$, the we have the Jordan decomposition
\beq\label{Jordan-decomp}
\frac{\rm i}2 \la K,\alpha\ra I_{2n}-\CB^{(n)}_3(\mu,\gamma)=S_{\mu,\gamma} \left(\begin{array}{cccc:cccc}
a_+ & -1 &   &   &   &   &   &   \\
 &   \ddots & \ddots  &  &     &   &   &     \\
 &  &  a_+   & -1& &  &   &      \\
 &   &  & a_+ & &  &  &    \\
\hdashline
 &  &  & & a_{-}  & -1   &  &  \\
 &  &  & &  & a_{-}   &  \ddots &   \\
 &  &  & &  &  & \ddots & -1      \\
 &  &  & &  &  &    &   a_{-}
	\end{array}
\right)S_{\mu,\gamma}^{-1},\eeq
with $a_\pm:=\frac{\rm i}2(\la K,\alpha\ra\pm2\mu)$, which shows that, for all $K\in \Z^d$,
\beq\label{esti-Diophantine}
\left\|\left(\frac{\rm i}2 \la K,\alpha\ra I_{2n}-\CB^{(n)}_3(\mu,\gamma)\right)^{-1}\right\|\lesssim_\mu 1+ |\la K,\alpha\ra-2\mu|^{-n} +  |\la K,\alpha\ra+2\mu|^{-n}\lesssim_\mu \kappa^{-n} (1+|K|)^{n\sigma}.\eeq
If $r\geq n\sigma+d+2$ (including $r=\infty,\omega$), then, with $v_K$ the solution to Eq. (\ref{eq_coeff}) for any $K\in \Z^d$, satisfying
$$|v_K|\lesssim_\mu \frac{|\tilde\ell |_{C^r(2\T^d, \R^{2n})}}{(1+|K|)^{r-\sigma n}}, $$
we obtain $v(\cdot)=\sum_{K\in \Z^d} v_K e^{\frac{\rm i}2 \la K, \cdot\ra}\in C^1(2\T^d, \R^{2n})$.
Thus by the translation as (\ref{translation-qp}) the system (\ref{nonhomo-qp}) is conjugated to $w'=\CB^{(n)}_3(\mu,\gamma) w$.

For $\CB=\CB^{(n)}_4(\gamma)$ with $\gamma=1$, $\frac{\rm i}2\la K,\alpha\ra I_{2n} -\CB$ is invertible for $K\in \Z^d\setminus\{0\}$,
and, for $a:=\frac{\rm i}2 \la K,\alpha\ra$,
\begin{eqnarray*}
& &\left(\frac{\rm i}2\la K,\alpha\ra I_{2n}-\CB^{(n)}_4(1)\right)^{-1}\\&=&\left(\begin{array}{cccc:cccc}
	a^{-1} & & &  &   &  &  &      \\
	-a^{-2} & a^{-1} & &  &   &  &  &      \\
	\vdots & \ddots  & \ddots & &  &  &  &      \\
	(-1)^{n+1}a^{-n} &  \dots & -a^{-2} & a^{-1} &  &  &  &   \\
	\hdashline
	a^{-2n}& -a^{-(2n-1)} & \cdots   & (-1)^{n+1}a^{-(n+1)}  & a^{-1} & a^{-2}   &  \dots  & a^{-n}   \\
	a^{-(2n-1)}&  -a^{-(2n-2)}  &  \cdots     &  (-1)^{n+1}a^{-n}  &   & \ddots & \ddots & \vdots   \\
	\vdots  & \vdots   &   & \vdots   &  &  & a^{-1} & a^{-2} \\
	 a^{-(n+1)}&  -a^{-n} &  \cdots      & (-1)^{n+1}a^{-2} &   &   &    &  a^{-1}
\end{array}
\right)
\end{eqnarray*}
which implies that
$$\left\|\left(\frac{\rm i}2\la K,\alpha\ra I_{2n}-\CB^{(n)}_4(1)\right)^{-1}\right\|\lesssim_{\CK} (1+|K|)^{2n\tau}, \qquad K\in \Z^d\setminus\{0\}.$$
If $r\geq 2n\tau+d+2$, then, with $v_K$ the solution to Eq. (\ref{eq_coeff}) for $K\in \Z^d\setminus\{0\}$, satisfying
$$|v_K|\lesssim \frac{|\tilde\ell |_{C^r(2\T^d, \R^{2n})}}{(1+|K|)^{r- 2n\tau}}, $$
we obtain the translation
$$\tilde z=y+ v(t\alpha),\qquad v(\cdot)=\sum_{K\in \Z^d\setminus\{0\}} v_K e^{\frac{\rm i}2 \la K, \cdot\ra}\in C^1(2\T^d, \R^{2n}),$$
such that the system (\ref{nonhomo-qp}) is conjugated to
$$\tilde y'=\CB^{(n)}_4(1) \tilde y+\tilde\ell_{0} ,\qquad  \tilde\ell_{0}=\frac{1}{(4\pi)^n}\int_{2\T^d}\tilde\ell(\theta) \,  d\theta, $$
Then, as shown in Section \ref{sec5}, the system (\ref{nonhomo-qp}) is finally conjugated to
$$w'=\CB^{(n)}_4(1)w+\iota_1 {\bf e}_1 ,\qquad  \iota_1=\la \tilde\ell_{0},  {\bf e}_1\ra.$$
The proof is similar for $\gamma=-1$.

For $\CB=\CB^{(n)}_5$, $\frac{\rm i}2\la K,\alpha\ra I_{2n} -\CB$ is invertible for $K\in \Z^d\setminus\{0\}$,
and
 $$\left(\frac{\rm i}2\la K,\alpha\ra I_{2n} -\CB^{(n)}_5\right)^{-1}=\left(\begin{array}{cccc:cccc}
	a^{-1} & & &  &   &  &  &      \\
	-a^{-2} & a^{-1} & &  &   &  &  &      \\
	\vdots & \ddots  & \ddots & &  &  &  &      \\
	(-1)^{n+1}a^{-n} &  \dots & -a^{-2} & a^{-1} &  &  &  &   \\
	\hdashline
	&  &   &  & a^{-1} & a^{-2}   &  \dots  & a^{-n}   \\
	&  &  &  &   & \ddots & \ddots & \vdots   \\
	&  &   &  &  &  & a^{-1} & a^{-2} \\
	 &  &   &  &   &   &    &  a^{-1}
\end{array}
\right) $$
which implies that
$\left\|\left(\frac{\rm i}2\la K,\alpha\ra I_{2n} -\CB^{(n)}_5\right)^{-1} \right\|\lesssim_{\CK} (1+|K|)^{n\tau}$.
Similar to the above case, if $r>n\tau+d+2$, then
the system (\ref{nonhomo-qp}) is finally conjugated to
$$w'=\CB^{(n)}_5 w+\iota_1 {\bf e}_1+\iota_2 {\bf e}_{2n} ,\qquad  \iota_1=\la \tilde\ell_{0},  {\bf e}_1\ra ,\qquad  \iota_2=\la \tilde\ell_{0},  {\bf e}_{2n}\ra.$$

For $\CB=\CB^{(n)}_3(\mu,\gamma)$ with $2\mu=\la K_*, \alpha\ra$, $K_*\in\Z^d\setminus \{0\}$, according to (\ref{Jordan-decomp}) and the fact that
$$|2\mu\pm\la K, \alpha\ra|= |\la K\pm K_*, \alpha\ra|\geq\frac{\CK}{|K \pm K_*|^\tau},\qquad K\in \Z^d\setminus\{\pm K_*\},$$
$\frac{\rm i}2\la K,\alpha\ra I_{2n}  -\CB^{(n)}_3(\mu, \gamma)$ is invertible for $K\in \Z^d\setminus\{\pm K_*\}$, and
$$
\left\|\left(\frac{\rm i}2\la K,\alpha\ra I_{2n} -\CB^{(n)}_3(\mu,\gamma)\right)^{-1}\right\|\lesssim 1+|\la K+ K_*, \alpha\ra|^{-n} +  |\la K- K_*, \alpha\ra|^{-n}\lesssim_{ \CK, K_*} (1+|K|)^{n\tau}.$$
If $r>n\tau+d+2$, then, with $v_K$ the solutions to Eq. (\ref{eq_coeff}) for $K\in \Z^d\setminus\{\pm K_*\}$, we obtain the translation
$$\tilde z=y+ v(t\alpha),\qquad v(\cdot)=\sum_{K\in \Z^d\setminus\{\pm K_*\}} \left(\frac{\rm i}2\la K,\alpha\ra I_{2n}-\CB_3^{(n)}(\mu, \gamma)\right)^{-1}\tilde\ell_K e^{\frac{{\rm i}}{2}  \la K, \cdot\ra}\in C^1(2\T^d, \R^{2n}),$$
such that the system (\ref{nonhomo-qp}) is conjugated to
$$y'=\CB^{(n)}_3(\mu,\gamma) y+f_{K_*}(t \alpha) ,\qquad  f_{K_*}(\cdot):=\tilde\ell_{K_*} e^{\frac{\rm i}2 \la K_*,\cdot\ra}+\tilde\ell_{-K_*} e^{-\frac{\rm i}2 \la K_*,\cdot\ra}.$$
Through the transformation $y= R_{\mu,\gamma}(t) \tilde y$
$$R_{\mu,\gamma}(t)=\cos\left(\frac{\la K_*,\alpha\ra}{2}  t\right) I_{2n}+\gamma\sin\left(\frac{\la K_*,\alpha\ra}{2}  t\right)\J^{\rm a}_{n},$$
we obtain $\tilde y'=\CB^{(n)}_3(0,\gamma)\tilde y+R_{\mu,\gamma}(t)^{-1}f_{K_*}(t \alpha)$.
As noted in Remark \ref{RE-resonant}, we can achieve reducibility to QNF5 for odd $n$ or to two independent QNF4's for even $n$.

Given that this reducibility has been established for any spatial dimension $n$, it applies to decomposable cases as well, thus concluding the proof of Theorem \ref{thm-qp}. \qed

\begin{remark} Similar to the time-periodic situation (see Remark \ref{RE-resonant}), $\CB=\CB^{(n)}_3(\mu,\gamma)$ with $2\mu=\langle K, \alpha\rangle$, $K\in\Z^d$, considered as the complete resonance case, is also excluded from the classification of QNFs in the time quasi-periodic context.\end{remark}

\appendix					
			
\section{QNFs and GRs of Sobolev norm in $1-$D and $2-$D situation}\label{app_2d}

For Eq. (\ref{orig-equ-1}) with $n=1$, the QNFs and GRs were given in \cite{LLZ2022}. Recall the QNF$k$ and GR$k$, $k=1,\cdots,6$, in Theorem \ref{main1} and \ref{main2}. Except $k=2$ and $k=6$ (which do not appear in the $1-$D situation), we have

\begin{itemize}
\item [($k=1$)] $\CR^{(1)}_{1}(Z)=\CQ_{\CB_1^{(1)}}(Z)$ with
$\CB_1^{(1)}=\begin{pmatrix}
-\lambda & 0 \\
0 & \lambda
\end{pmatrix}$, $\lambda>0$, and $g^{(1)}_1(s,t)=e^{\lambda s|t|}$.

\item [($k=3$)] $\CR^{(1)}_{3}(Z)=\CQ_{\CB_3^{(1)}}(Z)$ with
$\CB_3^{(1)}=\pm \begin{pmatrix}
0 & \mu \\
-\mu & 0
\end{pmatrix}$, $\mu>0$ Diophantine w.r.t. $\alpha$,
and $g^{(1)}_3(s,t)=1$.

\item [($k=4$)] $\CR^{(1)}_{4}(Z)=\CQ_{\CB_4^{(1)}}(Z)+\CL_{l^{(1)}_4}(Z)$ with
$\CB_4^{(1)}=\pm\begin{pmatrix}
0 & 0 \\
1 & 0
\end{pmatrix}$, $l_4^{(1)}=\iota_1 {\bf e}_1$, $\iota_1\in\R$, and $$g^{(1)}_4(s,t)=|t|^{s} +  |\iota_1|^s |t|^{2 s}.$$

\item [($k=5$)] $\CR^{(1)}_{5}(Z)=\CQ_{\CB_5^{(1)}}(Z)+\CL_{l^{(1)}_5}(Z)$ with
$\CB_5^{(1)}=\begin{pmatrix}
0 & 0 \\
0 & 0
\end{pmatrix}$, $l_5^{(1)}=\iota_1 {\bf e}_1$, $\iota_1\in\R$, and $$g^{(1)}_5(s,t)=1  +  |\iota_1|^s  |t|^{s}.$$
\end{itemize}


For Eq. (\ref{orig-equ-1}) with $n=2$, there is no QNF5 and GR5 defined in Theorem \ref{main1} and \ref{main2}.

%

\begin{itemize}
\item [($k=1$)] $\CR_1^{(2)}(Z)=\CQ_{\CB_1^{(2)}(\lambda )}(Z)$ with $\CB_1^{(2)}(\lambda )=\begin{pmatrix}
-\lambda & -1 & 0 & 0\\
0 & -\lambda & 0 & 0\\
0 & 0 & \lambda & 0\\
0 & 0 & 1 & \lambda
\end{pmatrix}$, $\lambda>0$,
and $g_1^{(2)}(s,t)=|t|^{s}e^{\lambda s|t|}$.

\item [($k=2$)] $\CR_2^{(2)}(Z)=\CQ_{\CB_2^{(2)}(\lambda_1,\lambda_2)}(Z)$ with
$\CB_2^{(2)}(\lambda_1,\lambda_2)=\begin{pmatrix}
				-\lambda_1 & -\lambda_2 & 0 & 0\\
				\lambda_2 & -\lambda_1 & 0 & 0\\
				0 & 0 & \lambda_1 & -\lambda_2\\
				0 & 0 & \lambda_2 & \lambda_1
			\end{pmatrix}$, $\lambda_1,\lambda_2>0$,
and $g_2^{(2)}(s,t)=e^{\lambda_1 s|t|}$.

\item [($k=3$)] $\CR_3^{(2)}(Z)=\CQ_{\CB_3^{(2)}(\mu,\gamma)}(Z)$
with
$\CB_3^{(2)}(\mu,\gamma)=\gamma\begin{pmatrix}
				0 & 0 & -1 & \mu\\
				0 & 0 & \mu & 0\\
				0 & -\mu & 0 & 0\\
				-\mu & 1 & 0 & 0
			\end{pmatrix}$, $\mu> 0$, $\gamma=\pm1$,
and $g_3^{(2)}(s,t)=|t|^{s}$.

\item [($k=4$)] $\CR_4^{(2)}(Z)=\CQ_{\CB_4^{(2)}(\gamma)}(Z)+\CL_{l_4^{(2)}}$ with
$\CB_4^{(2)}=\gamma\begin{pmatrix}
				0 & 0 & 0 & 0\\
				-1 & 0 & 0 & 0\\
				0 & 0 & 0 & 1\\
				0 & -1 & 0 & 0
			\end{pmatrix}$, $\gamma=\pm1$, $l_4^{(2)}=\iota {\bf e}_1$, $\iota\in\R$, and $g_4^{(2)}(s,t)= |t|^{3s} +  |\iota|^s |t|^{4 s}$.	
\item [($k=6$)] $\CR(Z)$ is decomposed into two independent operators on $L^2(\R)$, i.e.,
$$\CR_6^{(2)}(Z)=\CR^{(1)}_{k_1}(Z_1)+\CR^{(1)}_{k_2}(Z_2),\quad g_6^{(2)}(s,t)= g^{(1)}_{k_1}(s,t) +g^{(1)}_{k_2}(s,t),\quad k_1, k_2=1,3,4,5. $$
%
%
\end{itemize}
\begin{remark}
The $2-$D normal forms for the homogeneous part were given by Matsumoto-Ueki \cite{MU1998}.
\end{remark}

\section{Proofs of Proposition \ref{boundofschres} and \ref{daoshuforrho}}\label{app_Schro}

The proofs for these two propositions are quite similar. We will provide a detailed proof for Proposition \ref{boundofschres} and an outline for Proposition \ref{daoshuforrho}.
					
\noindent{\it Proof of Proposition \ref{boundofschres}.}
For $w=\begin{pmatrix}p \\ q \end{pmatrix}$ with $ p, q\in \R^{n}$. Then $\|p\|, \|q\|\leq \|w\|$.
For $u\in \CH^s$, $s>0$, we have, through (\ref{SchroRep-L2}), that
$\left(\rho(w)u\right)(x) = e^{2\pi {\rm i}\la q, x\ra +\pi {\rm i}\la p, q\ra}u(x+p)$, and hence
\begin{equation}\label{rhoxingzhi1}
\left\||X|^s \rho(w) u\right\|_{L^2}=\left(\int_{\R^n} |x|^{2s} \left|u(x+p)\right|^2 dx\right)^\frac12\lesssim \|p\|^s \|u\|_{L^2}+ \| |X|^s u\|_{L^2} .
\end{equation}
By a straightforward computation we have $\left(\CF\rho(w)u\right)(\xi)= e^{2\pi {\rm i} \la p,\xi\ra-\pi {\rm i} \la p, q\ra}\hat{u}(\xi-q)$, which implies
\begin{equation}\label{rhodaoshu1}
\||X|^s\CF\rho(w)u\|_{L^2}\lesssim \|q\|^s \|u\|_{L^2}+\| |X|^s \hat{u}\|_{L^2}.
\end{equation}
It follows from (\ref{rhoxingzhi1}) and (\ref{rhodaoshu1}) that
$\|\rho(w) u\|_{s}\lesssim \|w\|^s  \|u\|_{s}$ for $s>0$.
For the other direction note
$
\|u\|_s=\|\rho(-w(t)) \rho(w(t)) u\|_s  \lesssim\|w\|^s \|\rho(w(t)) u\|_s,$
it follows $\|\rho(w(t)) u\|_s  \gtrsim \|w\|^{-s} \|u\|_s$.

For $s< 0$,  from (\ref{dfforrho}) we have
$$\|\rho(w) u\|_{s} =  \sup_{\|v\|_{-s}\leq 1} |\la \rho(w) u, v\ra|
\leq\sup_{\|v\|_{-s}\leq 1}  \|u\|_{s}\|\rho(-w) v\|_{-s}
\lesssim \sup_{\|v\|_{-s}\leq 1}  \|u\|_{s}\cdot \|w\|^{-s} \|v\|_{-s}   \leq  \|w\|^{-s} \|u\|_{s}. $$
On the other hand, note $\rho(-w) \rho(w) =\Id$ on $\CH^{-s}$ for $-s> 0$,
$$\la \rho(-w)\rho(w)u, v\ra = \la u, \rho(-w) \rho(w) v\ra=\la u, v\ra, \qquad u\in \CH^{s},\quad v\in \CH^{-s}.$$
It follows that $\rho(-w)\rho(w)u=u$ for $u\in \CH^{s}$ with $s< 0$.
Therefore, for $u \in \CH^{s}$ with $s< 0$, $$\|u \|_{s}= \|\rho(-w)\rho(w) u \|_{s}\lesssim \|w\|^{-s}  \|\rho(w) u \|_{s},$$
which implies that $\|\rho(w) u \|_{s}\gtrsim \|w\|^{s}\|u \|_{s}$.\qed
					
					\medskip
					
{\it \noindent Proof of Proposition \ref{daoshuforrho}.}
Let $w(\cdot)=\begin{pmatrix}p(\cdot)\\ q(\cdot)\end{pmatrix}$ with $ p(\cdot), q(\cdot)\in C_b^1(\mathbb R, \R^{n})$.					
					
					For $s\in\N^*$, (\ref{rhode}) implies that, for $u\in\CH^s$,
					\begin{eqnarray*}														\partial_t\left(\rho(w(t)) u\right)(x)&=& \pi {\rm i}e^{2\pi {\rm i} \la q(t), x\ra+\pi {\rm i} \la p(t), q(t)\ra}(2\la q'(t), x\ra+   \la p'(t), q(t)\ra+  \la p(t), q' (t)\ra)u(x+p(t))\\
						& & +  \, 2\pi {\rm i} e^{2\pi {\rm i} \la q(t), x\ra+\pi {\rm i} \la p(t), q(t)\ra} (\la p'(t), D \ra u)(x+p(t)),  \label{daoshuforrho1}
\end{eqnarray*}
and hence, $\|\partial_t\left(\rho(w(t))u\right)\|_{s-1}\lesssim \| w\|^{s+1}_{C_b^1} \|u\|_s$, which follows $\partial_t \rho(w(t)) \in   B(\H^{s},\H^{s-1})$ for $s\in \N^*$. By interpolation, we have $\partial_t \rho(w(t)) \in   B(\H^{s},\H^{s-1})$ for $s\geq 1$.
					
For $s\in -\N^*$, we can prove that for $w(\cdot)\in C_b^1(\mathbb R, \R^{2n})$,
\begin{equation}\label{fulu4daoshu}
\la \partial_t\rho(w(t))u, v\ra= \la u, \partial_t\rho(-w(t)) v\ra,\qquad u\in \CH^{s},\qquad v\in \CH^{-s+1}.
\end{equation}
From (\ref{fulu4daoshu}), the conclusion for the case $s\in \N^*$,  and a similar discussion as the proof of Proposition \ref{boundofschres}, we can show that $\partial_t \rho(w(t)) \in   B(\H^{s},\H^{s-1})$ when $s\in -\N^*$. By interpolation, we have $\partial_t \rho(w(t)) \in   B(\H^{s},\H^{s-1})$ for $s\leq -1$. When $|s|<1$ the proof is similar. \qed

\section{Heuristic discussion from Baker-Campbell-Hausdorff formula}\label{Sec_BCH}

We will provide a heuristic proof for the explicit expressions of the solution in (\ref{solCase4}) and (\ref{solCase5}). For $\CR^{(n)}_k = \CQ_{\CB^{(n)}_k} + \CL_{l_k^{(n)}}$, $k=4,5$, the challenge in deriving the expression for the propagator $e^{2\pi{\rm i}t\CR^{(n)}_k}$ comes from the non-commutativity between the homogeneous part $\CQ_{\CB^{(n)}_k}$ and the linear part $\CL_{l_k^{(n)}}$. This non-commutativity prevents the direct application of the Metaplectic and Schr\"odinger representations separately. The primary approach is to express the propagator as $\exp\{2\pi{\rm i}t\CR^{(n)}_k\} = \exp\{2\pi{\rm i}t\CQ_{\CB^{(n)}_k}\} \exp\{2\pi{\rm i}t \CL_{\tilde l_k}\}$ with an appropriately chosen $\tilde l_k\in \R^{2n}$. This $\tilde l_k$ is determined using the Baker-Campbell-Hausdorff (B-C-H) formula.

Let $\CG$ be a Lie group with its associated Lie algebra $\mathfrak{g}$. Consider the exponential map $\exp : \mathfrak{g} \rightarrow \CG$. For elements $Q, L \in \mathfrak{g}$, we denote $G=\ln (\exp\{Q\} \exp\{L\})$, where the logarithmic function is defined such that
$\exp\{Q\} \exp\{L\}=\exp\{G\}$. This can also be expressed using the formal series $\ln(1 + z) := \sum_{k=1}^\infty \frac{(-1)^{k+1}}{k} z^k$.

\begin{Proposition}(B-C-H formula \cite{G56}) For $Q, L \in \mathfrak{g}$, the formal series $G=\ln  (\exp\{Q\} \exp\{L\})$ is a Lie series generated by $\{Q,L\}$. This means that $G$ can be expressed as a series $Q + L + \cdots$, where all terms of order $\geq 2$ can be represented as linear combinations of iterated commutators $[\cdot, \cdot]$, known as Lie elements. More precisely, $G$ has the explicit series representation:
\beq\label{formu-BCH}
G=\sum_{m\ge1}\sum_{p_i+q_i>0}\frac{(-1)^{m-1}}{m}\frac{[\overbrace{Q,\cdots,Q}^{p_1 \ {\rm times}},\overbrace{L,\cdots,L}^{q_1 \ {\rm times}},\cdots,\overbrace{Q,\cdots,Q}^{p_m \ {\rm times}},\overbrace{L,\cdots,L}^{q_m \ {\rm times}}]}{p_1!q_1!\cdots p_m! q_m! \sum_i (p_i+q_i)},\eeq
where $[X, Y] = \ad_X Y := XY - YX$ denotes the commutator of $X$ and $Y$, and $[X_1, X_2, \cdots, X_n] := [X_1, [X_2, \cdots, [X_{n-1}, X_n] \cdots]]$ following the right-normed convention.
\end{Proposition}

For the two non-commuting operators $Q=2\pi{\rm i} t \CQ_{\CA}(Z)$ and $L=2\pi{\rm i} \CL_{\ell_0(t)}(Z)$, where $\CA \in \sp(n, \R)$ is nilpotent and $\ell_0(t) \in \R^{2n}$ is to be determined, we observe that $\exp\{Q\}=\M(e^{t\CA})$ and $\exp\{L\}=\rho(\ell_0(t))$
based on the identities in (\ref{zizhisolutions}) and (\ref{SchroRep-infinitesimal}).
We aim to demonstrate that the formal series (\ref{formu-BCH}) with $Q$ and $L$ specified above actually terminates after a finite number of terms. The finiteness of this series has also been established in Lemma 7 of \cite{Ber2021}.

It is evident that all Lie elements resulting from $Q$ and $L$ are pseudo-differential operators of order less than $2$, given that $Q$ is quadratic and $L$ is linear with respect to the operators $(D, X)$.

\begin{lemma}\label{commutatorcompute}
For any $\CA\in\sp(n,\R)$ and $\ell_0\in\R^{2n}$, we have $\ad^m_{\CQ_\CA(Z)}(\CL_{\ell_0}(Z))=(2\pi{\rm i})^{-m}\CL_{\CA^m\ell_0}(Z)$.\end{lemma}
\proof By a straightforward computation, we have, for $1\leq j,k,l \leq n$,
\begin{eqnarray*}
& &[D_j D_k, X_l]=\frac{1}{2\pi{\rm i}}(\delta_{kl}D_j+\delta_{jl}D_k), \qquad [X_j D_k, X_l]=[D_k X_j, X_l]=\frac{1}{2\pi{\rm i}}\delta_{kl} X_j,\\
& &[X_j D_k, D_l]=[D_k X_j, D_l]=-\frac{1}{2\pi{\rm i}}\delta_{jl}D_k,\qquad [X_j X_k, D_l]=-\frac{1}{2\pi{\rm i}}(\delta_{kl}X_j+\delta_{jl}X_k),
\end{eqnarray*}
which follows that, for $\begin{pmatrix}
A & B \\
C & F
\end{pmatrix}\in{\rm gl}(2n,\R)$ and $\begin{pmatrix}a \\ b \end{pmatrix}\in\R^{2n}$,
$$
\left[\left\la Z, \begin{pmatrix}
A & B \\
C & F
\end{pmatrix}Z\right\ra,\left\la \begin{pmatrix}
a \\ b\end{pmatrix}, Z\right\ra\right]
=\frac{1}{\pi{\rm i}}\left\la\begin{pmatrix} A & B \\ C & F \end{pmatrix}\J_n\begin{pmatrix}
a \\ b
\end{pmatrix},Z\right\ra.
$$
Hence, for $\CA\in\sp(n,\R)$ and $\ell_0\in\R^{2n}$,
$$\ad_{\CQ_{\CA}(Z)}(\CL_{\ell_0}(Z))=-\frac12[\la Z, \CA\J_n Z\ra,\la \ell_0, Z\ra]=\frac{1}{2\pi{\rm i}}\CL_{\CA\ell_0}(Z).$$
More generally, for $m\in\N^*$, we have
$$\ad^m_{\CQ_{\CA}(Z)}(\CL_{\ell_0}(Z))=\frac{1}{(2\pi{\rm i})^m}\CL_{\CA^m\ell_0}(Z).\qed$$

Given that $Q=2\pi{\rm i} t \CQ_{\CA}(Z)$ is a pseudo-differential operator of order $2$ and $L=2\pi{\rm i} \CL_{\ell(t)}(Z)$ of order $1$, any term in the Baker-Campbell-Hausdorff series that includes more than two instances of $L$ results in an operator of negative order, and therefore is considered to be zero. We can categorize the non-zero Lie elements based on the number of $L$ terms they contain:
\begin{itemize}
\item (No $L$ term) The only term without $L$ is $Q$ since $[Q,Q]=0$.
\item (One $L$ terms) Terms with a single $L$ can either be of the form
$$
[\underbrace{Q,\cdots,Q}_{m \ {\rm times}},L]=\ad_{Q}^m L \quad {\rm or}  \quad [\underbrace{Q,\cdots,Q}_{m-1 \ {\rm times}},L,Q]=-[\underbrace{Q,\cdots,Q}_{m \ {\rm times}},L]=-\ad_{Q}^m L.$$
Given that $\CA$ is nilpotent, these terms are non-zero only if $\CA^m$ is non-zero, as per Lemma \ref{commutatorcompute}. Therefore, there are only a finite number of such terms due to the nilpotency of $\CA$.
\item (Two $L$ terms)
For terms containing two instances of $L$, considering the commutator $[L, \ad_Q^m L]$ which is of order $0$, the relevant form is $[L, \ad_Q^m L]$. Like the previous case, the nilpotency of $\CA$ ensures that there are only a finite number of such terms.
\end{itemize}
This analysis demonstrates that the Baker-Campbell-Hausdorff series, in this context, terminates after a finite number of terms, making it a finite sum due to the specific nature of the operators involved and the nilpotency of $\CA$.		

\medskip

To demonstrate the expression for $e^{2\pi{\rm i}t\CR^{(n)}_5}u$ as outlined in (\ref{solCase5}), and similarly for $e^{2\pi{\rm i}t\CR^{(n)}_4}u$ in (\ref{solCase4}), we consider the operator $\exp\{2\pi{\rm i}t\CR^{(n)}_5\}=\exp\left\{2\pi{\rm i} t \CQ_{\CB}(Z)+2\pi{\rm i}t \CL_{l}(Z)\right\}$, where $n \geq 3$ is odd, $\CB = \CB_5^{(n)}$, and $l = l_5^{(n)}$ as defined in Theorem \ref{main1}. Here, $\CB$ satisfies $\CB^n = 0$, ${\rm Rank}(\CB^{n-1}) = 2$, and $l = \iota_1 {\bf e}_{1} + \iota_2 {\bf e}_{2n}$, with $\iota_1, \iota_2 \in \R$.
Our goal is to find $\ell_0(t) \in \R^{2n}$ such that, by formally applying B-C-H formula (\ref{formu-BCH}) to
$Q=2\pi{\rm i} t \CQ_{\CB}(Z)$ and $L=2\pi{\rm i}\CL_{\ell_0(t)}(Z)$
(even though $Q$ and $L$ are unbounded operators and not matrices), we obtain:
$$\exp\left\{2\pi{\rm i} t \CQ_{\CB}(Z)\right\}\exp\{2\pi{\rm i}\CL_{\ell_0(t)}(Z)\}=\exp\{2\pi{\rm i}t\CR^{(n)}_5+2\pi{\rm i}f(t)\},$$
with $f(t)$ independent of $Z$ (hence commuting with $2\pi{\rm i}t\CR^{(n)}_5$).			
By expanding the series according to the B-C-H formula and considering the specific nature of $Q$ and $L$, we can express:
$$2\pi{\rm i}t\CR^{(n)}_5+2\pi{\rm i}f(t)=Q+\sum_{k=0}^{n-1} c_k \ad_{Q}^k(L)+\sum_{k=0}^{n-1} d_k [L,\ad_{Q}^k(L)],$$
where $c_k$ is the coefficient of $\ad_Q^k L$ in the B-C-H formula (\ref{formu-BCH}) with $c_0=1$ and $d_k$ is that of $[L,\ad_{Q}^k(L)]$.
Since $\ad_{Q}^k(L)$ and $[L,\ad_Q^k L]$ are pseudo-differential operators of order $1$ and $0$ respectively, we have
$$2\pi {\rm i}f(t)= \sum_{k=0}^{n-1} d_k [L,\ad_{Q}^k(L)],\qquad 2\pi{\rm i}t \CL_{l}(Z) =\sum_{k=0}^{n-1} c_k \ad_{Q}^k(L)=2\pi{\rm i} \CL_{\sum_{k=0}^{n-1} c_k t^k \CB^k\ell_0(t)}(Z).$$
Thus, we need to solve for $\ell_0(t)$ such that $\CE\ell_0(t) = t(\iota_1 {\bf e}_{1} + \iota_2 {\bf e}_{2n})$, where $\CE:= \sum_{k=0}^{n-1} c_k t^k\CB^k$.
By a straightforward computation, we see that $\CE$ is invertible with
$$\CE^{-1}=\sum_{k=0}^{n-1}c_k'(t\CB)^k=\left(\begin{array}{cccc:cccc}
						1 &   &  &  &  &  &  &     \\
						-c_1't & \ddots &   & &  &  &  &    \\
						\vdots &   \ddots  &  \ddots & &  &  &  &    \\
						c_{n-1}'(-t)^{n-1}& \cdots  &   -c_1't & 1 &  &  &  &   \\
						\hdashline
						&  &  &  & 1 & c_1't   & \cdots  &   c_{n-1}'t^{n-1}  \\
						&  &  &  &   &  \ddots &\ddots  & \vdots   \\
						&  &  &  & & & \ddots &  c_1't \\
						&  &  &  & & &  & 1
					\end{array}
					\right),$$
where $c_0'=1$ and $c_k'$ satisfies $\sum_{k=0}^mc_{m-k}c_{k}'=0$ for $m\in\N^*$.
Hence,
$$\ell_0(t)=t\CE^{-1}(\iota_1 {\bf e}_{1}+\iota_2 {\bf e}_{2n})
=(\iota_1 t, -\iota_1 c_1't^2, \cdots, -\iota_1c_{n-1}'(-t)^{n},
\iota_2c'_{n-1}t^n, \cdots,  \iota_2 c'_1 t^2, \iota_2 t)^*,$$
and, for $u\in L^2(\R^n)$, $\exp\{2\pi{\rm i}t\CR^{(n)}_5\} u$ can be written as
$$\exp\{2\pi{\rm i}t\CR^{(n)}_5\} u=
e^{-2\pi{\rm i}f(t)}e^{2\pi{\rm i} t \CQ_{\CB}(Z)}e^{2\pi{\rm i}\CL_{\ell_0(t)}(Z)} u=e^{-2\pi{\rm i}f(t)}\M(e^{t\CB})\rho(\ell_0(t))u.$$
Moreover, one can prove that $c_k'=\frac{(-1)^{k}}{(k+1)!}$ and $c_k$ is the $k^{\rm th}$ Bernoulli number \cite{G56}.

\section{The logarithm of matrix}

For any complex matrix $M \in \C^{n \times n}$, finding the solution $X \in \C^{n \times n}$ to the matrix equation $M = e^X$, known as a logarithm of $M$, is an intriguing problem. Notably, the logarithm of a matrix is not uniquely determined, a phenomenon that persists even within the realm of real matrices. However, a particular type of logarithm, termed the principal logarithm, is defined under certain conditions.

\begin{lemma}\label{principlelog}[Principal logarithm, Theorem 1.31 and 11.1 of \cite{Hig2008}]
For $M\in {\rm GL}(n,\C)$, if $M$ has no real negative eigenvalues, then there is a unique logarithm, denoted by $\ln(M)$ and called the principle logarithm of $M$, all of whose eigenvalues lie in the strip $\{z\in\C:-\pi< \Im z < \pi\}$. If $M$ is a real matrix, so is $\ln(M)$. Besides,
	$\displaystyle \ln(M)=(M-I_n)\int_0^1(s(M-I_n)+I_n)^{-1}ds$.
	
\end{lemma}
\indent Furthermore, if $M\in\Sp(n,\R)$, then we have
	\begin{lemma}\label{logofsymp}[Logarithm of symplectic matrices, Lemma 5.2.2 of \cite{MeyOff2017}]
		For $M\in\Sp(n,\R)$, if $M$ has no real negative eigenvalues, then $M$ has a real logarithm $X\in\sp(n,\R)$.
	\end{lemma}
	\begin{remark}\label{prilog}
In the detailed proof of Lemma \ref{logofsymp}, it is demonstrated that the eigenvalues of the logarithm $X$ are situated within the strip $\{z \in \C : -\pi < \Im(z) < \pi\}$. By virtue of the uniqueness property of the principal logarithm, this logarithm $X$ is identified as the principal logarithm of the matrix $M$.
	\end{remark}

\end{document}